\documentclass[10pt,bibliography=totoc]{article}
\usepackage{etex}
\usepackage[T1]{fontenc} 
\usepackage{latexsym}
\usepackage{empheq}
\usepackage{amsmath,amssymb,amsfonts,amsthm}
\usepackage{makecell}

\usepackage[american]{babel}
\usepackage{mathrsfs}
\usepackage{eqnarray}
\usepackage{graphicx}
\usepackage{tabularx}
\usepackage{caption}
\usepackage{subcaption}
\usepackage[table]{xcolor}
\usepackage{colortbl}
\usepackage{multicol}
\usepackage{multirow}
\usepackage{tikz}
\usepackage{bbm}
\usetikzlibrary{intersections,shapes.arrows}
\usetikzlibrary{decorations.pathreplacing,decorations.markings}
\usetikzlibrary{arrows,decorations.markings}
\usetikzlibrary{spy}
\usepackage{pgfplots}
\pgfplotsset{compat=1.14}
\usepgfplotslibrary{groupplots}
\usepgfplotslibrary{colorbrewer}
\usepgfplotslibrary{fillbetween}
\pgfplotsset{
    cycle list/Set1,
}
\usepackage{epsfig}
\usepackage{textcomp}
\usepackage{color,xcolor}
\usepackage{enumitem} 
\usepackage{marginnote} 
\usepackage{todonotes}
\presetkeys{todonotes}{size=\footnotesize, color=white, linecolor = black}{}
\usepackage{mathtools}
\usepackage{extarrows}
\usepackage{xspace}
\usepackage{microtype}
\usepackage{csquotes} 
\usepackage[linesnumbered,ruled,vlined]{algorithm2e}
\usepackage{float}
\usepackage{authblk}
\usepackage{tikz-cd}
\usepackage[text={6.25in,9in},centering]{geometry}
\usepackage{comment}  
\usepackage{empheq}
\usepackage{soul} 
\usepackage{bm}
\usepackage{cases}
\usepackage{units}

\allowdisplaybreaks

\title{Adaptive symplectic model order reduction of parametric particle-based Vlasov--Poisson equation}
\author[]{%
Jan S. Hesthaven\thanks{Chair of Computational Mathematics and Simulation Science (MCSS),
			  \'Ecole Polytechnique F\'ed\'erale de Lausanne (EPFL),
			  CH-1015 Lausanne, Switzerland.
Email: \texttt{Jan.Hesthaven@epfl.ch}}\;,
Cecilia Pagliantini\thanks{Department of mathematics, University of Pisa, Italy.
Email: \texttt{cecilia.pagliantini@unipi.it,\,c.pagliantini@tue.nl}}\;, and
Nicol\`o Ripamonti\thanks{Chair of Computational Mathematics and Simulation Science (MCSS),
			  \'Ecole Polytechnique F\'ed\'erale de Lausanne (EPFL),
			  CH-1015 Lausanne, Switzerland.
Email: \texttt{nicolo.ripamonti@epfl.ch}}
}
\date{}

\usepackage[hyphens]{url}
\usepackage[
	style=numeric,
	citestyle=numeric,sorting=nyvt,sortcites=false,
	giveninits=true,
	maxnames=10,
	url=false,eprint=false,doi=false,isbn=false,
	backref=false,
    backend=biber,
	natbib=true,
]{biblatex}
\renewbibmacro{in:}{} 
\AtEveryBibitem{\clearfield{note}} 
\AtEveryBibitem{\ifentrytype{book}{\clearfield{pages}}{}}
\usepackage[
	pdftex,
	bookmarks,
	colorlinks=true, 
	pdftitle={Dynamical MOR for VP},
	pdfauthor={},
	pdfsubject={},
	pdfkeywords={},
	urlcolor= blue,
	linkcolor= purple,
	citecolor= purple
]{hyperref}
\usepackage{cleveref}%

\usetikzlibrary{external}
\usepackage{standalone}

\newcommand{\Ybf}{\ensuremath{\mathbf{Y}}}
\newcommand{\xbf}{\ensuremath{\mathbf{x}}}
\newcommand{\ybf}{\ensuremath{\mathbf{y}}}
\newcommand{\Ecal}{\ensuremath{\mathcal{E}}}
\newcommand{\Hcal}{\ensuremath{\mathcal{H}}}

\newcommand{\Mcal}{\ensuremath{\mathcal{M}}}
\newcommand{\Ncal}{\ensuremath{\mathcal{N}}}
\newcommand{\Pcal}{\ensuremath{\mathcal{P}}}
\newcommand{\Rcal}{\ensuremath{\mathcal{R}}}
\newcommand{\Tcal}{\ensuremath{\mathcal{T}}}
\newcommand{\Ucal}{\ensuremath{\mathcal{U}}}
\newcommand{\Xcal}{\ensuremath{\mathcal{X}}}
\newcommand{\Ycal}{\ensuremath{\mathcal{Y}}}
\newcommand{\Zcal}{\ensuremath{\mathcal{Z}}}
\newcommand{\Tcalt}{\ensuremath{{\Tcal_{\tau}}}}
\newcommand{\norm}[1]{\|{#1}\|}

\newcommand{\Vd}[1]{\ensuremath{\mathcal{V}_{#1}}}

\newcommand{\J}[1]{\ensuremath{J_{#1}}}
\newcommand{\HamC}{\ensuremath{H}} 
\newcommand{\Ham}{\ensuremath{\Hcal}} 
\newcommand{\HamU}{\ensuremath{\Hcal_{U}}}
\newcommand{\HamUi}{\ensuremath{\Hcal_{U,i}}}

\newcommand{\Hamgen}[1]{\ensuremath{\Hcal_{U_{#1}}}}
\newcommand{\EE}{\ensuremath{\Ecal}}

\newcommand{\R}[2]{\ensuremath{\mathbb{R}^{{#1}\times{#2}}}}
\renewcommand{\r}[1]{\ensuremath{\mathbb{R}^{#1}}}

\newcommand{\V}{\ensuremath{\xi}}

\newcommand{\Nf}{\ensuremath{{2N}}} 
\newcommand{\Nr}{\ensuremath{{2n}}} 
\newcommand{\Nfh}{\ensuremath{N}} 
\newcommand{\Nrh}{\ensuremath{n}} 
\newcommand{\Np}{\ensuremath{p}} 
\newcommand{\Nps}{\ensuremath{p^*}} 
\newcommand{\ns}{\ensuremath{{n_s}}} 
\newcommand{\is}{\ensuremath{{i_s}}} 
\newcommand{\Nx}{\ensuremath{N_x}} 
\newcommand{\T}{\ensuremath{T}} 
\newcommand{\Ndmdt}{\ensuremath{r_{\tau}}} 
\newcommand{\Ndeim}{\ensuremath{d}} 
\newcommand{\nDEIMup}{\ensuremath{n_{\text{DEIM}}}} 
\newcommand{\DEIMfreq}{\ensuremath{k_{\text{DEIM}}}} 

\newcommand{\dt}{\ensuremath{\Delta t\,}} 

\newcommand{\dX}{\ensuremath{\dot{X}}}
\newcommand{\dV}{\ensuremath{\dot{V}}}

\DeclareMathOperator*{\argmin}{argmin}
\DeclareMathOperator*{\tol}{tol}
\DeclareMathOperator*{\St}{St}
\DeclareMathOperator*{\Simp}{Sp}

\DeclareMathOperator*{\dmd}{DMD}
\DeclareMathOperator*{\deim}{DEIM}
\DeclareMathOperator*{\dd}{DD}

\newcommand{\Idm}{\ensuremath{I}}

\newcommand{\prm}{\ensuremath{\eta}} 
\newcommand{\prmi}{\ensuremath{\eta_i}} 
\newcommand{\prmhsub}{\ensuremath{\Sprmh^{*}}} 
\newcommand{\Sprm}{\ensuremath{\Gamma}} 
\newcommand{\Sprmh}{\ensuremath{\Sprm_h}} 
\newcommand{\Vprm}{\ensuremath{V}_{\prm}}
\newcommand{\Sp}{\ensuremath{S_{\Np}}} 
\newcommand{\Spsub}{\ensuremath{S_{\Nps}}} 
\newcommand{\Mred}{\ensuremath{\mathcal{M}_{\Nr}^{\Np}}} 

\newcommand{\mmp}{\ensuremath{m_p}}
\newcommand{\Mq}{\ensuremath{M_q}}

\newcommand{\Ghp}{\ensuremath{G^{\,\Np}_{\Ham}(U,Z)}}
\newcommand{\Ghpi}[1]{\ensuremath{G^{\,\Np}_{\Ham}(U,Z_{#1})}}
\newcommand{\ghp}{\ensuremath{g^{\Np}_{\Ham}(U,Z)}}
\newcommand{\ghpi}[1]{\ensuremath{g^{\Np}_{\Ham}(U,Z_{#1})}}
\newcommand{\Ghps}{\ensuremath{G^{\,\Nps}_{\Ham}(U,Z)}}

\newcommand{\Ghddi}[1]{\ensuremath{G^{\dd}_{\Ham}(U,Z_{#1},t)}}
\newcommand{\ghdd}{\ensuremath{g^{\dd}_{\Ham}(U,Z,t)}}
\newcommand{\ghddi}[1]{\ensuremath{g^{\dd}_{\Ham}(U,Z_{#1},t)}}
\newcommand{\ghddvar}[1]{\ensuremath{g^{\dd}_{\Ham}\bigg(#1\bigg)}}
\newcommand{\ghddv}{\ensuremath{g^{\dd}_{\Ham}}}

\numberwithin{equation}{section}
\theoremstyle{plain}
  \newtheorem{theorem}{\sffamily Theorem}[section]
  \Crefname{proposition}{Proposition}{Propositions}
  \Crefname{lemma}{Lemma}{Lemmas}

  \theoremstyle{definition}
  \newtheorem{remark}[theorem]{\sffamily Remark}\Crefname{remark}{Remark}{Remarks}

  \newtheorem*{acknowledgment*}{Acknowledgment}

\makeatletter
\def\tikz@auto@anchor{%
    \pgfmathtruncatemacro\angle{atan2(\pgf@x,\pgf@y)-90}
    \edef\tikz@anchor{\angle}%
}
\makeatother 

\begin{document}

\maketitle

\begin{abstract}
High-resolution simulations of particle-based kinetic plasma models typically require a high number of particles and thus often become computationally intractable.
This is exacerbated in multi-query simulations, where the problem depends on a set of parameters.
In this work, we derive reduced order models for the semi-discrete Hamiltonian system resulting from a geometric particle-in-cell approximation of the parametric Vlasov--Poisson equations.
Since the problem's non-dissipative and highly nonlinear nature makes it reducible only locally in time, we adopt a nonlinear reduced basis approach where the reduced phase space evolves in time. This strategy allows a significant reduction in the number of simulated particles, but the evaluation of the nonlinear operators associated with the Vlasov--Poisson coupling remains computationally expensive.
We propose a novel reduction of the nonlinear terms that combines adaptive parameter sampling and hyper-reduction techniques to address this.
The proposed approach allows decoupling the operations having a cost dependent on the number of particles from those that depend on the instances of the required parameters.
In particular, in each time step, the electric potential is approximated via dynamic mode decomposition (DMD) and the particle-to-grid map via a discrete empirical interpolation method (DEIM).
These approximations are constructed from data obtained
from a past temporal window at a few selected values of the parameters to guarantee a computationally efficient adaptation.
The resulting DMD-DEIM reduced dynamical system retains the Hamiltonian structure of the full model, provides good approximations of the solution, and can be solved at a reduced computational cost.    
\end{abstract}


\section{Introduction}
The kinetic modeling of collisionless magnetized plasmas is based on the Vlasov--Maxwell equations, which describe the evolution of the distribution function of a collection of charged particles under the action of self-consistent electromagnetic fields.
Because of the high dimensionality of the phase space, 
the large separation of scales,
the inherent nonlinearity, and the infinitely many conserved quantities, the numerical treatment of the Vlasov–Maxwell equations, and of its electrostatic limit of Vlasov--Poisson equations, is a challenging task.

Arguably, the most widely used family of numerical methods for the solution of kinetic plasma models are Particle-In-Cell (PIC) methods \cite{LB91}. The idea of PIC schemes is to sample the distribution function in velocity space using a finite number of macro-particles that are evolved along their characteristics. The electromagnetic fields are discretized on a grid in the computational domain, and the macro-particles move through the grid according to the Lorentz force. 
%
To preserve key physical properties of the problem, such as conservation of total energy, PIC schemes have evolved into variational algorithms based on least action principles \cite{Low58,Lewis70,SQT12} or on the Hamiltonian formulation of the Vlasov--Maxwell and Vlasov--Poisson equations \cite{MW82,ES13}.
In parallel, several numerical methods \cite{SQT12,ES13} for kinetic plasma models have leveraged discrete differential forms and de Rham complexes for the geometric approximation of the electromagnetic fields through Maxwell's equations.
Combining these ideas of structure-preserving approximations has led to finite element PIC methods able to exactly satisfy physical constraints, like the Gauss laws, and to guarantee the preservation of the Hamiltonian structure of the problem.
Examples include the canonical \cite{Qin16} and
non-canonical \cite{Xiao15} symplectic particle-in-cell algorithms,
the Hamiltonian particle-in-cell method of \cite{He16}, and its generalization, the Geometric Electromagnetic PIC (GEMPIC) method \cite{KKMS17}.

The multiscale nature of plasmas implies that PIC codes require a significant amount of computational resources to resolve the shortest length scale and the fastest plasma frequency, and, thus, to yield stable and accurate numerical approximations.
Moreover, the slow convergence rate of particle-based methods necessitates the use of many particles to achieve sufficient accuracy, to capture, for example, time-dependent physical phenomena such as plasma instabilities. As a result, PIC methods can be prohibitively expensive in terms of computational cost.
This computational burden can become intractable
in the parametric case, when simulations for many input parameters are of interest.
This problem has been tackled from an algorithmic standpoint via the improvement of algorithms' structure and the use of suitable computational hardware.
In this work, we propose to address this computational issue through model order reduction.
Starting from the high-resolution geometric PIC approximation of the Vlasov--Poisson problem, the idea is to derive a low-dimensional surrogate model that can be
solved at a reduced computational cost and still provides accurate approximate solutions.

In recent years, several works have been devoted to the development of reduced representations of kinetic plasma dynamics based on model order reduction techniques, sparse grids and low-rank approximation \cite{Kor15,KS16,EL17,Eink1,Eink2,Eink21}.
A dynamical low-rank approximation of the infinite-dimensional Vlasov--Poisson problem
has been proposed in \cite{Eink1,Eink2,Eink21}
with the intent to optimize the number of degrees of freedom needed for a sufficiently accurate and stable approximation of the solution.
The continuous distribution function of the Vlasov--Poisson problem is expanded into a finite sum of low-rank factors, for which a new dynamical system is derived. Discretization in space and velocity of the resulting problem yields a low-rank approximation of the original dynamics.
Another class of methods has considered the model reduction of the ordinary differential equation ensuing from the semi-discretization of the Vlasov--Poisson problem. A semi-Lagrangian method to solve the problem on a sparse grid is proposed in \cite{KS16}. In \cite{Kor15}, time splitting is combined with a spline-based semi-Lagrangian
scheme where each step consists of linear combinations of low-rank approximations. In \cite{EL17}, the solution of the full-Eulerian time-dependent Vlasov--Poisson system is approximated using a tensor decomposition whose rank is adapted at each time step. To address the problem of the fluid closure for the collisionless linear Vlasov system,
an interpolatory order reduction is proposed in \cite{gillot2021model}.
In the context of a particle-based discretization of kinetic plasma models,
a dynamic mode decomposition (DMD) strategy has been proposed in \cite{nayak2020reduced}
to reconstruct the electric field within an Electromagnetic particle-in-cell (EMPIC) algorithm.
Although the proposed approach can effectively capture and extrapolate the electric field behavior around equilibria,
the computational burden associated with the high number of particles is not overcome.
Model order reduction --- in the number of particles --- of \textit{parametric} plasma models is, to the best of our knowledge, an open problem. Moreover, the conservation of the Hamiltonian structure has not been explored in reduced representations of plasma dynamics. This work aims at addressing these issues.

We focus on the parametric 1D-1V Vlasov--Poisson problem and consider its discretization in space and velocity via a geometric PIC method. The resulting dynamical system has a Hamiltonian form and, hence, corresponds to a symplectic flow.
Model order reduction of the Vlasov--Poisson problem poses some major challenges, and standard reduced basis techniques are prone to fail in terms of numerical stability, computational efficiency, and accuracy of the simulations.
The application \cite{tyranowski2019symplectic} of the symplectic reduced basis methods proposed in \cite{peng2016symplectic,afkham2017structure} to the Vlasov equation, with a fixed external electrostatic field, shows the importance of retaining the symplectic structure of Vlasov's equation in the reduced model. However, the multi-scale nature of the problem makes it difficult for a reduced order model to characterize, with sufficient accuracy, the plasma behavior using a small number of degrees of freedom. This implies that accurate reduced representations of the solution may require large approximation spaces that eliminates the benefits of model order reduction. To overcome this limitation, we focus on methods that adapt the approximate reduced space in time while preserving the geometric properties of the original problem.
This approach allows to accurately describe the plasma evolution with a considerably reduced number of particles without compromising the quality of the simulation. 
The bulk of the computational effort to solve the reduced dynamics is due to the nonlinearity of the particles-to-grid mapping, and thus the Hamiltonian, whose evaluation needs to be performed in the original high-dimensional space.
To alleviate these computational inefficiencies, we propose a strategy that 
approximates the reduced Hamiltonian gradient via a combination of hyper-reduction techniques and parameter sampling procedures. A reduction in the computational runtimes of the algorithm is achieved by decoupling the operations that depend on the number of particles from those that depend on the number of parameters while retaining an accurate representation of the plasma dynamics.
The resulting dynamical system preserves the symplectic structure of the problem, ensures the stability of the approximation, and exploits the local-in-time low-rank nature of the solution.

The remainder of this work is organized as follows.
In \Cref{sec:physical_model}, the parametric Vlasov--Poisson problem is introduced both in its classical Eulerian formulation and in the Hamiltonian formulation. Moreover, the semi-discrete approximation of the problem via a particle method coupled with a finite element discretization of the Poisson problem is described.
In \Cref{sec:MOR}, the model order reduction of the parametric dynamical system originating from the semi-discretization of the Vlasov--Poisson equation is considered. First, we describe a global-in-time symplectic reduced basis method, and then we consider a nonlinear approach, based on dynamical low-rank approximation, still preserving the symplectic structure of the flow field.
After discussing the computational complexity of the dynamical reduced basis algorithm, \Cref{sec:hyperred} is devoted to the DMD-DEIM structure-preserving approximation of
the nonlinear Hamiltonian gradient and the particles-to-grid mapping. The numerical temporal discretization of the resulting reduced model is discussed in \Cref{subsec:temporal_integration}.
Numerical experiments in \Cref{sec:numExp} on benchmark tests show that the proposed method can accurately reproduce the dynamics of particle-based kinetic plasma models
with significant speedups compared to solving the original system.
Concluding remarks are presented in \Cref{sec:conclusions}.

\section{The physical model}\label{sec:physical_model}
We consider the parametric Vlasov--Poisson problem with parameters that describe physical properties of the system. 
In particular, we focus on the study of the effect of parametrized initial distributions on the plasma dynamics.
Let us assume that the parameters range in a compact set $\Sprm\subset\r{q}$, with $q\geq 1$.
The plasma, at any time $t\in \Tcal\subset\mathbb{R}$, is described in terms of the distribution function $f(t,x,v;\prm)$ in the Cartesian phase space domain
$(x,v)\in\Omega:=\Omega_x\times\Omega_v\subset\mathbb{R}^2$.
We focus on single species plasmas but the proposed method can be trivially extended to the multi-species case.
Assume that $\Omega_x:=\mathbb{T}^d=\mathbb{R}^d/(2\pi\mathbb{Z}^d)$ is the $d$-dimensional torus
and $\Omega_v:=\mathbb{R}^d$.
For $\prm\in\Sprm$ fixed, we introduce the space
\begin{equation*}
\begin{aligned}
	& \Vprm:=\{f(t,\cdot,\cdot\,;\prm)\in L^2(\Omega):\;f(t,x,v;\prm)>0\,\mbox{ for all }\, (x,v)\in\Omega,\; f(t,\cdot,v;\prm)\sim e^{-v^2}\,\mbox{ as }\,|v|\to\infty\}.
\end{aligned}
\end{equation*}
The 1D-1V Vlasov--Poisson problem ($d=1$) reads:
For each $\prm\in\Sprm$ and $f_0(\prm)\in{\Vprm}\big|_{t=0}$, 
find $f(\cdot,\cdot,\cdot\,;\prm)\in C^1(\Tcal;L^2(\Omega))\cap C^0(\Tcal;\Vprm)$, and the electric field $E(\cdot,\cdot\,;\prm)\in C^0(\Tcal;L^2(\Omega_x))$ such that
\begin{equation*}
    \left\{
    \begin{aligned}
	& \partial_t f(t,x,v;\prm)+v\,\partial_x f(t,x,v;\prm) + \dfrac{q}{m} E(t,x;\prm)\,\partial_v f(t,x,v;\prm)=0, &\qquad\mbox{in}\;\Omega\times\Tcal,\\
	& \partial_x E(t,x;\prm)= q\int_{\Omega_v}f(t,x,v;\prm)\, dv,                       &\qquad\mbox{in}\;\Omega_x\times\Tcal,\\
	& f(0,x,v;\prm)=f_0(\prm),						  &\qquad\mbox{in}\;\Omega.
\end{aligned}\right.
\end{equation*}
Here $q$ is the charge and $m$ is the particle mass.
The boundary conditions for $f$ are periodic in space and
prescribed via the space $\Vprm$ in velocity.

Since the electric field can be written as the spatial derivative of the electric potential $\phi$, namely $E(t,x;\prm)=-\partial_x\phi(t,x;\prm)$,
the Vlasov--Poisson problem can be recast, for each $\prm\in\Sprm$, as
\begin{equation}\label{eq:VPpot}
    \begin{aligned}
	& \partial_t f(t,x,v;\prm)+v\,\partial_x f(t,x,v;\prm) - \dfrac{q}{m}
    \partial_x\phi(t,x;\prm)\,\partial_v f(t,x,v;\prm)=0, &\qquad\mbox{in}\;\Omega\times\Tcal,\\
	& -\partial_{xx} \phi(t,x;\prm)=\rho(t,x;\prm):=q\int_{\Omega_v}f(t,x,v;\prm)\, dv,                     &\qquad\mbox{in}\;\Omega_x\times\Tcal,
\end{aligned}
\end{equation}
where $\rho$ is the electric charge density. In the presence of a background charge, as in the numerical experiments, the right hand side of the Poisson problem reads $1-\int_{\Omega_v}f(t,x,v;\prm)\, dv$ and the formulation of the method presented in the following sections can be extended straightforwardly.

The Lagrangian and Hamiltonian formulation
of the Vlasov--Poisson and Vlasov--Maxwell equations
reveals a set of mathematical and geometric features that encode the physical properties of these systems.
The Vlasov--Poisson problem admits a Hamiltonian formulation, \emph{cf.} \cite[Sections 1 and 2]{CCFM17}, with a Lie--Poisson bracket, and Hamiltonian
given by the sum of the kinetic and electric energy as
\begin{equation}\label{eq:VPHam}
    \HamC(f,\prm) = \dfrac{m}{2} \int_{\Omega} v^2 f(t,x,v;\prm) \,dx\,dv + \dfrac12\int_{\Omega_x} |E(t,x;\prm)|^2\,dx.
\end{equation}
Eulerian-based discretizations of kinetic plasma models in Hamiltonian form, with general noncanonical Poisson brackets, do not appear to
inherit the phase space structure of the continuous problem, as has been observed in
\cite{Mo81,ScoWein94}. On the contrary, particle-in-cell methods have led to the geometric approximation of these models when coupled to the discretization of the electromagnetic fields via discrete differential forms \cite{Xiao15,Qin16,He16,KKMS17}.
For the structure-preserving approximation and reduction of the Vlasov--Poisson system \eqref{eq:VPpot} we rely on the Hamiltonian structure of its semi-discrete formulation obtained via particle-based methods as derived in the following.


\subsection{Geometric particle-based discretization}\label{sec:GEMPIC}

We consider a particle method for the approximation of the Vlasov equation, coupled with a $H^1$-conforming discretization of the Poisson problem for the electric potential.
In detail, the distribution function $f^s$ is approximated by the superposition of $\Nfh\in\mathbb{N}$ computational macro-particles as
\begin{equation}\label{eq:fapprox}
    f(t,x,v;\prm)\approx f_h(t,x,v;\prm) = \sum_{\ell=1}^\Nfh \omega_\ell\, S(x-X_\ell(t,\prm))\, \delta(v-V_\ell(t,\prm)),
\end{equation}
where $\omega_\ell\in\mathbb{R}$ is the weight of the $\ell$-th particle,
$\delta$ is the Dirac delta,
$S$ is a compactly supported shape function,
and, for each $\prm\in\Sprm$ and $t\in\Tcal$, $X(t,\prm)\in\mathbb{R}^{\Nfh}$ and $V(t,\prm)\in\mathbb{R}^{\Nfh}$ denote the vectors of the position and velocity of the macro-particles, respectively.
In this work, we assume homogeneous macro-particles' weights so that $m_{\ell}=m$ and $\omega_{\ell} =\omega$, for all $\ell=1,\ldots,\Nfh$.
The proposed method can be extended to the general case \emph{mutatis mutandis}.

The idea of particle methods is to derive the time evolution of the approximate distribution function $f_h$ by advancing the macro-particles along the characteristics of the Vlasov equation, i.e. the particles' positions and velocities satisfy the following set of ordinary differential equations
\begin{equation}\label{eq:fom}
    \begin{aligned}
        & \dX(t,\prm) = V(t,\prm), &\qquad \mbox{in }\Tcal, \\
        & \dV(t,\prm) = \dfrac{q}{m} E(t,X(t,\prm);\prm), &\qquad \mbox{in }\Tcal,
    \end{aligned}
\end{equation}
 under suitable initial conditions.
 Here and throughout, the dot denotes the derivative
 with respect to time.
The macro-particles move through a computational grid under the influence of electromagnetic
fields. The latter are self-consistently calculated from the positions of the particles on the grid via the Poisson equation \eqref{eq:VPpot}.
On a partition of the spatial domain $\Omega_x$, we consider a finite element discretization of the Poisson equation in the space
$\Pcal_k\Lambda^0(\Omega_x)\subset H^1(\Omega_x)$
of continuous piecewise polynomial functions of degree at most $k\geq 1$.
The semi-discrete variational problem reads: For each $\prm\in\Sprm$, find $\phi_h(\cdot\,;\prm)\in C^1(\Tcal;\Pcal_k\Lambda^0(\Omega_x))$ such that
\begin{equation}\label{eq:VF}
    a_h(\phi_h(t,\cdot\,;\prm),\psi) = g(\psi),\qquad \forall\,\psi\in\Pcal_k\Lambda^0(\Omega_x),
\end{equation}
where the bilinear form $a_h$ corresponds to the Laplace operator, while the
linear function $g$ is associated with the density $\rho$\,; thereby
\begin{equation*}
    a_h(\varphi,\psi) := \int_{\Omega_x}d_x\varphi(x)\, d_x\psi(x)\,dx,
    \qquad
    g(\psi) := q\int_{\Omega}f_h(t,x,v;\prm)\psi(x)\, dv\,dx,
    \quad\forall\psi,\varphi\in\Pcal_k\Lambda^0(\Omega_x).
\end{equation*}
Let $\{\lambda_i^0(x)\}_{i=1}^{\Nx}$ be a basis of the space $\Pcal_k\Lambda^0(\Omega_x)$, so that
the semi-discrete potential can be written as
$\phi_h(t,\cdot\,;\prm)=\sum_{i=1}^{\Nx}\Phi_i(t,\prm)\lambda_i^0(x)$, for any $t\in\Tcal$ and $\prm\in\Sprm$.
We introduce the matrices
$L\in\R{\Nx}{\Nx}$ and $\Lambda^0(X)\in\R{\Nfh}{\Nx}$, defined as
\begin{equation}
    L_{i,j}:=(d_x \lambda_j^0,d_x \lambda_i^0)_{L^2(\Omega_x)},
    \quad
    (\Lambda^0(X))_{\ell,i}:=\lambda^0_i(X_\ell),
    \quad i,j=1,\ldots,\Nx,\,\ell=1,\ldots,\Nfh.
\end{equation}
Then the time-dependent algebraic system ensuing from \eqref{eq:VF} reads
\begin{equation}\label{eq:PoissonSystem}
    L\Phi(t,\prm) = \Lambda^0(X(t,\prm))^\top \Mq =: \rho_h(X(t,\prm);\prm),
\end{equation}
where
$\Mq\in\r{\Nfh}$ is the vector of entries $(\Mq)_{\ell} = q\omega_\ell$, for $\ell=1,\ldots,\Nfh$.
The proposed discretization of the electromagnetic field allows to recast the characteristic equations \eqref{eq:fom} as a Hamiltonian system.
The phase space of Hamiltonian systems is characterized by a symplectic geometric structure. We denote with $\Vd{\Nf}\subset\r{\Nf}$ the phase space of \eqref{eq:fom} and we assume it is a $\Nf$-dimensional symplectic vector space.

Let $W(t,\prm)=[X(t,\prm);V(t,\prm)]\in\Vd{\Nf}$ denote the vector of all particle positions and velocities at a given time $t\in\Tcal$ and parameter value $\prm\in\Sprm$, obtained by concatenating the vectors $X(t,\prm)$ and $V(t,\prm)$. The latter are also known as generalized position and momentum in the symplectic formalism.
We denote with $W\in\Vd{\Nf}\subset \r{\Nf}$ the state vector collecting the positions and velocities of all particles of all species.
To simplify the notation, we define the Hamiltonian $\Ham$ resulting from the proposed discretization of \eqref{eq:VPHam} up to the constant $\mmp:=m\omega$.
First, we introduce the (nonlinear) electric energy $\EE:\mathbb{R}^{\Nf}\times\Sprm\rightarrow\mathbb{R}$ defined as
\begin{equation}\label{eq:HamEEfull}
\begin{aligned}
    \EE(X(t,\prm);\prm) := 
    \dfrac{\mmp^{-1}}{2}\int_{\Omega_x} |\partial_x \phi_h(t,x;\prm)|^2\,dx
    &= \dfrac{\mmp^{-1}}{2}\Phi(t,\prm)^\top L \Phi(t,\prm)\\
    &= \dfrac{\mmp^{-1}}{2} \Mq^\top\Lambda^0(X(t,\prm)) L^{-1}\Lambda^0(X(t,\prm))^\top\Mq.
\end{aligned}
\end{equation}
The Hamiltonian is the sum of the kinetic and electric energy and it is given, for any $W(t,\prm)\in\Vd{\Nf}$, by
\begin{equation}\label{eq:high_fidelity_Hamiltonian}
    \Ham(W(t,\prm))
    =\dfrac12 V(t,\prm)^\top V(t,\prm) + \EE(X(t,\prm);\prm).
\end{equation}
Differentiating the discrete Hamiltonian in \eqref{eq:high_fidelity_Hamiltonian} with respect to the vector $X$ of particles' positions and the vector $V$ of particles' velocities, results in the semi-discrete system in Hamiltonian form
\begin{align}
    \begin{pmatrix}
    \dX(t,\prm)\\
    \dV(t,\prm)
    \end{pmatrix}
    &=  J_{\Nf}
    \begin{pmatrix}
    \mmp^{-1}\mathrm{diag}(\Mq)\nabla\Lambda^0(X(t,\prm))L^{-1}\Lambda^0(X(t,\prm))^\top \Mq\\
    V(t,\prm)
    \end{pmatrix}.
    \label{eq:VPHamEqn}
\end{align}
Here $\nabla\Lambda^0(X)\in\R{\Nfh}{\Nx}$
is defined as
$(\nabla\Lambda^0(X))_{\ell,i}:=(d_x\lambda^0_i)(X_\ell)$, for $i=1,\ldots,\Nx$ and $\ell=1,\ldots,\Nfh$, and
$J_{\Nf}\in\R{\Nf}{\Nf}$, called Poisson tensor, is the block matrix 
\begin{equation}
    J_{\Nf} =
    \begin{pmatrix}
    0_{\Nfh} & I_{\Nfh} \\ 
    -I_{\Nfh} & 0_{\Nfh}
    \end{pmatrix},
\end{equation}
with $0_{\Nfh}$ and $I_{\Nfh}$ being the null and identity matrices of size $\Nfh$, respectively.
A generalization of this discretization to the case of Vlasov--Maxwell's equations leads to the GEMPIC method
introduced in \cite{KKMS17}.
Note that the proposed semi-discretization preserves the Hamiltonian, which corresponds to the discrete energy of the system, but not the momentum.
In PIC it does not appear possible to simultaneously conserve momentum and energy  \cite{LB91}.

\section{Model order reduction of the Vlasov--Poisson problem}\label{sec:MOR}

The parametric Hamiltonian system \eqref{eq:VPHamEqn} reads: For each $\prm\in\Sprm$ and for $W_0(\prm)=[X(0,\prm);V(0,\prm)]\in\Vd{\Nf}$, find $W(\cdot,\prm)\in C^{1}(\mathcal{T},\Vd{\Nf})$ such that
\begin{equation}\label{eq:VPHamEqnParam}
\left\{
\begin{aligned}
& \dot{W}(t,\prm) =  J_{\Nf} \nabla_{W} \Ham(W(t,\prm)), \qquad \text{in }\Tcal,\\[0.5ex]
& W(0,\prm) = W_0(\prm),
\end{aligned}\right.
\end{equation}
where the initial condition $W_0(\prm)\in\Vd{\Nf}$ is prescribed by the initial distribution $f_0(\prm)$.

\subsection{Projection-based symplectic reduced basis methods}\label{sec:POSD}
The parametric differential problem \eqref{eq:VPHamEqnParam} is said to be reducible if the solution set
$\mathcal{M}:=\{ W(t,\prm)\in\Vd{\Nf}\subset\r{\Nf}:\, t\in\mathcal{T},\; \prm\in\Sprm\}$, which collects
the solution of \eqref{eq:VPHamEqnParam} under variation of time and parameter, can be well approximated by a $\Nr$-dimensional subspace 
$\mathcal{M}_{\Nr}$, with $\Nrh\ll \Nfh$.
The reducibility of a problem is measured via the so-called Kolmogorov $n$-width \cite{Pinkus85}:
roughly speaking, the faster its decay as $n$ increases, the better $\Mcal$ can be approximated by a low dimensional space.

Within model order reduction techniques,
projection-based reduced basis methods construct
a modal approximation of the solution $W(t,\prm)$, for each $(t,\prm)\in\Tcal\times\Sprm$, of the form
$UZ(t,\prm)$,
$U\in \R{\Nf}{\Nr}$ is a basis spanning the reduced approximation space $\mathcal{M}_{\Nr}$, and
{$Z(t,\prm)\in\Mcal_{\Nr}\subset\mathbb{R}^{\Nr}$} are the expansion coefficients of the approximate reduced solution in the reduced basis.
The reduced solution is then obtained by solving a lower dimensional dynamical systems obtained by the Galerkin projection of the full model \eqref{eq:VPHamEqnParam} onto the reduced space.
If the basis $U$ is orthogonal, the reduced problem reads: For each $\prm\in\Sprm$ and $W_0(\prm)$, find $Z(t,\prm)\in\Mcal_{\Nr}$ such that
\begin{equation}\label{eq:generic_ODEs_RB}
    \begin{cases}
        \dot{Z}(t,\prm) = U^{\top}J_{\Nf}\nabla_{W} \Ham(UZ(t,\prm)),\qquad \text{in }\Tcal,\\
        Z(0;\prm) = U^{\top}W_0(\prm).
    \end{cases}
\end{equation}
The reduced basis $U$ is constructed from a collection of full order solutions at sampled values of time and parameters, via SVD-type methods, like POD, or greedy algorithms.

Note that \eqref{eq:generic_ODEs_RB} does not possess the Hamiltonian structure of the full model \eqref{eq:VPHamEqnParam}. In the setting of conservative systems, it is vital not only to provide accurate reduced representations of the solution but also to preserve the geometric structure underlying the physical properties and conservation laws of the dynamics. Moreover, violating intrinsic properties of the dynamics raises concerns about the reduced model's validity and reliability as a surrogate model.
To preserve the Hamiltonian structure of the problem after reduction, symplectic reduced basis techniques have been proposed in \cite{peng2016symplectic,afkham2017structure,HP18}. 
Despite symplectic reduced basis methods providing stable and structure-preserving surrogate models, these might not be considerably more efficient to solve than the original model. Indeed, the non-dissipative nature of Hamiltonian problems is associated with poor reducibility and a slowly decaying Kolmogorov $n$-width.
As a consequence, accurate reduced representations of the solution of \eqref{eq:VPHamEqnParam} may require large approximation spaces that
jeopardize the benefits of model order reduction.
To overcome this limitation, we focus on adaptive methods where the reduced space changes in time while preserving the geometric properties of the original problem.

{We remark that, in this work, we aim at tackling the challenges posed by problems that are not globally reducible \textit{in time}. We implicitly assume that, for a fixed time, the solution under variation of parameter can be accurately represented by a low-rank approximation. The more general case of local reducibility with respect to both time and parameter is out of the scope of the present work and a short discussion on possible research in this direction is presented in the conclusions.}


\subsection{Dynamical structure-preserving model order reduction}

For the model order reduction of the semi-discrete Vlasov--Poisson problem \eqref{eq:VPHamEqnParam}, we adopt the structure-preserving dynamical reduced basis method for Hamiltonian systems introduced in \cite{P19}.
This approach can be related to both classical reduced basis methods and dynamical low-rank approximation \cite{KL07}.
Differently from a traditional reduced basis approach, the method does not have a separation into an offline phase and an online phase and it does not rely on snapshots.
On the other hand, when the reduced basis is fixed, i.e. the approach becomes global,
we recover a symplectic reduced basis method as in
\cite{peng2016symplectic,afkham2017structure}.
As described in this section, the dynamical reduced models derived in \cite{P19} are obtained via a structure-preserving projection of the original evolution problem onto the tangent space of the reduced manifold, rather than on the reduced manifold itself. This is typical of dynamical low-rank approximation techniques and dynamically orthogonal (DO) field equations \cite{SL09,musharbash2017symplectic}.
We refer to \cite[Introduction]{P19} for a detailed discussion on the relationship among these approaches.

Assume we seek to solve the semi-discrete Vlasov--Poisson problem \eqref{eq:VPHamEqnParam} for $\Np$ parameters $\Gamma_h:=\{\prm_j\}_{j=1}^{\Np}\subset\Gamma$.
{The idea of the dynamical approach \cite{P19} is to find a low-rank approximation of the solution of problem \eqref{eq:VPHamEqnParam} at each time instant. For this reason it is of interest to look at the properties of the solution, for each fixed time, under variation of the $p$ parameters of interest.}
We {therefore incorporate the parameters of interest in the problem formulation by recasting} \eqref{eq:VPHamEqnParam} as an evolution equation in the matrix-valued unknown $\Rcal(t):=[ W(t,\prm_1) | \dots | W(t,\prm_{\Np})]\in\mathcal{V}_{\Nf}^{\Np}\subset\R{\Nf}{\Np}$, with
 $\mathcal{V}_{\Nf}^{\Np}:=\mathcal{V}_{\Nf}\times\dots\times\mathcal{V}_{\Nf}$.
 Given $\Rcal_{0}=[W_0(\prm_1)|\dots|W_0(\prm_{\Np})]$, we look for $\Rcal\in C^{1}(\mathcal{T},\mathcal{V}_{2N}^{\Np})$,  such that
 \begin{equation}\label{eq:full_param_HamSystem}
     \begin{cases}
         \dot{\Rcal}(t) = 
         J_{2N}\nabla_{\Rcal} \Ham^{\Np}(\Rcal(t)), \qquad \text{in }\Tcal,\\
         \Rcal(0)=\Rcal_{0}.
     \end{cases}
 \end{equation}
The Hamiltonian of \eqref{eq:full_param_HamSystem} is a vector-valued quantity $\Ham^{\Np}:\mathcal{V}_{2N}^{\Np}\rightarrow \mathbb{R}^{\Np}$ that collects the values of the Hamiltonian of \eqref{eq:VPHamEqnParam} for each parameter in $\Sprmh$, namely
$(\Ham^{\Np}(\Rcal(\cdot)))_{j}=\Ham(W(\cdot,\prm_j))$, for $j=1,\ldots,\Np$.
Moreover, the gradient $\nabla_{\Rcal}\Ham^{\Np}(\Rcal)\in\mathcal{V}_{2N}^{\Np}$ is defined as 
\begin{equation*}
(\nabla_{\Rcal}\Ham^{\Np}(\Rcal))_{\ell,j}=
\frac{\partial \Ham(W(\cdot,\prm_j))}{\partial W_{\ell}(\cdot,\prm_j)},
\quad \ell=1,\dots,\Nf,\; j=1,\dots,\Np.
\end{equation*}
The idea of dynamical low-rank approximations is to expand the full order solution in a truncated modal decomposition where both the basis and the expansion coefficients are time-dependent.
For all  $t\in\Tcal$, we approximate $\Rcal(t)$, in \eqref{eq:full_param_HamSystem}, as
$R(t) = U(t)Z(t)$, where {$U(t)\in\R{\Nf}{\Nr}$ is the time-dependent basis and $Z(t)\in\R{\Nr}{\Np}$ are the associated expansion coefficients}, with $Z(t):=[ Z_1(t)| \dots | Z_{\Np}(t)]$, and $Z_i(t):=Z(t,\prm_i)$.
The reduced space is then defined as
\begin{equation}\label{eq:dynManifold}
    \Mred:=\{ R\in\mathbb{R}^{\Nf\times \Np}: R=UZ \text{ with } U\in\mathcal{U}(\Nr,\mathbb{R}^{\Nf}), Z\in\mathcal{Z} \},
\end{equation}
where
\begin{equation*}
\mathcal{U}(\Nr,\mathbb{R}^{\Nf}):= \text{St}(\Nr,\mathbb{R}^{\Nf})\cap \text{Sp}(\Nr,\mathbb{R}^{\Nf}),
\end{equation*}
is the manifold of orthosymplectic matrices, being $\St(\Nr,\r{\Nf}):= \{ U\in\R{\Nf}{\Nr}: U^{\top}U = I_{\Nr}\}$ the Stiefel manifold and
$\Simp(\Nr,\r{\Nf}):=\{ U\in\R{\Nf}{\Nr}: U^TJ_{\Nf}U=J_{\Nr}\}$
the manifold of symplectic matrices.
The factor $Z$ in \eqref{eq:dynManifold} is assumed to belong to the space
\begin{equation}\label{eq:coeff_space}
    \mathcal{Z}:= \{ Z\in\mathbb{R}^{\Nr\times \Np}: \text{ rank}(ZZ^T+J_{\Nr}ZZ^TJ_{\Nr})=\Nr \},
\end{equation}
with $\Nrh\ll \Nfh$ and $\Nr<\Np$. The full-rank condition, \eqref{eq:coeff_space}, guarantees the uniqueness of $U$ in the decomposition in \eqref{eq:dynManifold}, for all $t\in \mathcal{T}$, once $Z$ is fixed.
Let $X_r^{i}(t):=U_X(t) Z_i(t)\in\mathbb{R}^{\Nfh}$ and $V_r^{i}(t):= U_V(t) Z_i(t)\in\mathbb{R}^{\Nfh}$ denote the reduced velocity and position vectors, respectively, associated with the parameter $\prm_i$, for $i=1,\ldots,\Np$.
The submatrices $U_X, U_V\in\mathbb{R}^{\Nfh\times\Nr}$ are defined as $(U_X)_{i,j}=(U)_{i,j}$, $(U_V)_{i,j}=(U)_{i+\Nfh,j}$ for any $ i=1,\dots,\Nfh$ and $j=1,\dots,\Nr$.
For fixed $U$ and for each parameter $\prm_i \in \Gamma_h$, the flow map of the coefficient equation is a canonical symplectic map with a Hamiltonian having the $i$-th entry equal to
\begin{equation}\label{eq:nonlinear_hamiltonian} 
    \HamUi(Z_i(t)) :=\Ham(U(t)Z_i(t))=\frac12  V_r^{i}(t)^\top V_r^{i}(t) + \EE_{U,i}(Z_{i}(t)),
\end{equation}
where the first part is quadratic in the coefficients $Z_i$, and $\EE_{U,i}:\mathbb{R}^{\Nf}\rightarrow\mathbb{R}$ is the nonlinear electric energy component \eqref{eq:HamEEfull} of the Hamiltonian function, i.e., for all $i\in\Sp$,
\begin{equation}\label{eq:HamEE}
    \EE_{U,i}(Z_{i}(t)) = \EE(X^i_r;\prmi)= \dfrac{\mmp^{-1}}{2} \Mq^\top\Lambda^0(X_r^{i}(t)) L^{-1}\Lambda^0(X_r^{i}(t))^\top\Mq.
\end{equation}
We introduce the matrices
\begin{equation*}
\Ghp:=[\Ghpi{1}|\dots | \Ghpi{\Np}]\in\R{\Nf}{\Np},
\end{equation*}
and
\begin{equation*}
\ghp:=[\ghpi{1}|\dots | \ghpi{\Np}]\in\R{\Nr}{\Np}
\end{equation*}
having as columns the $p$ instances of the gradient of the Hamiltonian and of the reduced Hamiltonian, respectively,
\begin{equation}\label{eq:grads}
    \Ghpi{i}:=\nabla_{UZ_i}\HamUi(Z_i)\in\mathbb{R}^{\Nf},
    \qquad
    \ghpi{i}:=U^\top\Ghpi{i}= \nabla_{Z_i}\HamUi(Z_i)\in\mathbb{R}^{\Nr}.
\end{equation}
where
\begin{equation}\label{eq:grH-full}
     \nabla_{UZ_i}\HamUi(Z_i(t)) =
     \begin{pmatrix}
     \mmp^{-1}\mathrm{diag}(\Mq)\nabla\Lambda^0(X_r^i(t))L^{-1}\Lambda^0(X_r^i(t))^\top \Mq\\
     V_r^i(t)
     \end{pmatrix}.
\end{equation}
Under these assumptions, a dynamical system for the reduced solution is characterized via the symplectic projection of the velocity field of the full dynamical system \eqref{eq:full_param_HamSystem} onto the tangent space of $\Mred$ at each state \cite{P19,musharbash2017symplectic}.
The resulting reduced dynamics is given in terms of evolution equations for the reduced basis and for the expansion coefficients as
\begin{subequations}\label{eq:UZred}
\begin{empheq}[left = \empheqlbrace\,]{align}
	&\dot{Z}(t) = \J{\Nr}\ghp = \J{\Nr}\nabla_Z \HamU(Z), & \text{in }\Tcal, \label{eq:UZred_coeff} \\ 
	&\dot{U}(t) = (\Idm_{\Nf}-UU^\top)\big(\J{\Nf}\Ghp Z^\top -
	\Ghp Z^\top\J{\Nr}^\top\big)
	S(Z)^{-1}, & \text{in }\Tcal, \label{eq:UZred_basis} \\ 
	&U(t_0)Z(t_0) = U^0 Z^0,& \label{eq:UZred_initcond}
\end{empheq}
\end{subequations}
where $S(Z)=ZZ^\top+\J{\Nr}^\top ZZ^\top\J{\Nr}$, and
the initial condition $U^0Z^0\in\Mred$ is computed via a truncated complex SVD of $\Rcal_0\in\R{\Nf}{\Np}$. More details on the symplectic projection and on the gauge conditions enforced on the tangent space of $\Mred$ can be found in \cite{P19}.

Equation \eqref{eq:UZred_coeff} describes the evolution of the coefficients $Z(t)$ and is a system of $\Np$ independent equations, each in $\Nr$ unknowns. It corresponds to the Galerkin projection of the full order Hamiltonian systems onto the space spanned by the columns of $U(t)$, as obtained with the global symplectic reduced basis method. 
Here, however, the basis $U$ changes in time, and the evolution problem \eqref{eq:UZred_basis}, for the basis $U$, is a matrix-valued problem in $\Nf\times \Nr$ unknowns on the manifold of ortho-symplectic rectangular matrices. Observe that the reduced basis depends on the parameters, but it is the same for all parameters in the set $\Sprmh$.

\section{Efficient treatment of nonlinear terms}\label{sec:hyperred}
In this section, we discuss the computational cost of the numerical solution of the reduced problem \eqref{eq:UZred}, and propose a novel algorithm for the efficient and structure-preserving treatment of the nonlinear operators.

The proposed reduction of the nonlinear terms is independent of the numerical time integrators used to solve the reduced dynamical system \eqref{eq:UZred}.
However, the algorithm can be optimized depending on the time integrator of choice.
We consider the structure-preserving
partitioned Runge--Kutta temporal integrators proposed in \cite{P19,hesthaven2020rank}.
In particular, the evolution of the 
basis $U$ is approximated with an explicit method, while a symplectic temporal integrator is employed for the evolution of the coefficients $Z$, and this latter will generally be an implicit scheme.
 Observe that we do not require that the stages of the RK integrators for the basis $U$ and coefficients $Z$ coincide.
  We will discuss the details and implementation of such schemes in \Cref{sec:numExp}.
 Although not strictly necessary, here we also assume that the first step of the partitioned RK method involves the evolution of the reduced basis; this assumption implies that we have the information on the Hamiltonian gradient at the beginning of the temporal interval (at least for some parameter values).

Let us split the temporal domain $\Tcal$ into sub-intervals $\Tcalt:=(t_{\tau-1},t_{\tau}]$, for any $\tau=1,\dots,N_{\tau}$, where $t_{0}=0$ and $\Delta t = t_{\tau}-t_{\tau-1}$ is the uniform time step.
For each temporal interval $\Tcalt$, the dynamical reduced basis method involves the following operations.
\begin{itemize}
\item The evolution of the basis $U$ requires
$O(\Nfh\Nrh\Np)+O(\Nx^2\Np)+O(\Nfh\Np c)+O(\Nfh\Nrh^2)+O(\Nrh^2\Np)$ flops,
where $c\in\mathbb{N}$ is the number of finite element basis functions whose support is contained in a given mesh element, and we recall that $\Nfh$ is the number of particles, $\Nrh$ the size of the reduced basis and $\Np$ denotes the number of parameter values.
Note that $c$ is a mild constant, and it is equal to $2$ for piecewise polynomial functions in 1D, as in the discretization discussed in \Cref{sec:GEMPIC}.
The computational costs of this step are distributed as follows
\begin{center}
\begin{tabular}{c|c}
Arithmetic complexity  & Operation\\
\hline
$O(\Nfh\Nrh\Np)$
    & computation of $V_r=U_VZ$ and $X_r=U_XZ$ \\
$O(\Nfh\Np c)$
    & assembly of $\nabla\Lambda^0(X_r^i)$ and of $\Lambda^0(X_r^i)$, for all $i\in\Sp$ \\
$O(\Nfh\Np c)+O(\Nx^2\Np)$
    & computation of $\mathrm{diag}(\Mq)\nabla\Lambda^0(X_r^i(t))L^{-1}\Lambda^0(X_r^i(t))^\top \Mq$\\
$O(\Np\Nrh^2)+O(\Nrh^3)$
    & construction and inversion of the matrix $S(Z)$\\
$O(\Nfh\Nrh\Np)+O(\Nfh\Nrh^2)$
    & matrix-matrix multiplications in the r.h.s. of \eqref{eq:UZred_basis}
\end{tabular}
\end{center}
The first three rows of the table correspond to the
assembly and evaluation of $\Ghp$ for all $\Np$ parameters in $\Sprmh$.
\item The integration, using an implicit time scheme, of the evolution equation \eqref{eq:UZred_coeff} for the $\Np$ vector-valued coefficients requires $O(\Nfh\Nrh\Np)+O(\Nx^2\Np)+O(\Nfh\Np c)$ flops:
\begin{center}
\begin{tabular}{c|c}
Arithmetic complexity  & Operation\\
\hline
$O(\Nfh\Nrh)$
    & computation of $U_V^\top V_r^i=U_V^\top U_VZ_i$\\
$O(\Nfh\Nrh)+O(\Nx^2)+O(\Nfh c)$
    & assembly and evaluation of $\Ghpi{i}$\\
$O(\Nfh\Nrh)$
    & computation of $U_X^\top \Ghpi{i}$
\end{tabular}
\end{center}
Each of the operation listed in the table needs to be performed for each parameter $\prm_i\in\Sprmh$, at each stage of the RK scheme, and at each iteration of the nonlinear solver.
\end{itemize}

The leading computational cost in both steps depends on the product of the number of particles $\Nfh$ and of the number of parameter $\Np$, both potentially large in multi-query simulations of high-dimensional problems.
This cost is associated with the remapping of the particles to the full dimensional space, in each temporal interval and for each parameter, and with the evaluation of the velocity field of the reduced flow.
Indeed, the sole knowledge of the expansion coefficients with respect to the reduced basis is not enough to compute the particles-to-grid mapping needed to evaluate the electric field and, hence, the Hamiltonian. Even in the reduced model \eqref{eq:UZred}, these operations require the reconstruction of the approximate particle positions, at a cost proportional to the size of the full model. This lifting to the high-dimensional space needs to be performed for each instance of the parameter, at each stage of the RK time integrator, and at every iteration of the nonlinear solver.
Analogous computational problems are common in model order reduction and emerge whenever non-affine and nonlinear operators are involved, \emph{cf.}  e.g. \cite[Chapters 10 and 11]{QMN16}.
In the presence of low-order polynomial nonlinearities, tensorial techniques \cite{cstefuanescu2014comparison} can be used to separate terms that depend on full spatial variables and on reduced coefficients to allow efficient computations of the nonlinear terms.
The non-polynomial nature of the nonlinearity in the gradient of \eqref{eq:nonlinear_hamiltonian} prevents us from using the aforementioned tensorial approach to accelerate the computation. For general nonlinear operators, several so-called \emph{hyper-reduction} techniques have been proposed to mitigate or overcome the computational bottleneck associated with the lifting.
The empirical interpolation method (EIM) \cite{barrault2004empirical} and its discrete counterpart (DEIM) \cite{chaturantabut2010nonlinear} are interpolatory techniques used to approximate the nonlinearity in the projection-based ROM, requiring the computations of only a few components of the original nonlinearity. 
While effective in the case where each component of the nonlinearity depends only on a few components of the input, they are not suited for the treatment of not component-wise nonlinear terms. Using a sparsity argument \cite{chaturantabut2010nonlinear} or the introduction of auxiliary variables \cite{de2016nonlinear}, DEIM has been adapted to deal with the approximation of the nonlinear terms at interpolation points that require the evaluation of the reduced solution on a limited number of neighboring mesh points, as it happens for high-order spatial discretization schemes with large stencils. 
However, the same strategy does not work for the treatment of the gradient of \eqref{eq:nonlinear_hamiltonian}, as the inverse of the discrete Laplacian operator is generally dense, and hence each of its entries requires the computation of $X_r^i$ for $p$ sampled parameters, making traditional approaches computationally impractical.
Moreover, traditional hyper-reduction techniques, like EIM and DEIM, applied to a gradient vector field do not result in a gradient field, which means that the geometric structure of the Hamiltonian dynamics is compromised in the hyper-reduction process.

To achieve computational efficiency in the simulation of \eqref{eq:UZred} without compromising its geometric structure,
we propose a strategy that 
approximates the reduced Hamiltonian gradient via a combinated hyper-reduction technique and sampling procedure. A reduction in the computational runtimes of the algorithm is achieved by decoupling the operations that depend on $\Nfh$ from those that depend on $\Np$, while retaining an accurate representation of the plasma dynamics.
There are several challenges that we need to face in the development of such techniques:
(i) the preservation of the Hamiltonian structure of the dynamics;
(ii) the lack of information on the full model solution and nonlinear operators, traditionally collected in an offline phase via snapshots;
and (iii) the lack of a sparsity pattern in the nonlinear Hamiltonian gradient, i.e., the fact that each entry of the electric energy vector \eqref{eq:HamEE} depends on all $\Nfh$ computational particles.

\subsection{Parameter sampling} \label{sec:pSampl}
The reduced dynamics \eqref{eq:UZred} involve $\Np$ evolution equations for the expansion coefficients, one per parameter value, and one evolution problem for the matrix-valued reduced basis. Since, at each time, the reduced basis is the same for all parameters, one can reduce the computational cost required for its evolution by sampling over the parameter space and constructing a reduced basis for only a subsample of parameters, but which remains accurate for all other parameters in $\Sprmh$.
This corresponds to a reduction in parameter space.
Let us denote by $\Np$ the cardinality of the set $\Sprmh$ and assume that the parameters in $\Sprmh$ are indexed from the set $\Sp:=\{1,\ldots,\Np\}$.
Let us consider a subset $\prmhsub$ of $\Sprmh$ of size
$\Nps\ll\Np$. Define $\Spsub\subset \Sp$ to be the set of indices corresponding to the parameters in the selected subset so that
$\prmhsub = \left \{ \prmi \in \Sprmh | \; i\in \Spsub \right \}\subset\Sprmh$.
The idea of the proposed sampling approach is to replace, in the evolution of the basis \eqref{eq:UZred_basis}, the matrix of the expansion coefficients $Z(t)\in\R{\Nr}{\Np}$ by the matrix obtained via the concatenation of the columns of $Z$ with indices in $\Spsub$. Following the discussion at the beginning of the section, this approximation leads to a computational complexity for the basis evolution of the order of $O(\Nfh\Nrh\Nps)+O(\Nx^2\Nps)+O(\Nfh\Nps c)+O(\Nfh\Nrh^2)+O(\Nrh^2\Nps)$.
To preserve the accuracy of the method, we must ensure that the chosen subset $\prmhsub$ is representative of the entire parameter set $\Sprmh$.
For the sake of simplicity, in this work, we set it at $t=0$, and we keep it fixed over time. Starting from $\prmhsub=\emptyset$, the set $\prmhsub$ is constructed using a greedy algorithm that, at each iteration, adds to the index subset $\Spsub$ the index $i$ that satisfies
\begin{equation}\label{eq:pSampl}
    \max_{i\in\Sprmh\setminus\prmhsub}\min_{j\in\Sprmh}\| Z^0_i-Z^0_j\|_2,
\end{equation}
until a user-defined threshold value of the cost function or a maximum number of iterations.
Possible research directions to improve the selection strategy would be to adapt in time the set $\prmhsub$ to capture significant changes in the behavior of the solution relative to the parameters or to modify the cost function to incorporate errors in the evaluation of physical quantities, such as the electric field. 
However, as \eqref{eq:pSampl} is an NP-complete problem \cite{shitov2021column}, it is not currently known if it is possible to find an optimal parameter selection strategy with a polynomial cost in $\Np$. Thus, in the current form, while affordable if performed only once at $t=0$, the selection strategy may become computationally expensive for large $\Np$ and could compromise the efficiency of the proposed dynamical RB method if repeated in each temporal sub-interval.
Suboptimal, yet more efficient, algorithms could be adopted by framing the problem into the more general column subset selection problem (CSSP) \cite{boutsidis2009improved}, that consists in finding an optimal subset of columns of a given matrix that minimizes the residual of the projection of the given matrix onto the selected column subset.
Other parameter reduction strategies, like active subspaces \cite{Constantine15},
might also be envisioned.

 Concerning the expansion coefficients, for which one differential equation per parameter needs to be solved, subsampling is not an option.


\subsection{DMD-DEIM approximation of the Hamiltonian gradient}

In this section, we develop a reduction algorithm where, in each temporal interval $\Tcalt$, DMD is used for the hyper-reduction of the electric potential $\Phi(X_r^i(t))=L^{-1}\Lambda^0(X_r^{i}(t))^\top \Mq$ in \eqref{eq:HamEE}, while a DEIM strategy is developed to approximate the component $\Lambda^0(X_r^{i}(t))$ of the particles-to-grid mapping.
Note that in \eqref{eq:nonlinear_hamiltonian}, the quadratic term involving the particles' velocity represents a linear contribution in the gradient of the reduced Hamiltonian and, hence, does not require any hyper-reduction.

\subsubsection{Dynamic Mode Decomposition of the electric potential}\label{sec:DMD}
Dynamic mode decomposition is an equation-free data-driven approach, proposed in \cite{schmid2010dynamic,schmid2011application}, that
uses only data measurements of a given dynamical system to approximate the dynamics and predict future states.
The idea is to decompose the problem into a set of coherent spatial structures, known as DMD modes, and associate correlated data to specific Fourier modes that capture temporal variations. DMD was initially employed as a spectral decomposition method for complex fluid flows \cite{rowley2009spectral}. More recently, it has proved successful in a wide range of settings such as background/foreground separation in real-time video \cite{grosek2014dynamic}, characterization of dynamic stall \cite{mariappan2014analysis},
and analysis of the propagation of infectious diseases \cite{proctor2015discovering}.
DMD hinges on the theory of Koopman operators \cite{koopman1931hamiltonian}, which allows to represent the flow of a nonlinear dynamical system via an infinite-dimensional linear operator on the space of measurement functions.
Without explicit knowledge of the operator describing the dynamics, DMD computes a least-squares regression of data measurements to an optimal finite-dimensional linear dynamical system that approximates the infinite-dimensional Koopman operator.
In this subsection, we first describe the classical DMD algorithm following \cite{DMDbook16}.
Next, we introduce a sliding-window based DMD formulation for the hyper-reduction of the electric potential in the dynamical reduced model \eqref{eq:UZred} of the Vlasov--Poisson problem.

Consider a general nonlinear dynamical system: Find $\ybf:\Tcal\rightarrow \r{\ell}$, for $\ell\in\mathbb{N}$ such that
\begin{equation}\label{eq:NLsys}
    \left\{
    \begin{aligned}
        & \dot{\ybf}(t)= F(t,\ybf(t)), & \quad t\in\Tcal,\\
        & \ybf(t_0)=\ybf_0. &
    \end{aligned}\right.
\end{equation}
Assume that we have as data measurements exact values or approximations of the state at different time instants, namely
\begin{equation}\label{eq:DMDdata}
    \Ybf = 
    \begin{bmatrix}
    \ybf_{0} & \ybf_{1} & \dots & \ybf_{\tau-1}
    \end{bmatrix}\in\R{\ell}{\tau},
    \qquad
    \Ybf' = 
    \begin{bmatrix}
    \ybf_{1} & \ybf_{2} & \dots & \ybf_{\tau}
    \end{bmatrix}\in\R{\ell}{\tau},    
\end{equation}
where $\ybf_k=\ybf(k\Delta t)$ and $\Delta t$ is the uniform time step.
In the DMD method, data measurements are used to approximate the nonlinear dynamics \eqref{eq:NLsys}
by a locally linear system $\dot{\ybf}=A\ybf$,
where $A\in\R{\ell}{\ell}$ is the matrix that best fits the measurements in a least-square sense, i.e., $A=\argmin_B{\norm{\Ybf'-B\Ybf}_F}$, with $\norm{\cdot}_F$ being the Frobenius norm. Then, $A$ is given by $A=\Ybf'\Ybf^{\dagger}$,  where $\dagger$ denotes the Moore-Penrose pseudoinverse.

From the linear approximation of the dynamics, the DMD algorithm computes a low-rank eigendecomposition of
the matrix $A$ by extracting its $\Ndmdt$ largest eigenvalues $\Lambda^A$ and corresponding eigenvectors $\Theta^A=[\theta^A_1 \dots \theta^A_{\Ndmdt}]\in\R{\ell}{\Ndmdt}$. The resulting DMD approximation of the state $\ybf(t)$, for $t>\tau \Delta t$, reads
\begin{equation}\label{eq:DMD_approximation}
    \ybf(t) \approx \ybf_{\dmd}(t) = 
    \Theta^{A}\left( \Pi \odot e^{\Omega(t-\tau \Delta t)}\right)
    = \sum^{\Ndmdt}\limits_{j=1} \theta_j^A \pi_j
    e^{\omega_j(t-\tau\Delta t)},
\end{equation}
where, for any $j=1,\ldots,\Ndmdt$, $\omega_j:=\ln (\Lambda^A_j) / \Delta t$ is the $j$-th entry of the vector $\Omega\in\r{\Ndmdt}$, while $\pi_j$
is the $j$-th entry of the vector $\Pi=(\Theta^A)^{\dagger}\ybf_{0}\in\r{\Ndmdt}$ containing the coordinates of the initial condition $\ybf_{0}$ with respect to the DMD modes.

If the size of the matrix $A$
is large, $A$ might be severely ill-conditioned and not directly tractable. In this situation, a different version of the DMD algorithm, proposed in \cite{tu2013dynamic}, projects the data into a low-rank subspace
instead of deriving $A$ directly from the data, as described in Algorithm \ref{alg:DMD_algorithm}. Moreover, given the sensitivity of the DMD algorithm to the duration and sampling of the series $\Ybf$ and $\Ybf'$, \cite{dylewsky2019dynamic} proposes a sliding-window approach where the measurement data are not taken in the whole temporal interval but only in the sampling window
$[t_{\tau-\T},t_{\tau}]$ of length $\T\Delta t$, with $\T\in\mathbb{N}$. The rationale is that if the system is time-varying and the incoming data is harvested in a streaming fashion, it may be beneficial to accuracy and memory storage to consider only the most recent data. The only computational overhead is the computation of the DMD modes and weights in the DMD approximation \eqref{eq:DMD_approximation} as new data are collected. This cost may be mitigated by efficient online updates of the eigenvalues and eigenvectors of A \cite{hemati2014dynamic} or by means of incremental SVD algorithms \cite{matsumoto2017fly}.

\begin{algorithm}[t]
\SetAlgoLined
\SetKwData{Left}{left}\SetKwData{This}{this}\SetKwData{Up}{up}
\SetKwFunction{Union}{Union}\SetKwFunction{FindCompress}{FindCompress}
\SetKwInOut{Input}{Input}\SetKwInOut{Output}{Output}
\caption{DMD algorithm}
\label{alg:DMD_algorithm}
\Input{ $\Ybf$, $\Ybf'$, $\tol$.}
\Output{ $\Lambda$, $\Theta$.}
\BlankLine
Compute the truncated SVD of $\Ybf$, $\Ybf=U\Sigma V^{\top}$, using \emph{tol} as tolerance for singular values selection\;\label{step:1}
Define $A_{\tol}=U^{\top}\Ybf'V\Sigma^{-1}$\;\label{step:2}
Compute the eigendecomposition of $A_{\tol}$: $A_{\tol}W=W\Lambda$\;\label{step:3}
Reconstruct the eigendecomposition of $A$ by defining is eigenvectors as $\Theta=\Ybf'V\Sigma^{-1}W$.\label{step:4}
\end{algorithm}

In the context of kinetic plasma PIC simulations, a DMD strategy has been used
in \cite{nayak2020reduced} to detect and track equilibrium states.
The aforementioned method relies on snapshots of the high-fidelity simulation until an equilibrium is detected and, after this time, the solution is extrapolated via the DMD modes.
Here, we propose to employ a DMD strategy in a different way, namely to hyper-reduce the self-consistent electric potential
$\Phi(X_r^i(t))=L^{-1}\Lambda^0(X_r^i(t))^\top \Mq\in\r{\Nx}$ that enters the reduced Hamiltonian \eqref{eq:nonlinear_hamiltonian} for each parameter $\prmi\in\Sprmh$.
The idea is to extract low-dimensional dynamical features from a time-series of the electric potential and use them, as part of the DMD algorithm in \eqref{eq:DMD_approximation}, to extrapolate the value of $\Phi(X_r^i(t))$ needed for the computation
of the reduced Hamiltonian \eqref{eq:nonlinear_hamiltonian} in each temporal interval $\Tcalt$.
In details, let $\Phi_{\tau}^{i}$, for a fixed parameter with index $i \in \Spsub$ and $\tau=1,\dots,N_{\tau}$, be the approximation of $\Phi(X_r^i(t))$ at $t=t_{\tau}$ for each parameter in the subset $\prmhsub$.
Since the first step of the temporal integrator involves the evolution of the reduced basis,
these quantities are computed while assembling the right hand side of the basis evolution equation \eqref{eq:UZred_basis}.
For each $\Tcalt$, we collect the time-discrete approximations of the electric potential obtained in a {past} time window of length {$(\T+1)\Delta t$};
\begin{equation}\label{eq:DMD_potential_Y_time}
  \Ybf_i = \begin{bmatrix} \Phi_{\tau-\T-1}^{i} & \Phi_{\tau-\T}^{i} & \cdots & \Phi_{\tau-2}^{i} \end{bmatrix},
  \qquad
  \Ybf'_i = \begin{bmatrix} \Phi_{\tau-\T}^{i} & \Phi_{\tau-\T+1}^{i} & \cdots & \Phi_{\tau-1}^{i} \end{bmatrix},
\end{equation}
where $\Ybf_{i}, \Ybf_{i}'\in\mathbb{R}^{\Nx\times \T}$ for $i=1,\dots, \Nps$.
Extracting the dominant modes from each
realization of the electric potential associated with a fixed parameter is a cumbersome task. 
To the best of our knowledge, DMD-based methods for the model order reduction of parametric problems have not been developed.
In our setting, the dependence on the parameter comes from the state, and it is propagated via the parametric initial distribution $f_0$. This suggests that, instead of extracting the DMD modes for each fixed parameter $\prm_i\in \prmhsub$, we can incorporate the parameter in the DMD procedure to approximate the dynamics of the electric potential for all parameters. A similar approach can be found in \cite{sayadi2015parametrized} when dealing with bifurcation parameters in thermo-acoustic systems.
For each parameter index $i\in\Sp$, the datasets $\Ybf_i$ and $\Ybf'_i$ are concatenated column-wise to form two global datasets $\Ybf$ and $\Ybf'$, i.e.
\begin{equation}\label{eq:DMD_potential_Y_time_parameter}
  \Ybf = \begin{bmatrix} \Ybf_1 & \Ybf_2 & \cdots & \Ybf_{\Nps} \end{bmatrix},
  \qquad
  \Ybf' = \begin{bmatrix} \Ybf_1^{\prime} & \Ybf_2^{\prime} & \cdots & \Ybf_{\Nps}^{\prime} \end{bmatrix} ,
\end{equation}
with $\Ybf,\Ybf'\in\R{\Nx}{\Nps\T}$.  This procedure is justified by the absence of an explicit dependence of the electric potential on the parameter.
Following Algorithm \ref{alg:DMD_algorithm}, we generate the DMD eigenvectors $\Theta\in\R{\Nx}{\Ndmdt}$ and eigenvalues $\Lambda\in\R{\Ndmdt}{\Ndmdt}$ of the linear approximation of the dynamics for the problem of interest. The resulting DMD approximation of the self-consistent electric potential reads
\begin{equation}\label{eq:DMD_approx_potential_subset}
    \Phi(X_r^i(t)) \approx \Phi^i_{\dmd}(t) = \Theta\left( \Pi_i \odot e^{W(t-(\tau-1)\Delta t)}\right), \qquad \forall\, i\in\Spsub, \forall\, t\in \Tcalt,
\end{equation}
with $W\in\r{\Ndmdt}$ is the vector of entries $\omega_j=\ln(\Lambda_j)/\Delta t$, for any $j=1,\ldots,\Ndmdt$, and $\Pi_i:=\Theta^{\dagger}\Phi^i_{\tau-1}\in\r{\Ndmdt}$ for any $i\in\Spsub$.

Assuming a smooth dependence of the DMD coordinates $\Pi_i$ on the parameter, interpolation techniques can be used to recover the DMD coordinates for parameters not included in $\prmhsub$, similarly to the POD with interpolation (PODI) \cite{bui2003proper}. In this work, we adopt the radial basis interpolation \cite{broomhead1988radial}, with a Gaussian kernel, as interpolation algorithm to reconstruct the DMD coordinates $\Pi_i$ for $i\in \Sprmh\setminus\prmhsub$.
The computational cost of the interpolation step is negligible as compared to the cost of \Cref{alg:DMD_algorithm}, as we comment on at the end of the section.
For ease of the notation, we use the same symbol $\Pi_i$ to represent the interpolated DMD coefficients for all $i\in \Sp$. Knowing $\Pi_i$ for all $i\in\Sp$, the electric potential is reconstructed using \eqref{eq:DMD_approx_potential_subset}.
The resulting sampling error can be controlled by enriching the subset of parameters $\prmhsub$ and by optimal placement of the location of the parameters with indices in $\Spsub$ in the parameter space.

Using the DMD estimate of the potential $\Phi$, the Hamiltonian function \eqref{eq:nonlinear_hamiltonian} is approximated as
\begin{equation}\label{eq:Hamiltonian_DMD}
    \HamUi^{\dmd}(Z_i)=\dfrac12  V_r^{i}(t)^\top V_r^{i}(t) +
        \dfrac{\mmp^{-1}}{2} \Mq^\top\Lambda^0(X_r^{i}(t)) \Phi^i_{\dmd}(t), \qquad \forall\, i\in \Sp,\;\forall\, t\in \Tcalt.
\end{equation}
\begin{remark}
If the Hamiltonian function depends explicitly on the parameter, the approach outlined above is not legitimized because a non-parametric operator would be used to approximate the parametric potential. An alternative strategy would require an approximation of the form \eqref{eq:DMD_approx_potential_subset} for each parameter realization $\prm_i$, with $i\in S_{\Nps}$, using different $W_i$ and $\Theta_i$ for each $i$. The resulting DMD approximations of the potential $\Phi^i_{\dmd}(t)$ could then be directly interpolated on $\Sp$ or, as suggested in \cite{druault1999developpement}, interpolated based on physical concepts as in PODI.
\end{remark}

The DMD approximation based on the matrix $A$, instead of its projection $A_{\tol}$ defined in \Cref{alg:DMD_algorithm},
would result in a computational complexity
$O(\Nx^3)$. Although this cost might still be
tractable in one dimension, it becomes prohibitive when considering the Vlasov--Poisson problem in higher dimension. The method described in \Cref{alg:DMD_algorithm} is, therefore, the preferred choice.

The computational cost of the proposed DMD strategy, applied to the electric potential, reduces to the cost needed to perform \Cref{alg:DMD_algorithm} from the datasets $\Ybf,\Ybf'\in\R{\Nx}{\Nps\T}$ in \eqref{eq:DMD_potential_Y_time_parameter}.
The truncated SVD decomposition of $\Ybf$ in \Cref{step:1}
has arithmetic complexity $O(\Nx\Nps\T\Ndmdt)$, where $\Ndmdt$ is the number of retained modes \cite{Martinsson11}.
Observe that, if the number $\Ndmdt$ of truncated modes is chosen based on a tolerance to control the magnitude of the neglected singular values, \Cref{alg:DMD_algorithm} computes the full SVD of $\Ybf$ and then performs the truncation. This variant of \Cref{step:1} has computational complexity
$O(\Nx(\Nps\T)^2)$, under the assumption that {the number $\T$ of chosen DMD samples} and the number of parameter subsamples $\Nps$ satisfy $\Nx>\Nps\T$.
The eigendecomposition of $A_{\tol}\in\R{\Ndmdt}{\Ndmdt}$ in \Cref{step:3} costs $O(\Ndmdt^3)$.
Finally, the matrix-matrix multiplications to compute $A_{\tol}$ and $\Theta$ in \Cref{step:2} and \Cref{step:4}, respectively, require
$O(\Nx\Ndmdt^2)+O(\Nx\Nps\T\Ndmdt)$ operations.
The computational cost to compute $\Phi^i_{\dmd}$ for every parameter $\prmi\in\Sprmh$ -- including sampling parameters and reconstructed parameters -- is $O(\Nx\Ndmdt\Np)$.
The leading cost is, therefore, $O(\Nx\Ndmdt\Np) + O(\Nx\Nps\T\Ndmdt)$, with the last term replaced by $O(\Nx(\Nps)^2\T^2)$ for a naive implementation of the truncated SVD. This cost is linear in $\Nx$, does not depend on the number $\Nfh$ of particles,
and only the computation of the DMD coordinates $\Pi_i$ depends on the number of parameters $\Np$.

\subsubsection{Discrete Empirical Interpolation Method for reduction in the number of particles}\label{sec:deim}
The DMD approach described in the previous section allows to derive an approximate electric potential that can be evaluated independently on the number of particles. However, the evaluation of the electric energy component of the approximate reduced  Hamiltonian \eqref{eq:Hamiltonian_DMD} still requires the particles-to-grid mapping for $\Lambda^0(X_r^{i}(t))$, for each value of the parameter, at each stage of the temporal solver, and at each iteration of the nonlinear solver. The computational cost of this step is a major bottleneck of the algorithm. To overcome this computational burden, we propose hyper-reduction of the approximate reduced  Hamiltonian \eqref{eq:Hamiltonian_DMD} with a DEIM-based strategy. 

The DEIM approach is a discrete variant of the empirical interpolation method
(EIM) introduced in \cite{barrault2004empirical} to approximate nonlinear functions via a combination of projection and interpolation. DEIM constructs carefully selected interpolation
indices to specify an interpolation-based projection so that the complexity of
evaluating the nonlinear term becomes proportional to the (small) number of selected
spatial indices.
The gist of the DEIM approximation of a given nonlinear operator $F:\r{\Nfh}\rightarrow\r{\Nfh}$ is to replace it with its oblique projection into a lower-dimensional space obtained from the most relevant components of $F$.
Starting from snapshots of the nonlinear function $F$, DEIM
constructs a lower-dimensional space as the span of the left singular vectors associated with the $\Ndeim\ll\Nfh$ largest singular values of the snapshot matrix.
If $\Psi=[\psi_1 \dots \psi_{\Ndeim}]\in\mathbb{R}^{\Nfh\times \Ndeim}$ denotes the DEIM basis, 
the DEIM approximation of $F$ reads
\begin{equation}\label{eq:DEIM_example}
    F(\xbf) \approx \Psi (P^{\top}\Psi)^{-1}P^{\top}F(\xbf),\qquad\forall\,\xbf\in\r{\Nfh},
\end{equation}
where $P\in\R{\Nfh}{\Ndeim}$ is defined as $P:=[\mathbf{e}_{\ell_1} \dots \mathbf{e}_{\ell_{\Ndeim}}]$, with $\mathbf{e}_k=\{0,1\}^{\Nfh}$ being the $k$-th canonical vector. We denote by $I_{\deim}$ the set of indices corresponding to the DEIM sampling points so that $I_{\deim}=\{ \ell_k \}_{k=1}^{\Ndeim}$.  The DEIM sampling indices are iteratively chosen according to the greedy procedure proposed in \cite[Algorithm 1]{chaturantabut2010nonlinear}. The $k$-th step, for $1\leq k\leq \Ndeim$, of the algorithm consists in choosing the index of the $k$-th DEIM vector $\psi_k$ that maximizes the DEIM projection error of $\psi_k$ into the current DEIM basis $[\psi_1 \dots \psi_{k-1}]\in\R{\Nfh}{(k-1)}$.
Under the assumption that each entry of the nonlinear vector-valued function $F$ depends only on a few entries of the input vector $\xbf$, the computation of $P^\top F(\xbf)$ requires only $\Ndeim$ evaluations of the nonlinear function leading to a computational cost independent of the original dimension~$\Nfh$.

The application of the classical DEIM procedure for the hyper-reduction of the nonlinear Hamiltonian gradient \eqref{eq:grH-full} is challenged by several factors.
As stated at the beginning of the section, applying the DEIM interpolation directly to the right-hand side of the coefficients evolution equations \eqref{eq:UZred_coeff} arising from the dynamic reduced basis approach, would not result in a structure-preserving approximation. Moreover,
the classical DEIM algorithm hinges on the availability of snapshots of the full model nonlinear operator of interest collected in the offline phase. In our dynamical model order reduction approach, there is no offline phase and, therefore, snapshots are not available.

We consider the Hamiltonian splitting in \eqref{eq:nonlinear_hamiltonian}, and the
approximation of the reduced electric energy \eqref{eq:HamEE} resulting from DMD, namely
\begin{equation*}
    \EE_{U,i}^{\dmd}(Z_{i}(t),t) := \dfrac{\mmp^{-1}}{2} \Mq^\top\Lambda^0(X_r^{i}(t)) \Phi^i_{\dmd}(t), \qquad \forall i\in \Sp,\; \forall\, t\in \Tcalt,
\end{equation*}
where $\Phi^i_{\dmd}$ is defined in \eqref{eq:DMD_approx_potential_subset}.
Approximating directly the vector $\nabla \EE_{U,i}^{\dmd}$ by a DEIM interpolation, as in \eqref{eq:DEIM_example}, would not preserve the geometric structure of the problem because it is not possible to define explicitly an Hamiltonian gradient from the interpolated vector field.
We propose a DEIM approximation of the reduced electric energy via hyper-reduction of the
term $\Lambda^0(X_r^{i}(t))\in\R{\Nfh}{\Nx}$, which otherwise would require the evaluation of the finite element basis functions at each particle position.

Let us introduce the function
$\Ncal_i(X_r^{i}(t),t):=\Lambda^0(X_r^{i}(t))\Phi^i_{\dmd}(t)\in\r{\Nfh}$;
approximated using a DEIM approach in each temporal interval $\Tcalt$ as follows.
First, we consider {samples} of the nonlinear term
associated with the electric potential $\Phi$ \eqref{eq:PoissonSystem} at $\Nps$ instances of the parameter and over a temporal window of length {$(\T+1)\Delta t$}.
{These are obtained from simulations of the reduced dynamical system in the past window of interest.}
The {data} matrix $\Ybf\in\R{\Nfh}{\Nps(\T+1)}$ is defined as
\begin{equation}\label{eq:Ydeim}
\Ybf = \begin{bmatrix} \Ybf_1 & \Ybf_2 & \cdots & \Ybf_{\Nps} \end{bmatrix},\quad
\Ybf_i:=\begin{bmatrix} \Lambda^0(X_r^{i}(t_{\tau-\T-1}))\Phi^i_{\tau-\T-1} &  \cdots & \Lambda^0(X_r^{i}(t_{\tau-1}))\Phi^i_{\tau-1} \end{bmatrix}.
\end{equation}
Note that the terms $\Lambda^0(X_r^{i}(t_{\tau-j-1}))$ and $\Phi^i_{\tau-j-1}$, for $i\in\prmhsub$ and $j=0,\ldots,\T$, are available from the evolution equation \eqref{eq:UZred_basis} for the reduced basis
solved at previous time steps.
The DEIM basis matrix $\Psi_{\tau} \in\R{\Nfh}{\Ndeim}$ is obtained by taking the first $\Ndeim$ left singular vectors of the {matrix} $\Ybf$, where the value $\Ndeim$ is fixed at the beginning of the simulation and might differ for different problems. We will comment on this in the numerical experiments in \Cref{sec:numExp}.
Denoting with $P_{\tau}\in\R{\Nfh}{\Ndeim}$ the matrix corresponding to the DEIM indices obtained as describe above, the nonlinear term $\Ncal_i(X_r^i(t),t)$ is approximated by
\begin{equation*}
   \Psi_\tau^{\top} (P_{\tau}^{\top} \Psi_{\tau})^{-1} P_{\tau}^{\top} \Ncal_i (X_r^{i}(t),t), \qquad \forall\, i\in \Sp,\; \forall\, t\in \Tcalt.
\end{equation*}
Observe that, although the basis $\Psi_{\tau}$ is constructed from the parameter subsample $\prmhsub$, the nonlinear term $\Ncal_i$ is approximated by its DEIM projection onto the DEIM space for all instances of the parameter, i.e., for all $i\in\Sp$.

To reduce the computational burden associated with the computation of the DEIM sampling points $P_{\tau}$ in each temporal interval $\Tcalt$, we follow an update strategy similar to the one proposed in \cite{peherstorfer2015online}. All the interpolation indices in the set $I_{\deim}$ are computed using the standard DEIM greedy method; not at all time steps but only every $\DEIMfreq>1$ time steps. In other temporal intervals, we proceed as follows. Assume we have computed the set of DEIM indices $I_{\deim}^{\tau-1}$ in the temporal interval $\Tcal_{\tau-1}$, then, in the following interval $\Tcalt$, we update only the indices in the subset $I^{*}\subset I^{\tau-1}_{\deim}$ of cardinality $\nDEIMup$ given by
\begin{equation}\label{eq:indices_adaptive_DEIM}
    I^{*}=\underset{
    \substack{
    I\subset I_{\deim}^{\tau-1},\\
    \text{dim}(I)=\nDEIMup
    }
    }
    {\text{argmax}}
    \sum_{k\in I}
    (\psi_k^{\tau})^\top\psi_k^{\tau-1}\,,
\end{equation}
where $\psi_k^{\tau}$ denotes the $k$-th vector of the DEIM basis $\Psi_{\tau}$ at time $t_{\tau}$. The remaining $\Ndeim-\nDEIMup$ indices in $I_{\deim}^{\tau-1}\backslash I^{*}$ are inherited by $I_{\deim}^{\tau}$. The rationale for the choice of $I^*$ is to only update the indices associated with the DEIM basis vectors at $t_{\tau-1}$ that have undergone the largest rotations in the DEIM basis update from $\Psi_{\tau-1}$ to $\Psi_{\tau}$.

The resulting approximate reduced Hamiltonian, associated with the parameter $\prm_i\in\Sprmh$, reads
\begin{equation}\label{eq:Hdd}
    \Ham^{\dd}_{U,i}(Z_i(t),t) = \dfrac12  V_r^{i}(t)^\top V_r^{i}(t) + \dfrac{\mmp^{-1}}{2} \Mq^\top
        \Psi_{\tau} (P_{\tau}^\top\Psi_{\tau})^{-1}P_{\tau}^{\top} \Lambda^0(X_r^{i}(t))\Phi^i_{\dmd}(t), \quad \forall i\in \Sp,\; \forall t\in\Tcalt.
\end{equation}
Observe that the multiplication of
the matrix $\Lambda^0(X_r^i)$, of the finite element basis functions evaluated at the particles' position,
by the DEIM sampling matrix $P^{\top}_{\tau}$,
corresponds to evaluating the finite element basis functions only on a subset of $\Ndeim\ll \Nfh$ particles.
Hence, this operation represents a substantial reduction in the number of particles.

The computational cost of the DEIM algorithm can be summarized as follows.
The computation of the {data matrix} in \eqref{eq:Ydeim} only involves the multiplications of the terms $\Lambda^0(X_r^{i}(t_{\tau-j-1}))$ and $\Phi^i_{\tau-j-1}$, for $i\in\prmhsub$ and $j=0,\ldots,\T$. Indeed, since these terms are available from the solution of the reduced basis evolution at previous time steps, there is no cost associated with their assembly, at least at this stage of the proposed model order reduction algorithm.
The matrix-matrix multiplications require $O(\Nfh\Nps\T c)$, where $\Nps$ is the dimension of the subset of sampling parameters, {$\T$ is the number of samples in the temporal window},
and $c$ is a mild constant that depends only on the support of the finite element basis functions.
The truncated SVD decomposition of the {matrix} $\Ybf\in\R{\Nfh}{\Nps(\T+1)}$ has arithmetic complexity $O(\Nfh\Nps\T \Ndeim)$, where $\Ndeim$ is the number of DEIM modes. 
%
The computational cost required to assemble the interpolation matrix $P_{\tau}\in\R{\Nfh}{\Ndeim}$ using \cite[Algorithm 1]{chaturantabut2010nonlinear} only depends on $\Ndeim$. This cost is further reduced by updating the indices according to the strategy described above and inspired by the adaptive sampling of  \cite{peherstorfer2015online}.
Hence, the leading computational cost of the DEIM algorithm is $O(\Nfh\Nps\T d)$.

\subsection{DMD-DEIM reduced dynamics and computational complexity}\label{sec:DMD_DEIM_complexity}

From \eqref{eq:Hdd}, the Hamiltonian gradient $\Ghpi{i}$ in \eqref{eq:grads} is approximated as
\begin{equation}\label{eq:grHDD-full}
     \Ghddi{i}:=\nabla_{UZ_i}\Ham^{\dd}_{U,i}(Z_i(t),t) =
     \begin{pmatrix}
     \mmp^{-1}\mathrm{diag}\big(\nabla\Lambda^0(X_r^i(t))\Phi^i_{\dmd}(t)\big)P_{\tau} (P_{\tau}^\top\Psi_{\tau})^{-\top} \Psi_{\tau}^{\top} \Mq\\
     V_r^i(t)
     \end{pmatrix},
\end{equation}
for all $i\in \Sp$, and $t\in\Tcalt$.
Similarly, the approximation of the gradient of the reduced Hamiltonian $\ghpi{i}$ reads
\begin{equation}\label{eq:grHDD}
\begin{aligned}
    \ghddi{i}&=U(t)^\top \Ghddi{i}= \nabla_{Z_i}\Ham^{\dd}_{U,i}(Z_i(t),t)\\
    & =    U_V(t)^{\top}V_r^{i}(t) + 
    \mmp^{-1}U_X(t)^{\top}\mathrm{diag}\big(\nabla\Lambda^0(X_r^i(t))\Phi^i_{\dmd}(t)\big)P_{\tau} (P_{\tau}^\top\Psi_{\tau})^{-\top} \Psi_{\tau}^{\top} \Mq,
\end{aligned}
\end{equation}
for all $i\in \Sp$, and $t\in\Tcalt$.

The reduced dynamical system \eqref{eq:UZred} is approximated by
replacing the gradient of the reduced Hamiltonian $\ghpi{i}\in\mathbb{R}^{\Nf}$ with its DMD-DEIM approximation $\ghddi{i}$ from \eqref{eq:grHDD} in the evolution equations of the expansion coefficients.
The DMD-DEIM reduced dynamics reads
\begin{subequations}\label{eq:UZredhyp}
\begin{empheq}[left = \empheqlbrace\,]{align}
	&\dot{Z}(t) = \J{\Nr}\ghdd, & \text{in }\Tcal, \label{eq:UZhypred_coeff} \\ 
	&\dot{U}(t) = (\Idm_{\Nf}-UU^\top)(\J{\Nf}\Ghps Z^\top{-}
	\Ghps Z^\top\J{\Nr}^\top)
	S(Z)^{-1}, & \text{in }\Tcal, \label{eq:UZhypred_basis} \\ 
	&U(t_0)Z(t_0) = U^0 Z^0,&
\end{empheq}
\end{subequations}
Note that this approximate reduced model retains the geometric structure of the full model.

We analyze the computational cost to assemble and evaluate the right hand side of the DMD-DEIM reduced model \eqref{eq:UZredhyp} at each time instance.
We then compare the results with the ones at the beginning of \Cref{sec:hyperred} corresponding to the reduced model.
The evolution of the reduced basis requires $O(\Nfh\Nrh\Nps)+O(\Nx^2\Nps)+O(\Nfh\Nps c)+O(\Nfh\Nrh^2)+O(\Nrh^2\Nps)$ flops owing to the parameter sampling discussed in \Cref{sec:pSampl}. This cost also includes the assembly of the quantities
$\Lambda^0(X_r^i)$ and $\Phi^i$ needed in the DMD and DEIM algorithms.
%
The computational cost required to assemble and evaluate the velocity field of the flow characterizing the evolution of the coefficients reduces to the cost of the evaluation of the gradient \eqref{eq:grHDD} of the DMD-DEIM Hamiltonian at each time instant $t$. This includes:
\begin{enumerate}[label=\textnormal{(\arabic*)}]
    \item The cost to compute the linear part $U_V^{\top}V_r$ of \eqref{eq:grHDD},
    for all instances of the parameter in $\Sprmh$, is
    $O(Nn^{2}) + O(pn^{2})$.
    Note that the cost $O(Nn^2)$ required to assemble the matrix $U_V^{\top}U_V$ is performed once per stage of the RK time integrator, while the matrix-matrix product $(U_V^{\top}U_V)Z$,
    at cost $O(pn^{2})$,
    has to be performed, at every RK stage, and for each iteration of the nonlinear solver.
    \label{cost_linear}
    %
    \item The cost to compute the DMD approximation of the electric potential $\Phi^i_{\dmd}\in\r{\Nx}$, for any $i\in\Sp$,
    is $O(\Nx\Nps\T\Ndmdt)+O(\Nx\Ndmdt\Np)$, as shown in \Cref{sec:DMD}.
    This cost is linear in $\Nx$, and does not depend on the number $\Nfh$ of particles.
    The evaluation of $\Phi^i_{\dmd}$, that requires $O(\Nx\Ndmdt\Np)$ flops, needs to be performed at each stage of the RK scheme and at each iteration of the nonlinear solver. The other cost $O(\Nx\Nps\T\Ndmdt)$ is accounted for once per time step.
    %
    \item The cost to run the DEIM algorithm is
    $O(\Nfh\Nps\T d)+O(\Nfh\Nps\T c)$, as described in \Cref{sec:deim}.
    The matrix-matrix multiplication $(P_{\tau}^\top\Psi_{\tau})^{-\top}\Psi_{\tau}^{\top}M_q$ costs
    $O(\Nfh d) + O(d^3)$.
    These operations are performed once per time step. 
    \label{cost_item_1}
    \item The computation of the nonlinear time-dependent part of \eqref{eq:grHDD} for all parameters $\prm_i\in\Sprmh$, namely
    $U_X^{\top}\mathrm{diag}\big(\nabla\Lambda^0(X^i_r(t))\Phi^i_{\dmd}(t)\big)P_{\tau}$, requires $O(pdn) + O(\Np d c)$ operations. This includes the cost of the matrix-matrix multiplications
    and the cost $O(\Np d c)$ to assemble $\nabla\Lambda^0(X_r^i)$ for $d$ particles.
    \label{cost_item_4}
\end{enumerate}
To summarize, {we report, in the following table, the leading arithmetic complexity of the various steps required to solve the DMD-DEIM reduced system \eqref{eq:UZredhyp} in a fixed time interval $\Tcal_{\tau}$. Here $\ns$ denotes the number of stages of the RK time integrator and $\is$ the number of iteration of the nonlinear solver in the implicit timestepping, see also \Cref{subsec:temporal_integration} for further details.
\begin{center}
\begin{tabular}{c|c}
Leading arithmetic complexity  & Operation\\
\hline
$O(\Nfh\Nrh\Nps\ns)+O(\Nx^2\Nps\ns)+O(\Nfh\Nrh^2\ns)$
    & evolution of the reduced basis \\
$O(\Nfh\Nrh^2\ns)$
    & computation of $U_V^\top U_V$ in \eqref{eq:grHDD}\\
\hline
$O(\Nx r_{\tau}\Nps\T)$
    & construction of the DMD approximation as in \Cref{sec:DMD}\\
$O(\Nfh\Nps\T d)+ O(\Nfh\Nps\T c)$
    & construction of the DEIM approximation as in \Cref{sec:deim}\\
\hline
$O(\Np\Nrh^2\ns\is)+O(\Nx\Ndmdt\Np\ns\is)$
    & evolution of \eqref{eq:grHDD} for all parameters $\prm_i\in\Sprmh$
\end{tabular}
\end{center}
The first two rows of the table report the operations required once per
stage of the RK time integrator, while the third and forth rows refer to the operations required
once per time interval. Note that these operations are shared by all parameters, resulting in a computational cost independent of the size $\Np$ of the parameter set $\Sprmh$, but only dependent on the number of parameter subsamples $\Nps$.
The last row reports the arithmetic complexity of the parameter-dependent computations.
These need to be performed at each RK stage and iteration of the nonlinear solver, but their computational cost is independent of the number of particles $\Nfh$. Therefore,} the DMD-DEIM approximation allows a complete separation of the costs involving the number $\Nfh$ of particles and the number $\Np$ of parameters, both potentially large.

\section{Numerical temporal integration of the DMD-DEIM system}\label{subsec:temporal_integration}
For the numerical time integration of the DMD-DEIM reduced dynamics \eqref{eq:UZredhyp},
we adopt the partitioned RK method of order 2 proposed in
\cite{hesthaven2020rank}.
The idea is to combine a symplectic temporal integrator for the evolution \eqref{eq:UZhypred_coeff} of the expansion coefficient, with a time discretization of the basis evolution \eqref{eq:UZhypred_basis} able to preserve the ortho-symplectic constraint. For the latter, we adopt the tangent method proposed in \cite{P19}, summarized next.

In each temporal sub-interval $\Tcalt=(t_{\tau-1},t_{\tau}]$, given the approximate reduced basis $U_{\tau-1}$,
the method constructs a local retraction $\Rcal_{\tau-1}$ from the tangent space $T_{U_{\tau-1}}\Ucal$ into $\Ucal$ so that $U(t)=\Rcal_{\tau-1}(\V(t))$ for some $\V$ in the tangent space at $U_{\tau-1}$. For the computation of $U_{\tau}$, the idea is to evolve $\V(t)$ in the tangent space and then recover $U$ via the retraction.
The reduced problem \eqref{eq:UZredhyp} in each temporal sub-interval $\Tcalt$ is re-written in terms of the variable $Z(t)$ and $\V(t)$ as: Given $Z_{\tau-1}$ and $\V_{\tau-1}=0$, find $Z(t)$ and $\V(t)$ such that
\begin{subequations}\label{eq:PRK-ZV}
\begin{empheq}[left = \empheqlbrace\,]{align}
	&\dot{Z}(t) = \J{\Nr}\,\ghddv\big(\Rcal_{\tau-1}(\V(t)),Z(t),t\big), & \text{in }\Tcalt, \\ 
	&\dot{\V}(t) = \Ycal(\V(t),Z(t)), & \text{in }\Tcalt.
\end{empheq}
\end{subequations}
The velocity field $\Ycal:T_{U_{\tau-1}}\Ucal\times\Zcal\rightarrow \R{\Nf}{\Nr},$ describing the local flow on the tangent space at $U_{\tau-1}$, is
\begin{equation*}
    \Ycal(\V,Z) := -U_{\tau-1}(\Rcal^\top_{\tau-1}(\V)U_{\tau-1}+\Idm_{\Nr})^{-1}
    (\Rcal_{\tau-1}(\V)+U_{\tau-1})^{\top} \Upsilon(\V,Z)+ \Upsilon(\V,Z) - U_{\tau-1}\Upsilon^\top(\V,Z) U_{\tau-1},
\end{equation*}
where $\Upsilon(\V,Z)$ is given by
\begin{equation*}
\Upsilon(\V,Z):=\big(2\Xcal(\Rcal_{\tau-1}(\V),Z)
- (W U^\top_{\tau-1}-U_{\tau-1} W^\top)\Xcal(\Rcal_{\tau-1}(\V),Z)\big)
(U_{\tau-1}^{\top}\Rcal_{\tau-1}(\V)+\Idm_{\Nr})^{-1},
\end{equation*}
with $2W:=(2I_{\Nf}-U_{\tau-1}U_{\tau-1}^{\top})\V$
and $\Xcal:\R{\Nf}{\Nr}\times \Zcal\rightarrow\R{\Nf}{\Nr}$ being the velocity field of the approximate basis evolution in \eqref{eq:UZhypred_basis}, i.e.
\begin{equation*}
     \Xcal(U,Z) := (\Idm_{\Nf}-UU^\top)\big(\J{\Nf}\Ghps Z^\top{-}
	\Ghps Z^\top\J{\Nr}^\top\big)(ZZ^\top+\J{\Nr}^\top ZZ^\top\J{\Nr})^{-1}.
\end{equation*}
The retraction is defined according to \cite[Proposition 5.6]{P19} as
$\Rcal_{\tau-1}(\V)=\text{cay}(WU_{\tau-1}^\top - U_{\tau-1}W^\top)U_{\tau-1}$
where cay denotes the Cayley transform.
We refer the reader to \cite{P19,hesthaven2020rank} for further details regarding the formal derivation of \eqref{eq:PRK-ZV}. 
Note that, with the algorithm proposed in \cite[Section 5.3.1]{P19}, the computation of the retraction $\Rcal$ and the
assembly of the operator $\Ycal$ have arithmetic complexity $O(\Nfh\Nrh^2)$.

The partitioned Runge--Kutta scheme applied to \eqref{eq:PRK-ZV} reads
\begin{subequations}
\begin{empheq}[]{align}
	& Z_{\tau} = Z_{\tau-1}+\dt\sum\limits_{\ell=1}^\ns b_{\ell} k_{\ell}, \label{eq:pRK_Z}\\
	& \V_{\tau} = \dt\sum\limits_{\ell=1}^\ns \widehat{b}_{\ell} \widehat{k}_{\ell},\qquad
	U_{\tau} = \Rcal_{U_{\tau-1}}(\V_{\tau}),
	\label{eq:pRK_U}\\
	& \quad k_1 = \J{\Nr}\,\ghddvar{U_{\tau-1},
	    Z_{\tau-1} + \dt\sum\limits_{j=1}^{\ns} a_{1,j}k_j,t_{\tau-1}},
	    \qquad \widehat{k}_1 = \Xcal\bigg(U_{\tau-1},
	    Z_{\tau-1} + \dt\sum\limits_{j=1}^{\ns} a_{1,j}k_j\bigg), \label{eq:pRK_Z_k1}\\
	& \quad k_{\ell} = \J{\Nr}\,\ghddvar{\Rcal_{U_{\tau-1}}\big(\dt\sum\limits_{j=1}^{\ell-1} \widehat{a}_{\ell,j}\widehat{k}_j\big),
	    Z_{\tau-1} + \dt\sum\limits_{j=1}^{\ns} a_{\ell,j}k_j,t_{\tau-1}+c_{\ell}\dt},   \quad \ell=2,\ldots,\ns, \label{eq:pRK_Z_ki}\\
	& \quad \widehat{k}_{\ell} = \Ycal\bigg(\dt\sum\limits_{j=1}^{\ell-1} \widehat{a}_{\ell,j}\widehat{k}_j,
	    Z_{\tau-1} + \dt\sum\limits_{j=1}^{\ns} a_{\ell,j}k_j\bigg), \qquad\qquad\quad\quad\!\!\! \ell=2,\ldots,\ns, \label{eq:pRK_U_ktilde}
\end{empheq}\label{eq:discrete_dynROM}
\end{subequations}%
where $\{ a_{\ell,j},b_j,c_j \}$ and  $\{ \widehat{a}_{\ell,j},\widehat{b}_j \}$ are the set of coefficients corresponding to the implicit midpoint rule and the explicit midpoint method, respectively, \emph{cf.} \cite[Appendix A]{hesthaven2020rank}.
Note that the Hamiltonian \eqref{eq:high_fidelity_Hamiltonian} is separable, i.e.
the gradient of the electric energy determines the evolution of the state variable and the gradient of the kinetic energy describes the dynamics of the momentum.
The separability of the full order Hamiltonian is, however, not inherited by the reduced model, as seen in \eqref{eq:nonlinear_hamiltonian}. This precludes the explicit integration of \eqref{eq:pRK_Z}. Even though explicit numerical integrators for non-separable Hamiltonian, based on Hamiltonian extensions, have been recently proposed \cite{tao2016explicit}, further investigations are required to assess their accuracy in the framework of partitioned Runge--Kutta schemes.

For comparison purposes, in the numerical experiments, we will solve the full order model \eqref{eq:VPHamEqn}.
The St\"ormer-Verlet scheme \cite[Section I.1.4]{HaLuWa06} is the most popular symplectic integrator for separable Hamiltonian systems and yields the following system of equations
\begin{subequations}
    \begin{empheq}[]{align}
    & X_{\tau}^i = X_{\tau-1}^i + \Delta t \left( V_{\tau-1}^i+\frac{\Delta t}{2} E_h(X_{\tau-1}^i;\prmi) \right),\\
    & V_{\tau}^i = V_{\tau-1}^i + \frac{\Delta t}{2} \left(E_h(X_{\tau-1}^i;\prmi)+ E_h(X_{\tau}^i;\prmi)\right),
    \end{empheq}\label{eq:sv_full_order}
\end{subequations}
to be solved for each of the $p$ parameters $\eta_i\in\Gamma_h$.
Here $E_h$ denotes the approximate electric field (up to constants), i.e. $E_h(X_{\tau}^i;\prmi)=-\mmp^{-1}\mathrm{diag}(\Mq)\nabla\Lambda^0(X_{\tau}^i)L^{-1}\Lambda^0(X_{\tau}^i) \Mq$.

We summarize the main steps of the proposed method in \Cref{alg:DMD-DEIM}.

\begin{algorithm}[H]
\SetAlgoLined
\SetKwData{Left}{left}\SetKwData{This}{this}\SetKwData{Up}{up}
\SetKwFunction{Union}{Union}\SetKwFunction{FindCompress}{FindCompress}
\SetKwInOut{Input}{Input}\SetKwInOut{Output}{Output}
\caption{DMD-DEIM dynamical model order reduction}
\label{alg:DMD-DEIM}
\Input{ $f(0,x,v;\prm)$, $\{\prm_i\}_{i=1}^{\Np}$, $\Gamma_h$, $\Nps$, $\T$, $\DEIMfreq$.}
\BlankLine
Inizialize particles' positions $\{X_0^i\}_{i=1}^{\Np}$ and velocities $\{V_0^i\}_{i=1}^{\Np}$ via inverse transform sampling from $f(0,x,v;\prm_i)$. Define $W_0(\prm_i)=[X_{0}^i ; V_{0}^i]$\;
Compute $U^0$ and $Z^0$ via truncated complex SVD of $\Rcal_{0}=[W_0(\prm_1)|\dots|W_0(\prm_{\Np})]$\;
Define the parameter subset $\prmhsub\subset\Gamma_h$, of size $\Nps\ll\Np$, and the indices subset $\Spsub$ as in \Cref{sec:pSampl}\;
Compute the potential $\Phi^i_{0}$, for any $i\in \Spsub$, by solving the Poisson problem \eqref{eq:PoissonSystem}.\;
Initialize $\Ybf^i_{\dmd}:=[\Phi^i_{0}]$ and $\Ybf_{\dmd}^{\prime,i}:=[\:]$.
Initialize $\Ybf_{\deim}^i:=[\Lambda^0(X_r^{i}(t_{0}))\Phi^i_{0}]$\;
\textbf{for} \emph{$\tau=1,\dots,N_{\tau}$} \textbf{do} \\
\quad \textbf{if} \emph{$\tau\leq T$}\\
\quad\quad
Solve the system \eqref{eq:UZred_coeff}-\eqref{eq:UZred_basis}
in $\Tcalt$ with the partitioned RK described in \Cref{subsec:temporal_integration}\;
\quad\quad 
For any $i\in \Spsub$, update $\Ybf_{\dmd}^i=[\Ybf_{\dmd}^i\quad\Phi^i_{\tau}]$ and $\Ybf_{\dmd}^{\prime,i}=[\Ybf_{\dmd}^{\prime,i}\quad\Phi^i_{\tau}]$\;
\quad\quad
For any $i\in \Spsub$, update $\Ybf_{\deim}^i:=\big[\Ybf_{\deim}^i\quad \Lambda^0(X_r^{i}(t_{\tau}))\Phi^i_{\tau}\big]$.\\
\quad \textbf{elseif} $\tau>T$\\
\quad\quad At each stage of the RK time integrator:\\
\quad\qquad Solve the evolution equation \eqref{eq:UZhypred_basis} for the reduced basis as in \eqref{eq:pRK_U}\;
\quad\qquad
For any $i\in \Spsub$, update $\Ybf_{\dmd}^i=[\Ybf_{\dmd}^i\quad\Phi^i_{\tau}]$ and $\Ybf_{\dmd}^{\prime,i}=[\Ybf_{\dmd}^{\prime,i}\quad\Phi^i_{\tau}]$\;
\quad\qquad
Remove the outdated terms at $\tau-\T-1$ from $\Ybf_{\text{DMD}}^i$ and at $\tau-\T$ from $\Ybf_{\text{DMD}}^{\prime,i}$\;
\quad\qquad
For any $i\in \Spsub$, update $\Ybf_{\deim}^i:=\big[\Ybf_{\deim}^i\quad \Lambda^0(X_r^{i}(t_{\tau}))\Phi^i_{\tau}\big]$\;
\quad\qquad Remove the outdated terms at $\tau-\T-1$\;
\quad\qquad Compute $\Phi_{\dmd}$ using 
\eqref{eq:DMD_approx_potential_subset} as described in \Cref{sec:DMD}\;
\quad\qquad Approximate $\Lambda^0(X_r^{i})\Phi^i_{\dmd}$, for any $i\in\Sp$, as described in \Cref{sec:deim}\;
\quad\qquad Build $\ghddv$ as in \eqref{eq:grHDD}\;
\quad\qquad
Solve the evolution equation \eqref{eq:UZhypred_coeff} for the coefficients as in \eqref{eq:pRK_Z}.\\
\qquad\textbf{end if}\\
\textbf{end for}
\end{algorithm}


\section{Numerical experiments}\label{sec:numExp}

\subsection{Implementation and numerical study}\label{sec:implementation_details}
In this section, we apply the proposed structure-preserving dynamical model order reduction approach to several periodic electrostatic benchmark problems. 
In all the examples, computational macro-particles are loaded from a perturbed initial distribution given by
\begin{equation}\label{eq:initial_distribution}
    f(0,x,v;\prm) = f_v(v;\prm)f_x(x;\prm),
\end{equation}
where $f_v(v;\prm)$ is the initial velocity distribution, and $f_x(x;\prm):=1+\alpha\cos( kx)$ is the initial perturbation, with $k$ as the wavenumber and $\alpha$ as the amplitude of the perturbation. We have chosen physical units such that the particle mass and particle charge are normalized to one for electrons, i.e., $q=-1$ and $m=1$, and the weight $w$ of the computational macro-particles is set to $\Nfh^{-1}$. To reduce the statistical noise \cite{hockney1968characteristics,krommes2007nonequilibrium} introduced by the particle discretization \eqref{eq:fapprox} of the initial condition \eqref{eq:initial_distribution}, particles are loaded following a quiet start procedure based on a quasirandom sequence of samples. In detail, particles' positions and velocities are initialized by evaluating the inverse cumulative distribution function of $f(0,x,v;\prm)$ at the points defined by the Hammersley sequence \cite{sydora1999low} of length $\Nfh$. The distribution is defined over $\Omega:=\Omega_x\times\Omega_v$, with $\Omega_v=[-10,10]$ for all numerical experiments and $\Omega_x$ specified for each example. The quasirandom Hammersley sequence is characterized by a discrepancy value proportional to $\Nfh^{-1}$, whereas for a random distribution, the discrepancy is proportional to $\Nfh^{-1/2}$. 
Since the discrepancy measures the highest and lowest densities of points in a sequence, the Hammersley sequence guarantees that the particles are almost evenly distributed, and a significant noise reduction in the electrostatic field is therefore achieved.
{The $\Np$ test parameters are sampled from the set $\Sprm$ also using a quasirandom Hammersley sequence.}

The DMD-DEIM reduced order model \eqref{eq:UZredhyp} is numerically integrated in time, {in the interval $\Tcal=(0,t_f]$,} according to the scheme described in \Cref{subsec:temporal_integration}, resulting in the system of equations  in \eqref{eq:discrete_dynROM}.
The full order model is solved using the St\"ormer-Verlet scheme \eqref{eq:sv_full_order}.
The same {uniform time} step $\Delta t$ and number of time steps $N_{\tau}{:=t_f/\Delta t+1}$ are considered for the numerical integration of the two models.
In the following, we adopt the notation $W_{\tau}^i:=[X_{\tau}^i ; V_{\tau}^i]\in\mathbb{R}^{\Nf}$ and $R^i_{\tau}:=[X^i_{r,\tau} ; V^i_{r,\tau}]=U_{\tau}Z_{i,\tau}\in\mathbb{R}^{2\Nfh}$ to denote the numerical solutions of the discrete full order model \eqref{eq:sv_full_order} and the discrete reduced order model \eqref{eq:discrete_dynROM} for the parameter $\prm_i$ at time $t_{\tau}$, respectively.

The reducibility of the considered benchmark tests is studied in terms of the decay of the singular values of the snapshots matrices
\begin{equation}\label{eq:global_snapshots}
    S_X = \begin{bmatrix}X_{0}^1 & \cdots & X_{0}^{\Np} &
    \cdots & X_{N_{\tau}}^1 & \cdots & X_{N_{\tau}}^{\Np}\end{bmatrix}
    \quad \text{and} \quad
    S_V = \begin{bmatrix}V_{0}^1 & \cdots & V_{0}^{\Np} &
    \cdots & V_{N_{\tau}}^1 & \cdots & V_{N_{\tau}}^{\Np}\end{bmatrix},
\end{equation}
collecting the position and velocity components of $W_{\tau}^i$. As in traditional reduced basis methods, the space spanned by the selected snapshots is assumed to be representative of the solution set.
The behavior of the singular values of $S_X$ and $S_V$ is also compared to the decay of the singular values of the \emph{local} snapshots matrices 
\begin{equation}\label{eq:local_snapshots}
    S_X^{\tau} = \begin{bmatrix}X_{\tau}^1 & \cdots & X_{\tau}^{\Np}\end{bmatrix}
    \quad \text{and} \quad
    S_V^{\tau} = \begin{bmatrix}V_{\tau}^1 & \cdots & V_{\tau}^{\Np}\end{bmatrix}, \quad \forall\, \tau=1,\dots,N_{\tau},
\end{equation}
to assess the applicability and the benefits of the dynamical approach over standard global reduction methods. 
For the local snapshots matrices, we compute the ordered singular values $\{\sigma^{\tau}_{X,j}\}_j$ of $S_X^{\tau}$, for each $\tau=1,\ldots,N_{\tau}$, and normalize them with respect to the maximum singular value $\sigma^{\tau}_{X,1}$. Then, for each $j$, we consider the average over time, i.e., $\sum_{\tau} \sigma^{\tau}_{X,j}$, and the maximum over time, i.e., $\max_{\tau} \sigma^{\tau}_{X,j}$. The same study is carried out for the matrix $S_V^{\tau}$.
A further indicator of the reducibility properties of the problem is given by the numerical rank of $S_X^{\tau}$ and $S_V^{\tau}$, defined as the number of singular values larger than a user-defined threshold tolerance.
In the following, different tolerances are considered.

The accuracy of the DMD-DEIM-ROM is evaluated by computing,
for each $\tau=1,\ldots,N_{\tau}$, the relative errors
\begin{equation}\label{eq:relative_errors}
    \varepsilon_{\text{rel},X}(t_{\tau}) = \dfrac{\norm{S_{X}^{\tau} - X_{r,\tau}}_F}
    {\norm{S_{X}^{\tau}}_F},
    \quad \text{and} \quad     
    \varepsilon_{\text{rel},V}(t_{\tau}) = \dfrac{\norm{S_{V}^{\tau} - V_{r,\tau}}_F}
    {\norm{S_{V}^{\tau}}_F},
\end{equation}
where $X_{r,\tau},V_{r,\tau}\in\mathbb{R}^{\Nfh\times\Np}$ are the position and velocity components of the discrete reduced order solution $R_{\tau}=U_{\tau}Z_{\tau}\in\mathbb{R}^{\Nf\times\Np}$, respectively.
We study the error in the position and velocity of the particles separately because they are characterized by different scales of absolute error.
The relative errors \eqref{eq:relative_errors} are compared to the target values given by the projection errors
\begin{equation}\label{eq:relative_errors_complex_SVD}
    \varepsilon_{\text{rel},X}^{\text{Target}}(t_{\tau}) = \dfrac{\norm{S_{X}^{\tau} - S_{X,\text{cSVD}}^{\tau}}_F}
    {\norm{S_{X}^{\tau}}_F},
    \quad \text{and} \quad     
    \varepsilon_{\text{rel},V}^{\text{Target}}(t_{\tau}) = \dfrac{\norm{S_{V}^{\tau} - S_{V,\text{cSVD}}^{\tau}}_F}
    {\norm{S_{V}^{\tau}}_F},
\end{equation}
where $S_{X,\text{cSVD}}^{\tau}, S_{V,\text{cSVD}}^{\tau}\in\mathbb{R}^{\Nfh\times\Np}$ are the position and velocity components of the projection of the snapshots onto the space spanned by the ortho-symplectic basis of size $\Nf\times \Nr$ obtained from the Complex SVD \cite{peng2016symplectic} 
of the matrix $S_X^{\tau}+\imath S_V^{\tau}$.
Since the Complex SVD provides the ortho-symplectic basis that minimizes the projection error in the snapshots \cite[Theorem 4.6]{peng2016symplectic},
comparing \eqref{eq:relative_errors} and \eqref{eq:relative_errors_complex_SVD} allows to test the approximability properties of the reduced basis constructed in the dynamical approach.

Moreover, we analyze the evolution of the electric field energy \eqref{eq:HamEEfull} for the reduced order approximation, i.e., $\EE(X_{r,\tau}^{i};\eta_i)$, and for the full order solution, i.e., $\EE(X_{\tau}^{i};\eta_i)$, for each instance $\prmi$ of the parameter.
This term gives information of the macroscopic behavior of the plasma, and it is also the one affected by more levels of approximation.

Finally, the efficiency of the proposed approach is investigated by comparing the running times required for the integration over a single time step of the fully-discrete DMD-DEIM reduced model \eqref{eq:discrete_dynROM} and of the discrete full order model \eqref{eq:sv_full_order}. The running time for the full model is obtained by summing the times required for each instance of the parameter $\prm_i\in\Sprmh$. The comparison focuses on the scalability in the approximation of parametric problems as the size $\Np$ of the parameter set increases, a typical scenario in a multi-query context. 
To compare the efficiency of the different methods, we analyze the runtime  required for integration over a single time interval for all parameter values considered. 
The values reported were obtained as the average of the runtimes obtained in the first $25$ time intervals.
For the dynamical reduced basis method, we also analyze separately the contributions to the computational cost due to the basis evolution \eqref{eq:pRK_U}-\eqref{eq:pRK_U_ktilde} and  to the coefficients evolution \eqref{eq:pRK_Z}-\eqref{eq:pRK_Z_k1}-\eqref{eq:pRK_Z_ki}, in line with the theoretical findings of \Cref{sec:DMD_DEIM_complexity}.

In the construction of the DMD-DEIM reduced model, we consider a tolerance equal to $10^{-5}$ in the computation of the DMD modes in \Cref{alg:DMD_algorithm}. Moreover, the nonlinear system, \eqref{eq:pRK_Z_k1} and \eqref{eq:pRK_Z_ki}, describing the evolution of the increments $k_{\ell}$, $\ell=1,\dots, \ns$, is solved using the fixed point iteration method. As a stopping criterion for the nonlinear solver, we check when the relative norm of the update to $k_{\ell}$ is smaller than the threshold value $10^{-9}$.
All numerical simulations are performed using Matlab  on computer nodes with Intel Xeon
E5-2643 (3.40GHz). The code and the data supporting the findings of this study are available upon request.


\subsection{Weak Landau damping of 1D Langmuir waves}
\begin{figure}[H]
    \centering
    \includegraphics{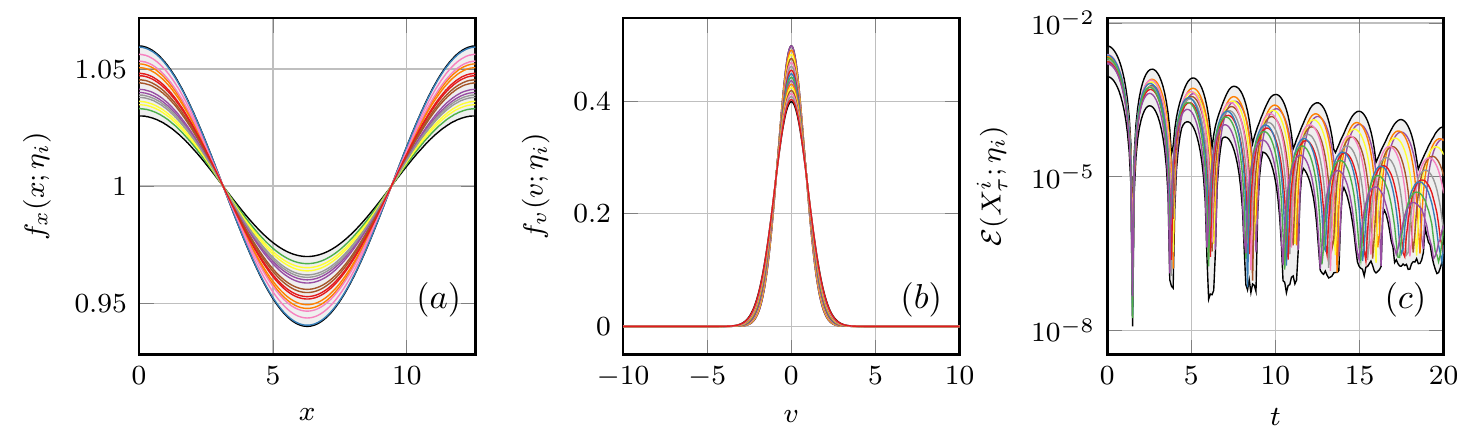}
    \caption{LD: $(a)-(b)$ Initial position and velocity distributions for selected values of the parameter in $\Sprmh$. $(c)$ Exponential time decay of the electrostatic energy $\EE(X^i_{\tau};\prmi)$ obtained from the full model solution, for selected values of $\prmi$ in $\Sprmh$.
    Since not all parameters in $\Sprmh$ are reported, the black lines in each subplot are used to mark the region where the plotted quantity is contained, for any value of the parameter in $\Sprmh$.}
    \label{fig:parametric_initial_condition_energy_LD}
\end{figure}

The first application we consider is the study of the damped propagation of small amplitude plasma waves, also known as Landau damping (LD). The resonance between physical particles and the propagating wave generates damping of the electric field energy, without particle collisions. This process is used in particle accelerators to prevent coherent beams oscillations that could cause potential instabilities \cite{hereward1977landau}. The initial 
condition is given by \eqref{eq:initial_distribution}
with the velocity distribution function
\begin{equation}\label{eq:LD_initial_velocity_distribution}
    f_v(v;\prm) = \dfrac{1}{\sqrt{2\pi} \sigma}\exp\left (-\frac{v^2}{2\sigma^2}\right ),
\end{equation}
where the amplitude of the perturbation $\alpha$ and the standard deviation $\sigma$ of the velocity Maxwellian are the study parameters $\prm=(\alpha,\sigma)$, with $\prm\in\Sprm=[0.03,0.06]\times [0.8,1]$, and the perturbation wavenumber $k$ is fixed to $0.5$. Following the sampling procedure described in Section \ref{sec:implementation_details}, we solve the Landau damping problem for $\Np=300$ different realizations of the parameter $\prm$. In Figure \ref{fig:parametric_initial_condition_energy_LD}(a) and \ref{fig:parametric_initial_condition_energy_LD}(b) the initial position and velocity distributions are shown for several of the selected parameter values. 
We consider periodic boundary conditions on the physical space domain $\Omega_x:=\left(0,\frac{2\pi}{k}\right)$ with a uniform neutralizing background charge. For the numerical solution of the full order model, we use $\Nx=32$ piecewise linear basis for the Poisson solver, and $\Nfh=5\times 10^{4}$  macro-particles for the approximation of the solution density.
A uniform time step $\Delta t=0.0025$ has been adopted for the evolution of particles' positions and velocities over the time interval $\Tcal=(0,20]$.

\begin{figure}[H]
    \centering
    \includegraphics{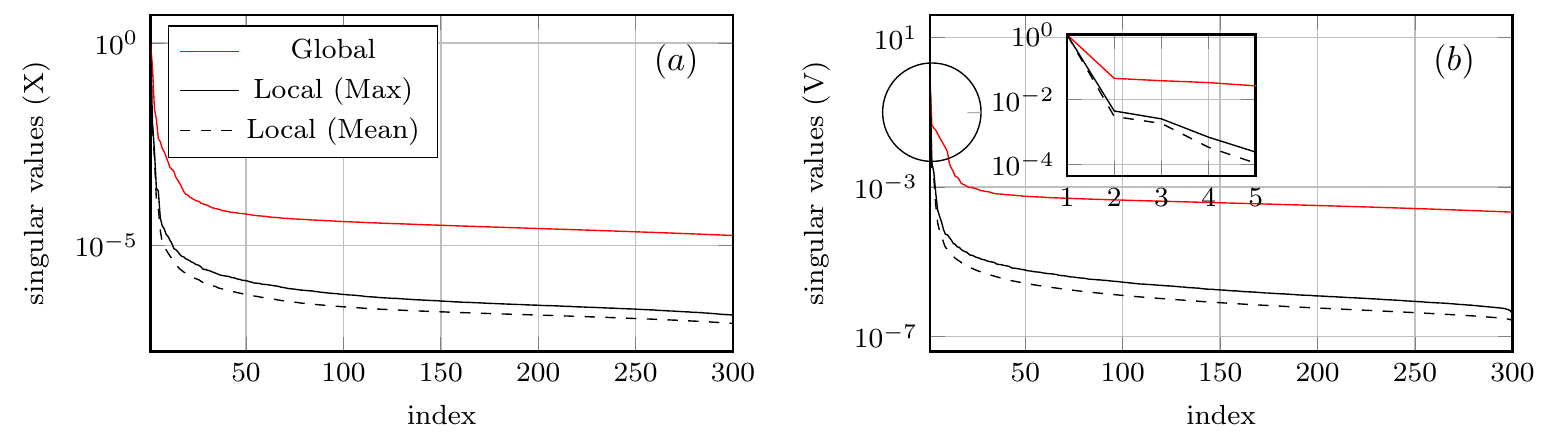}
    \caption{LD: Singular values of the global snapshots matrices $S_X$ and $S_V$ compared to the maximum and time average (in $\tau$) of the singular values of the local matrices $S_X^{\tau}$ and $S_V^{\tau}$.}
    \label{fig:singular_value_decay_LD}
\end{figure}

In Figure \ref{fig:singular_value_decay_LD}, the decay of the singular values of the global snapshots matrices $S_X$ and $S_V$, normalized with respect to the corresponding largest singular value, are compared to the maximum and averages over $\tau$ of the normalized singular values of the local counterparts $S_X^{\tau}$ and $S_V^{\tau}$, computed as described in \Cref{sec:implementation_details}. Concerning the particles position, although a plateau of the singular values can be seen for both global and local matrix, the initial decay is sharper in the local case with singular values that are two orders of magnitude smaller than in the global case, suggesting a more efficient representation using a local low-rank model. This gap increases when considering the particles velocity, suggesting that a \emph{global} reduced basis approach would not be effective in reducing the computational cost of the Landau damping simulation. 

\begin{figure}[H]
    \centering
    \includegraphics{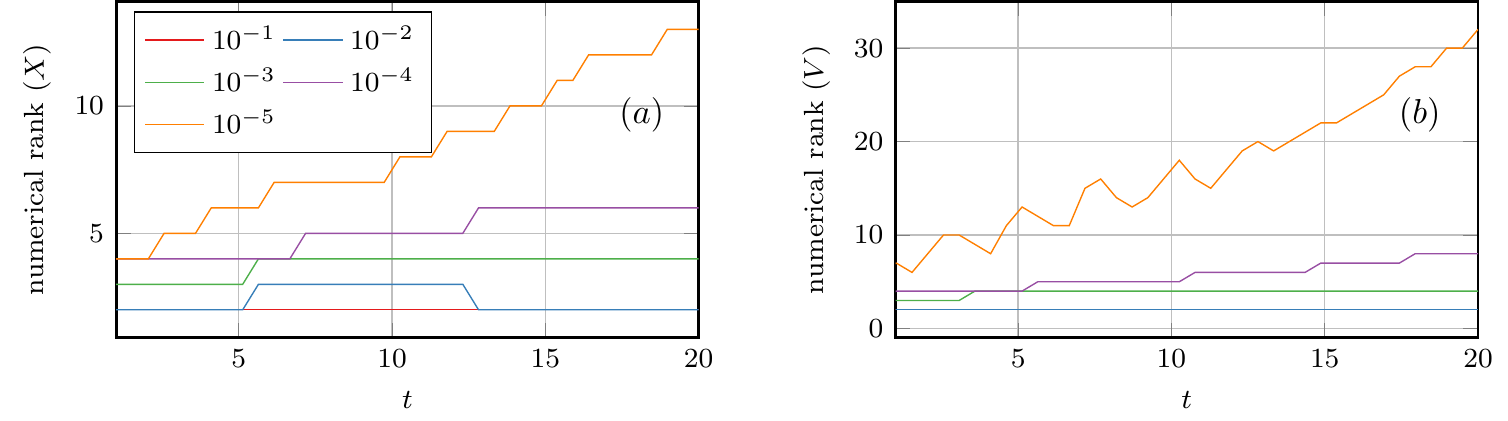}
    \caption{LD: Numerical rank of $S_X^{\tau}$ in (a) and $S_V^{\tau}$ in (b), as a function of $\tau$. Different colors are associated with different values of the threshold, according to the legend.}
    \label{fig:local_rank_position_velocity_LD}
\end{figure}

\begin{figure}[H]
    \centering
    \includegraphics{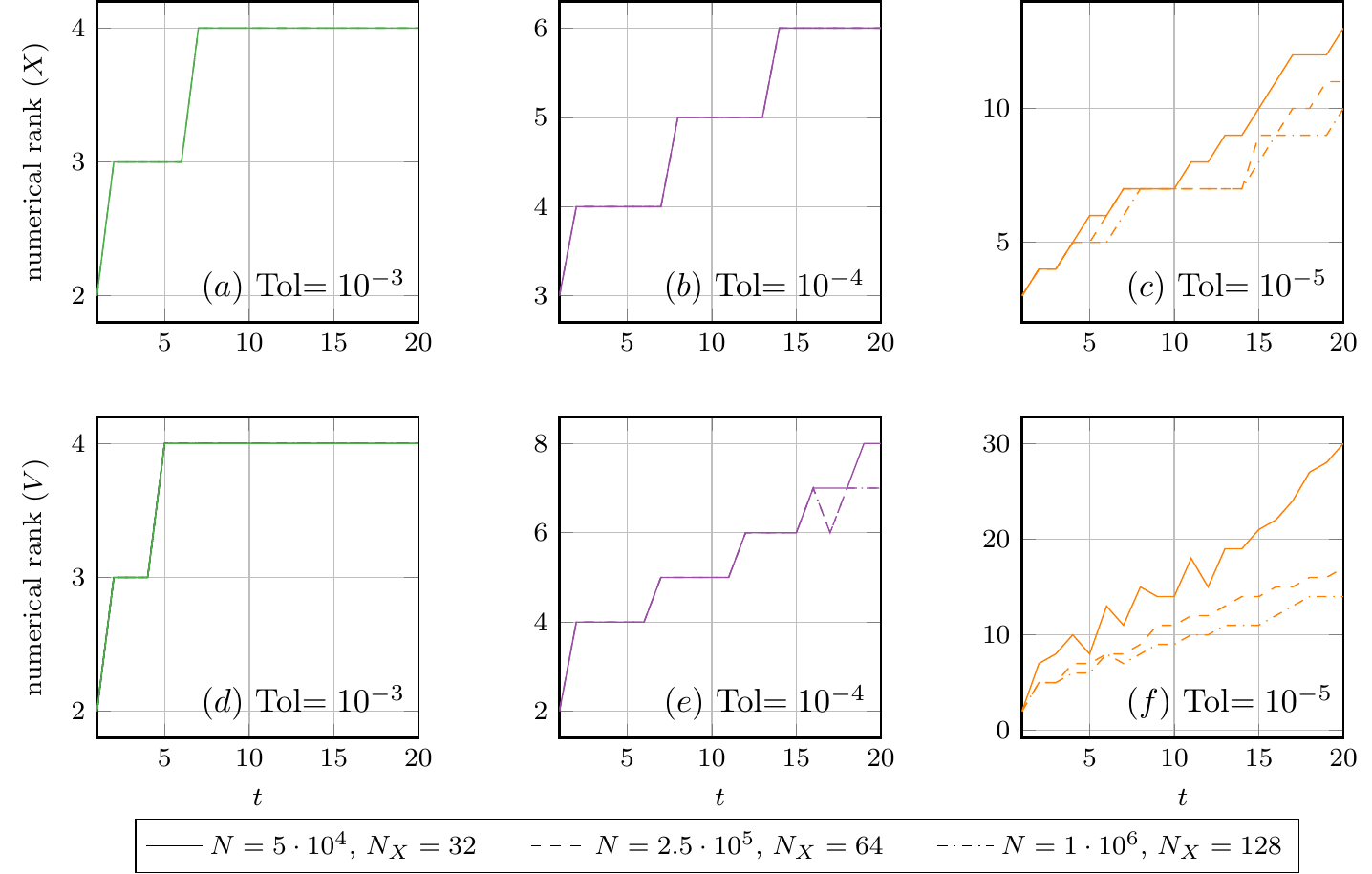}
    \caption{LD: Evolution of the numerical ranks of $S_X^{\tau}$ $(a)-(c)$ and of $S_V^{\tau}$ $(d)-(f)$ for different threshold indicated by Tol. In each subfigure, the rank behavior for different values of $\Nx$ and $\Np$, in the discretization of the full order model, is compared.
    }
    \label{fig:local_rank_position_velocity_LD_comparison}
\end{figure}

In Figure \ref{fig:local_rank_position_velocity_LD}, we report the numerical rank of the matrices $S_X^{\tau}$ and $S_V^{\tau}$ as a function of $\tau$ and for different values of the threshold. The numerical rank remains constant and below $4$ for tolerances larger than $10^{-3}$, grows to $8$ for a tolerance of $10^{-4}$, and reaches a maximum of $13$, for positions and $32$ for velocities, when the tolerance is set to $10^{-5}$.
The increase in the solution complexity over time is partially due to the accumulation of statistical noise associated with the discretization of $f(t,x,v;\prm)$ by macro-particles. In Figure \ref{fig:local_rank_position_velocity_LD_comparison}, in support of this conclusion, we note that as the average number of particles per cell increases during the initial particle loading phase, the numerical ranks of $S_X^{\tau}$ and $S_V^{\tau}$, at fixed tolerance, decrease.

This behavior of the numerical rank suggests that evolving the basis, but keeping the rank of the approximation constant, is sufficient to accurately approximate the solution of the full order model, at least in this test case.

Concerning the reduced dynamical model, we consider $\Nr=4$ as the reduced manifold dimension. For the DEIM reduction described in \Cref{sec:deim}, $d=32$ interpolation points have been used to reduce the approximation error, and $\nDEIMup=12$ DEIM indices are updated at each time step, for the sake of efficiency, according to \eqref{eq:indices_adaptive_DEIM}. All DEIM indices are recomputed every $\DEIMfreq=3$ time steps {using \cite[Algorithm 1]{chaturantabut2010nonlinear}. We refer to \Cref{sec:deim} for further details.}

\begin{figure}[H]
    \centering
    \includegraphics{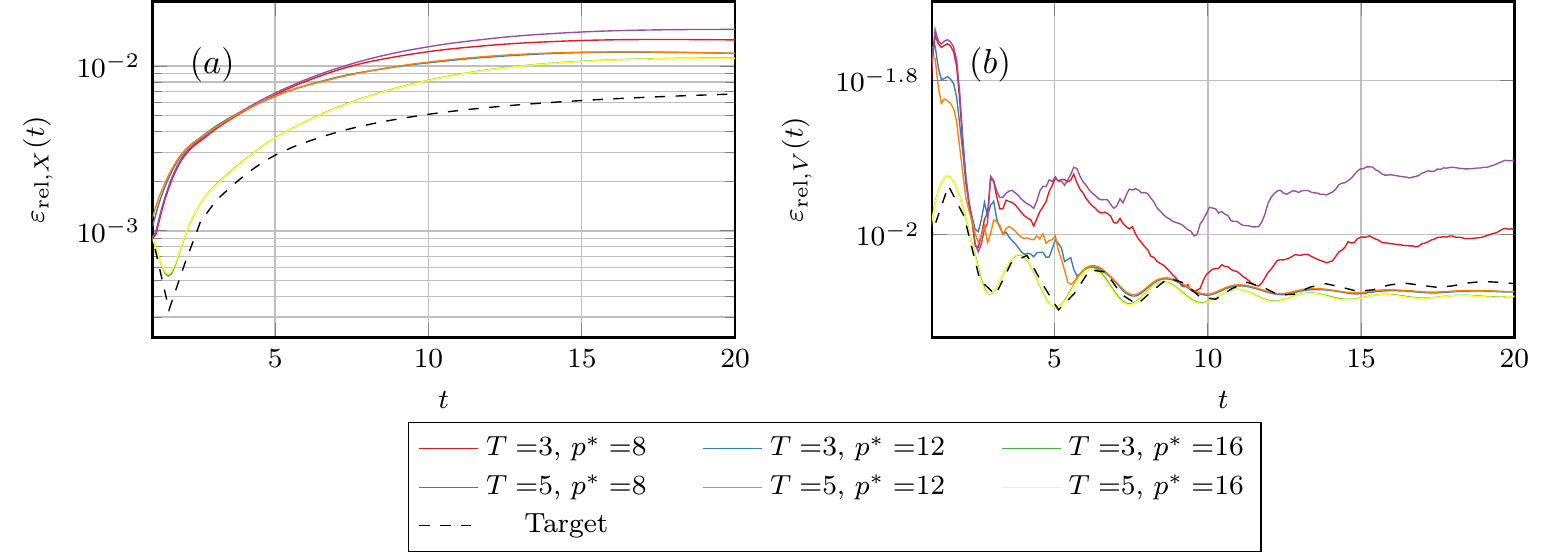}
    \caption{LD: Evolution of the position $(a)$ and velocity $(b)$ relative errors, as defined in \eqref{eq:relative_errors}, for different choices of $\Nps$ and $\T$. These errors are compared to the target values given by the position component $\varepsilon_{\text{rel},X}^{\text{Target}}$ and the velocity component $\varepsilon_{\text{rel},V}^{\text{Target}}$ of the relative projection errors defined in \eqref{eq:relative_errors_complex_SVD}.
    The target reduced basis has dimension $4$ and is computed, for each time step, using the Complex SVD algorithm, as described in \Cref{sec:implementation_details}.}
    \label{fig:time_error_position_velocity_LD}
\end{figure}

For this test case, we include a numerical study of the evolution of the approximation errors $\varepsilon_{\text{rel},X}$ and $\varepsilon_{\text{rel},V}$ in \eqref{eq:relative_errors} under variations of the size $\Nps$ of the subset of the parameters used to evolve the basis effectively, according to Section \ref{sec:pSampl}, and the length {$(\T+1)\Delta t$} of the time window adopted to harvest the self-consistent electric potential $\Phi(X_r^i)$ for the DMD extrapolation step, as described in \Cref{sec:DMD}. In particular, we consider $\Nps\in\{8,12,16\}$ and $\T\in\{3,5\}$ and the results are shown in Figure \ref{fig:time_error_position_velocity_LD}. In all tested combinations of $T$ and $\Nps$, the error is proportional to the best approximation error, both in position and velocity.
As $\Nps$ increases, both errors decrease: this is expected since a more refined sampling of the subset $\prmhsub$ results in a more accurate representation of the dynamics of evolution of the bases in \eqref{eq:UZhypred_basis}. We also note that, for $\Nps=16$, the error $\varepsilon_{\text{rel},V}$ is, at several time instances, smaller than the target value. This performance can be explained by the fact that the optimality of the Complex SVD algorithm concerns the projection of the entire state $[S_X^{\tau}\, S_V^{\tau}]$ and not of its components individually.
We observe that the {DMD number of samples, $\T+1$,} has no impact on the error when $\Nps$ is large, and a small accuracy degradation is even registered for $\Nps=8$ when $\T=5$ is chosen over $\T=3$. The optimal choice of $\T$ remains an open problem: as pointed out in \cite{dylewsky2019dynamic}, it should capture slow and fast scales of the local dynamics, but a rigorous optimization strategy would require a study of the multi-scale properties of the solution to the Vlasov--Poisson equation for each of the parameter realization considered. However, we stress that the results are relatively robust concerning this parameter.

\begin{figure}[H]
    \centering
    \includegraphics{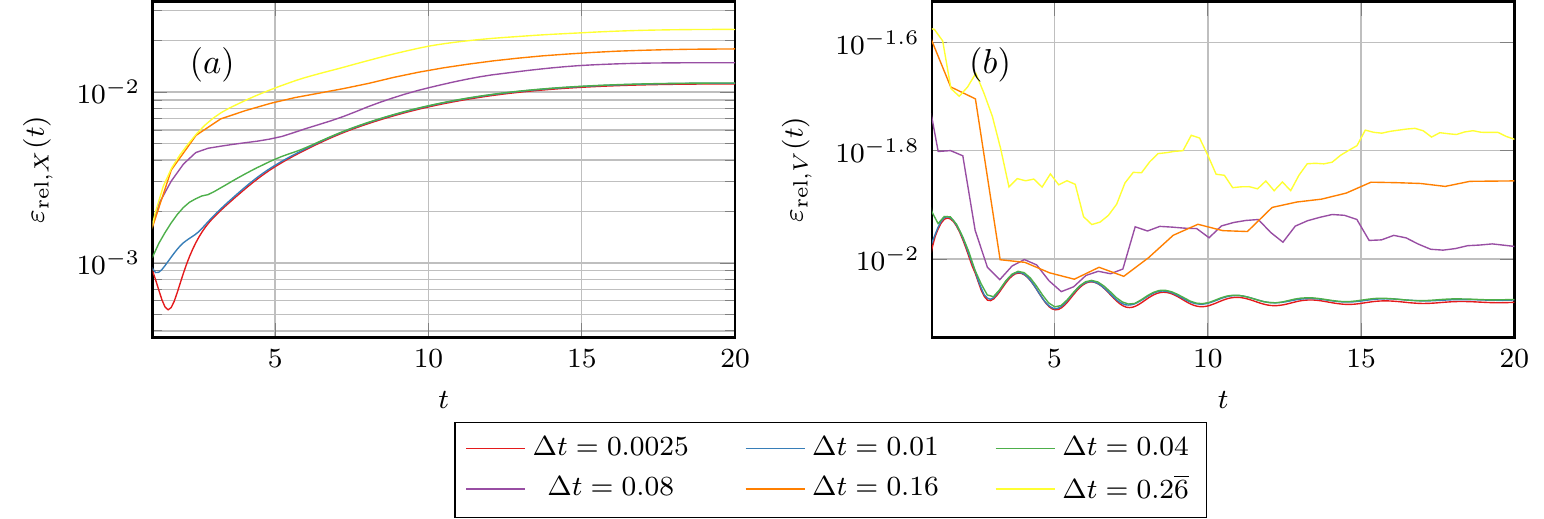}
    \caption{{LD:  Evolution of the position $(a)$ and velocity $(b)$ relative errors, as defined in \eqref{eq:relative_errors}, for different values of the time step $\Delta t$. The hyperparameters are $\Nps=12$, $\T=3$, and $\nDEIMup=1$.}
    }
    \label{fig:error_time_step}
\end{figure}

We remark that in this test case the time step $\Delta t$ is rather small, $\Delta t=0.0025$. The rationale is that a small time step allows us to gauge the error introduced by the reduction and hyper-reduction without pollution from other sources of errors, such as the temporal discretization. However, larger time steps can be chosen without a significant decrease of accuracy, as shown in \Cref{fig:error_time_step}. There we report the evolution of the relative errors $\varepsilon_{\text{rel},X}$ and $\varepsilon_{\text{rel},V}$ for values of the time step $\Delta t \in\{0.0025, 0.01, 0.04, 0.08, 0.16, 0.2\overline{6}\}$
in the case $T=3$, $p^*=12$, and $k_{DEIM}=1$. 
It can be observed that the timestepping is not the dominant source of error up to the value $\Delta t=0.04$, and that, for $\Delta t > 0.04$, increasing the time step does not significantly affect accuracy.

\begin{figure}[H]
    \centering
    \includegraphics{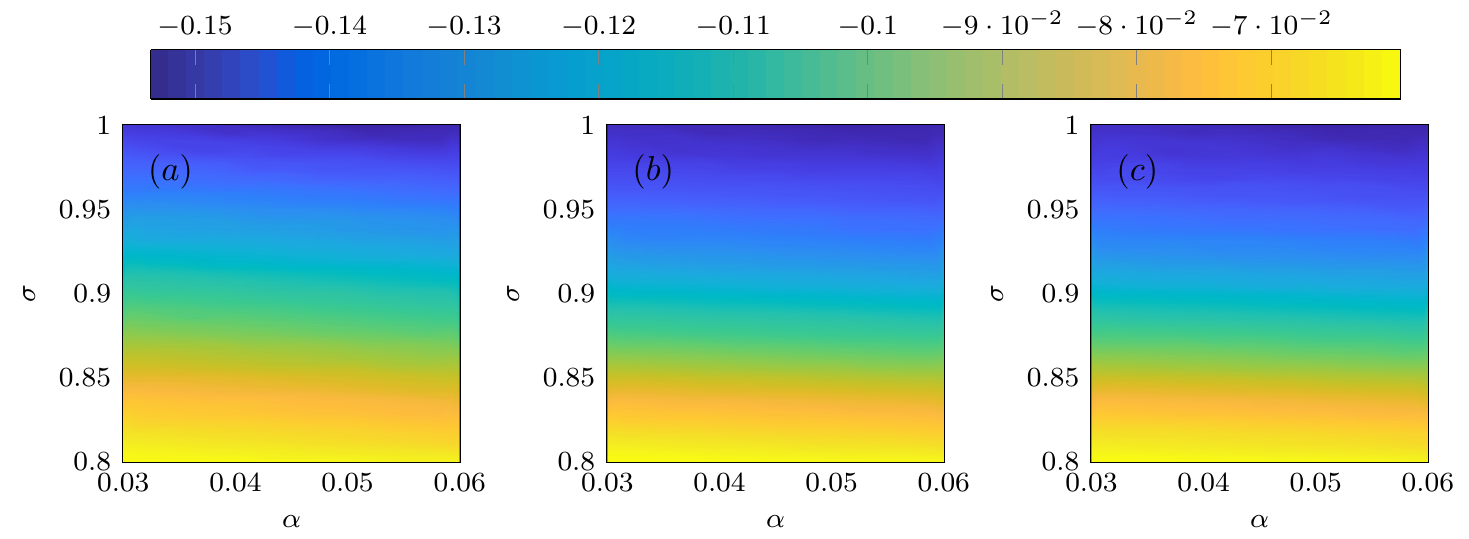}
    \caption{LD: Damping rates of the exponential time decay of
    the electric energy $\EE(X^i;\prm)$, defined in \eqref{eq:HamEEfull}, as a function of the two-dimensional parameter $\prm=(\alpha,\sigma)$. The plots refer to $(a)$ the full order model; $(b)$ the dynamical reduced model with $\T=3$ and $\Nps=16$; and $(c)$ the dynamical reduced model with $\T=5$ and $\Nps=8$.
    }
    \label{fig:image_energy_damping_LD}
\end{figure}

Landau theory \cite{berge1969landau} 
establishes that, for small perturbations of the initial analytical data of the form \eqref{eq:LD_initial_velocity_distribution}, the electric energy $\EE(X^i;\prm_i)$ decays (in time) exponentially with a damping factor that depends on the standard deviation $\sigma$ of the Maxwellian distribution $f_v$ but is independent of the amplitude of the perturbation $\alpha$. 
{In \Cref{fig:image_energy_damping_LD} we report the damping rate of the electric energy for each value of the parameter $\prm=(\alpha,\sigma)\in\Sprm$:
for the $\Np$ test parameters the damping rate is computed as the slope of the peaks of the electric energy (excluding the first one); all other values are generated via linear interpolation/extrapolation on a uniform grid of $90000$ points.}
As can be seen from \Cref{fig:image_energy_damping_LD}, the {theoretical} dependence of the damping rate on the considered parameters is captured by the reduced model numerical solution.

\begin{figure}[H]
    \centering
    \includegraphics{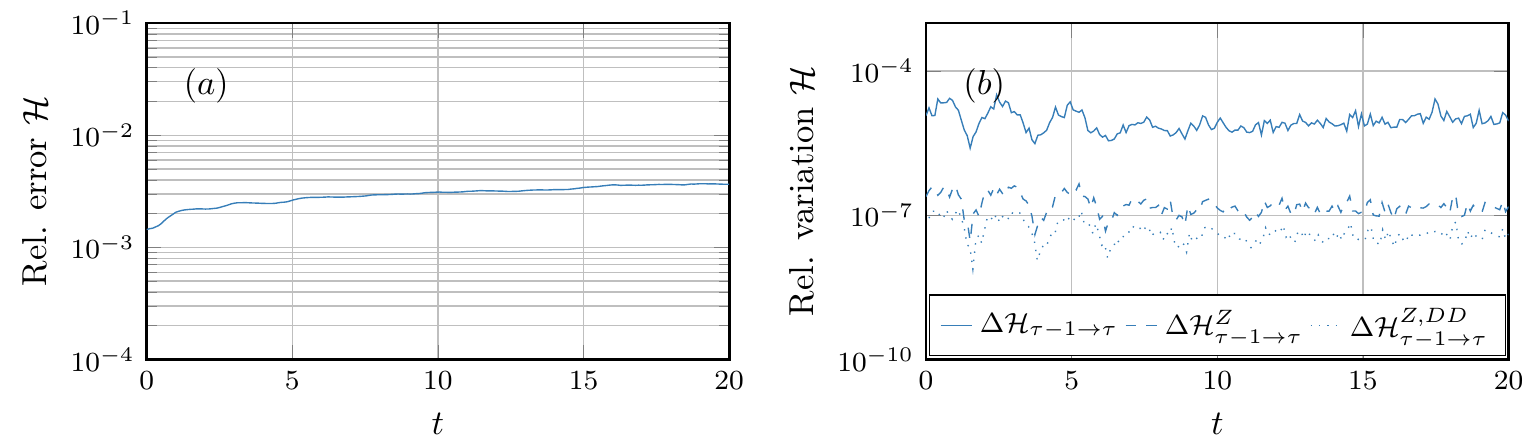}
    \caption{LD: $(a)$ Evolution of the relative error \eqref{eq:rel_error_Hamiltonian} of the Hamiltonian. 
    $(b)$ Evolution of the components $\Delta\Ham_{\tau-1\rightarrow\tau}$, $\Delta\Ham_{\tau-1\rightarrow\tau}^{Z}$ and $\Delta\Ham_{\tau-1\rightarrow\tau}^{Z,\dd}$ of the error bound \eqref{eq:error_ham_tot}, \eqref{eq:error_ham_tot1} in the local conservation of the reduced Hamiltonian.
    The values of the hyper-parameters are set to $\Nps=12$ and $\T=3$, respectively.}
    \label{fig:hamiltonian_comparison_LD}
\end{figure}

In Figure \ref{fig:hamiltonian_comparison_LD}(a), we report the evolution of the relative error of the Hamiltonian \eqref{eq:high_fidelity_Hamiltonian} computed in the reduced and full model solutions, i.e.
\begin{equation}\label{eq:rel_error_Hamiltonian}
    \dfrac{\norm{\Ham([X_{\tau},V_{\tau}])-\Ham(U_{\tau} Z_{\tau})}_2}{\norm{\Ham([X_{\tau},V_{\tau}])}_2},\qquad \forall\,\tau=1,\ldots,N_{\tau}.
\end{equation}
It is observed that the error is bounded and grows only slowly over time. The reason why the Hamiltonian is not exactly preserved is twofold: the numerical temporal integrator is symplectic but not Hamiltonian-preserving, and the reduced model possesses a Hamiltonian structure but with an approximate Hamiltonian function.
{As pointed out in \Cref{sec:hyperred}, applying DEIM directly to the nonlinear Hamiltonian vector field yields an approximation that, in general, does not possess a gradient structure and hence fails to provide a Hamiltonian dynamics.}
To better understand {the aforementioned} two sources of errors, we consider the error in the Hamiltonian at two consecutive time instances of the solution.
For the DMD-DEIM reduced model discretized with the partitioned Runge--Kutta method
described in Section \ref{subsec:temporal_integration}, it holds
\begin{subequations}
    \begin{empheq}[]{align}
    \Delta\Ham_{\tau-1\rightarrow\tau} := \norm{\Ham(U_{\tau}Z_{\tau})-\Ham(U_{\tau-1}Z_{\tau-1})}_2
    &\leq \norm{\Ham(U_{\tau}Z_{\tau})-\Ham(U_{\tau-\frac{1}{2}}Z_{\tau})}_2 \label{eq:error_ham_first_basis}\\
    &+ \norm{\Ham(U_{\tau-\frac{1}{2}}Z_{\tau-1})-\Ham(U_{\tau-1}Z_{\tau-1})}_2 \label{eq:error_ham_second_basis}\\ 
    &+ \norm{\Ham(U_{\tau-\frac{1}{2}}Z_{\tau})-\Ham(U_{\tau-\frac{1}{2}}Z_{\tau-1})}_2 \label{eq:error_ham_coeff}.
    \end{empheq}\label{eq:error_ham_tot}
\end{subequations}
The first two terms \eqref{eq:error_ham_first_basis} and \eqref{eq:error_ham_second_basis} depend on the numerical time integration of the basis equation \eqref{eq:pRK_U}, while the last term \eqref{eq:error_ham_coeff}
also depends on the DMD-DEIM approximation of the Hamiltonian. In particular, it holds
\begin{subequations}
    \begin{empheq}[]{align}
    \Delta\Ham_{\tau-1\rightarrow\tau}^{Z}:=
    \norm{\Ham(U_{\tau-\frac{1}{2}}Z_{\tau})-\Ham(U_{\tau-\frac{1}{2}}Z_{\tau-1})}_2
    &\leq \norm{\Ham(U_{\tau-\frac{1}{2}}Z_{\tau})-\Hamgen{\tau-\frac{1}{2}}^{\dd}(Z_{\tau},t_{\tau})}_2 \label{eq:error_ham_DMDDEIM_first}\\
    &+\norm{\Hamgen{\tau-\frac{1}{2}}^{\dd}(Z_{\tau-1},t_{\tau-1})-\Ham(U_{\tau-\frac{1}{2}}Z_{\tau-1})}_2, \label{eq:error_ham_DMDDEIM_second}\\
    &+\norm{\Hamgen{\tau-\frac{1}{2}}^{\dd}(Z_{\tau},t_{\tau})-\Hamgen{\tau-\frac{1}{2}}^{\dd}(Z_{\tau-1},t_{\tau-1})}_2, \label{eq:error_ham_int} 
    \end{empheq}\label{eq:error_ham_tot1} 
\end{subequations}
where $\Ham^{\dd}_U$ is defined in \eqref{eq:Hdd}.
The first two terms \eqref{eq:error_ham_DMDDEIM_first} and \eqref{eq:error_ham_DMDDEIM_second}
depend on the approximation of the Hamiltonian introduced by the DMD-DEIM method, while the last term \eqref{eq:error_ham_int}, that we dub $\Delta\Ham_{\tau-1\rightarrow\tau}^{Z,\dd}$, in only associated with the numerical time integrator of the coefficient equation.
In Figure \ref{fig:hamiltonian_comparison_LD}(b) we report the time evolution of $\Delta\Ham_{\tau-1\rightarrow\tau}$, $\Delta\Ham_{\tau-1\rightarrow\tau}^{Z}$ and $\Delta\Ham_{\tau-1\rightarrow\tau}^{Z,\dd}$, for the hyper-parameters $\T=3$ and $\Nps=12$. We can observe that the DMD-DEIM method provides a good approximation of the Hamiltonian (dashed line).

\begin{figure}[H]
    \centering
    \includegraphics{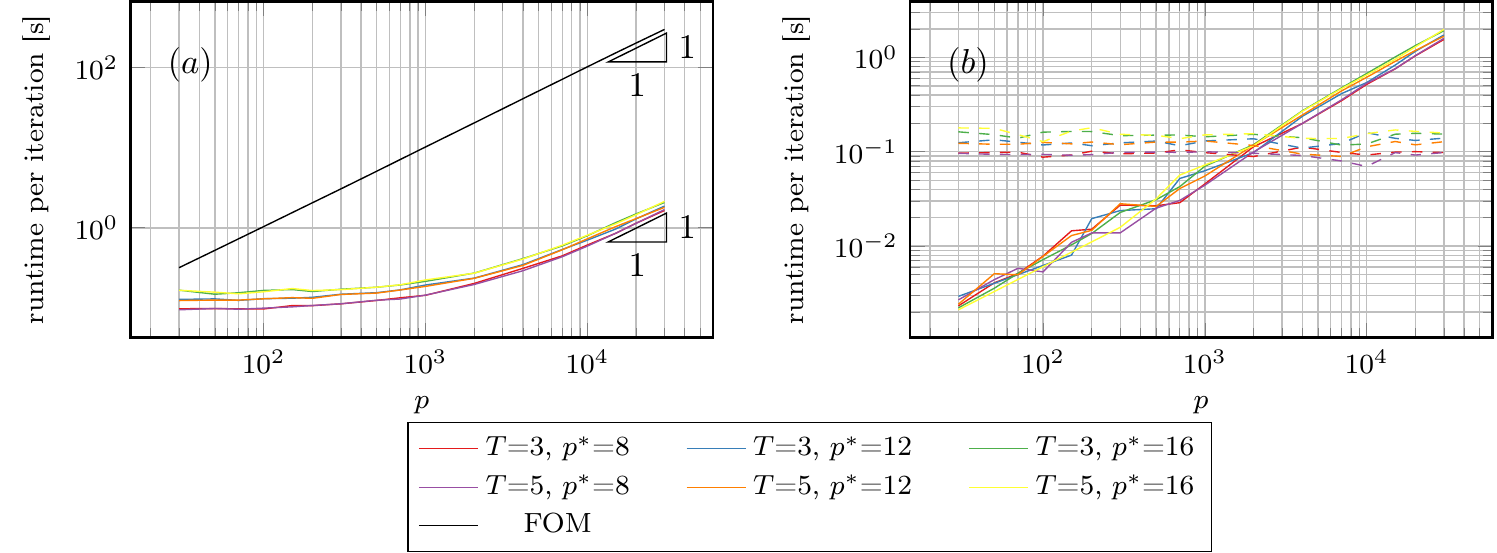}
    \caption{LD: $(a)$ Comparison of the runtime {per time step} (in seconds) between the full order solver and the dynamical reduced basis approach for different hyper-parameter configurations, as function of the parameter sample size $\Np$. $(b)$ Separation of contributions to the {runtime per time step} of the reduced model due to basis evolution \eqref{eq:pRK_U} (dashed lines) and coefficients evolution \eqref{eq:pRK_Z} (continuous line).}
    \label{fig:running_time_LD}
\end{figure}

To study the algorithm efficiency, we investigate the runtime as a function of the number $\Np$ of tested parameters.
The proposed approach outperforms the full order solver, as shown in Figure \ref{fig:running_time_LD}(a), and the gap widens as the value of $\Np$ increases.
Depending on the choice of hyper-parameters of the reduced model, the algorithm speed-up varies between $1.9$ and $3.3$ when $\Np=30$ and between $46$ and $71$ when $\Np=1000$.
For $\Np\geq 2000$, the evolution of the expansion coefficients \eqref{eq:UZhypred_coeff} becomes computationally more demanding than the evolution of the reduced basis \eqref{eq:UZhypred_basis}, as shown in Figure \ref{fig:running_time_LD}(b), and the overall computational cost of the reduced model begins to grow approximately linearly as a function of $\Np$. Thus, for values of $\Np$ larger than $2000$, the ratio between the time required to integrate the full model and the time to integrate the reduced model remains constant, with speed-ups ranging between 141 and 183, depending on the values of the hyper-parameters.
We also remark that the computational cost to evolve the reduced basis (dashed lines in \Cref{fig:running_time_LD}(b)) is independent of $\Np$ because it only depends on the number $\Nps$ of subsamples.

\begin{figure}[H]
    \centering    \includegraphics{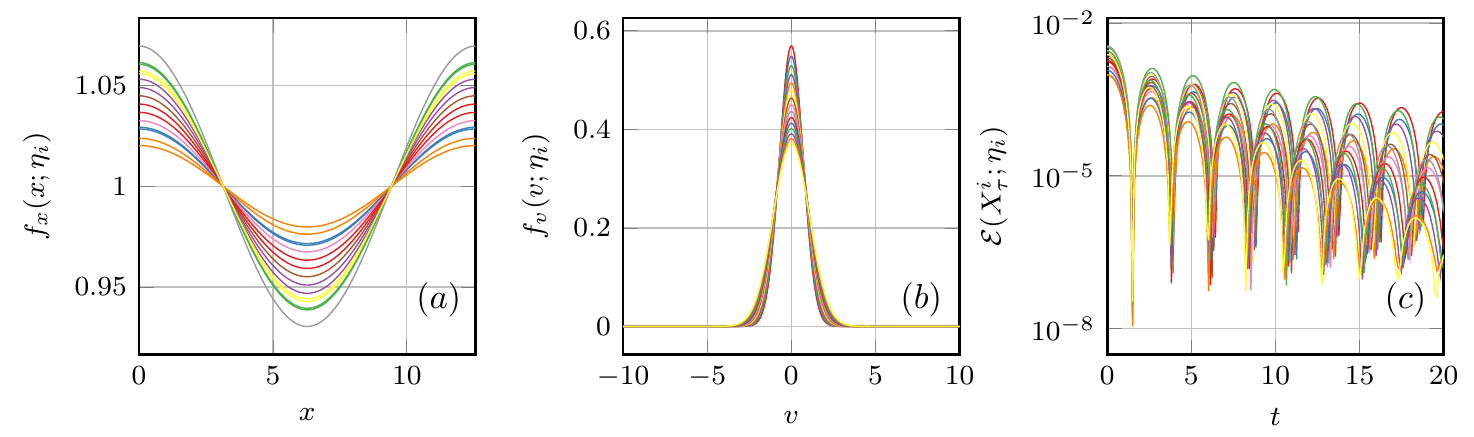}
    \caption{LD: $(a)-(b)$ Initial position and velocity distributions for selected values of the parameter $\prmi$ in $\Sprm_{2,h}$. $(c)$ Exponential time decay of the electrostatic energy $\EE(X^i_{\tau};\prmi)$ obtained from the full model solution, for selected values of $\prmi$ in $\Sprm_{2,h}$.}\label{fig:initial_condition_extended_parameter}
\end{figure}

\begin{figure}[H]
    \centering
    \includegraphics{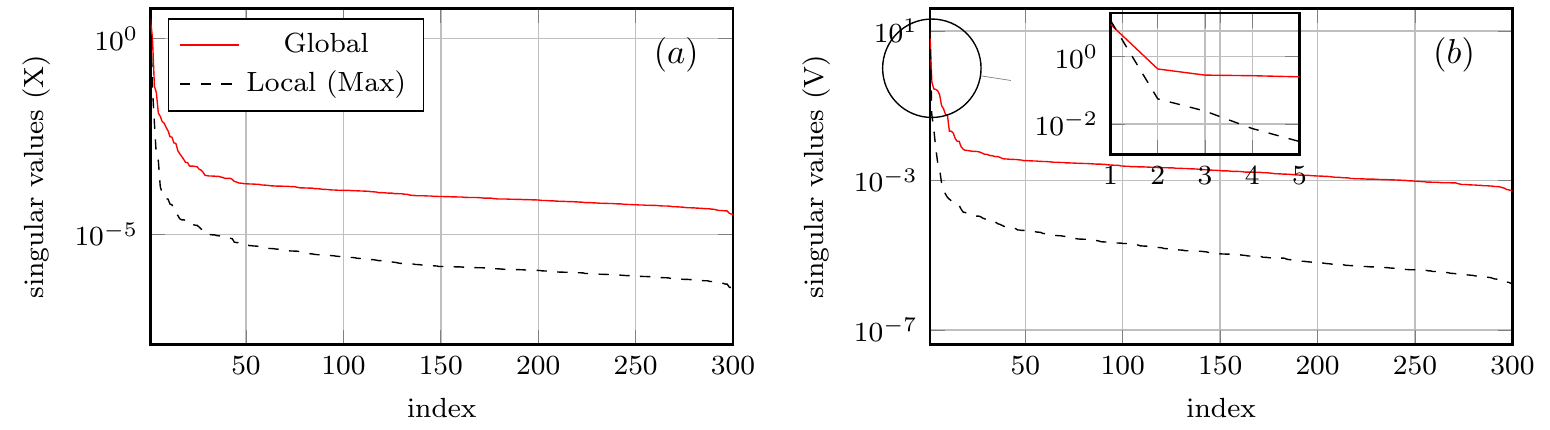}
    \caption{LD: Singular values of the global snapshots matrices $S_X$ and $S_V$ compared to the singular values of the local matrices $S_X^{\tau}$ and $S_V^{\tau}$ for $\Sprm=\Sprm_2$.}\label{fig:singular_value_extended_parameter}
\end{figure}

The choice and the number $\Nps$ of subsamples has also an effect on the accuracy of the reduced order approximation. In the following test, we numerically study the
convergence of the algorithm for different numbers of sampling points and for different parameter domains.
In particular, we compare the performances of the algorithm when the parameter $\prm=(\alpha,\sigma)$ belongs to the domain $\Sprm\in\{\Sprm_1,\Sprm_2\}$, where $\Sprm_2=[0.02,0.07]\times [0.7,1.1]$ and $\Sprm_1=[0.03,0.06]\times [0.8,1]$ is the domain considered above.
Although the length of the 1D parameter intervals used to define $\Sprm_1$ and $\Sprm_2$ might not appear very different, what matters for the numerical approximation is the reducibility of the solution as a function of the parameter, that is the approximability of the parameter-to-solution map.
To gauge the reducibility of the solution under variation of the parameter we perform a qualitative (\Cref{fig:initial_condition_extended_parameter}) and a quantitative assessment (\Cref{fig:singular_value_extended_parameter}).
More in details, analogously to Section 6.2, we consider $\Np=300$ parameter samples in the domain $\Sprm$.
In \Cref{fig:initial_condition_extended_parameter} we report the initial distribution of position (plot (a)) and velocity (plot (b)) and the time evolution of the electrostatic energy (plot (c)) for different samples of the parameters in $\Sprm_2$. Comparing with \Cref{fig:singular_value_decay_LD}, where $\Sprm=\Sprm_1$, we observe a greater variability of all quantities, with a variation of the maximum and minimum values of the electrostatic energy of one order of magnitude at time $t=20$.

The increased complexity is also certified by the decay of the singular values of the global and local snapshots matrix, as shown in \Cref{fig:singular_value_extended_parameter}: the singular values decay more slowly for $\Sprm=\Sprm_2$ compared to the case $\Sprm=\Sprm_1$, see \Cref{fig:singular_value_decay_LD}. In particular, looking at \Cref{fig:singular_value_extended_parameter}(b) (zoomed box), we can see that, for the local snapshots matrix, the fourth singular value reaches the value of around $10^{-4}$ for $\Sprm=\Sprm_1$, while it is equal to $10^{-2}$ for $\Sprm=\Sprm_2$.

We study the evolution of the error for the two intervals of parameters as a function of the number of sampling points $\Nps$.
We select the size $\Nr$ of the reduced space so that the errors when $\Np=\Nps$ are comparable in magnitude in the two cases: This entails taking $\Nrh=4$ for $\Sprm=\Sprm_1$ and $\Nrh=8$ for $\Sprm=\Sprm_2$. The fact that $\Sprm_2$ requires a larger reduced space to achieve an accuracy comparable to the one obtained for $\Sprm_1$ is expected, given the lower reducibility of the problem (see Figure~\ref{fig:singular_value_extended_parameter}).
All hyper-parameters (i.e. $N$, $d$, $n_{\deim}$, $k_{\deim}$, etc.) are set as above, $T=5$ and $\Delta t=0.0025$.

\begin{figure}[H]
    \centering    \includegraphics{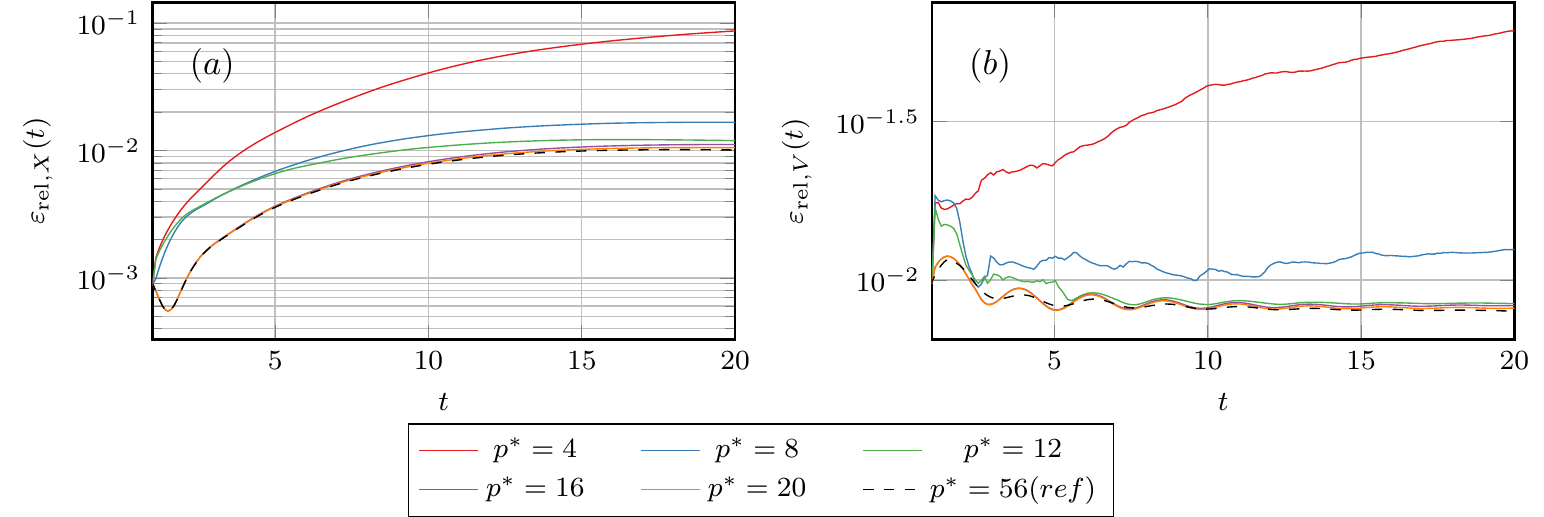}
    \caption{LD: Evolution of the relative errors in position (a) and velocity (b)  for different
choices of $\Nps$ and $\Nrh=4$ for the parameter interval $\Sprm_1 = [0.03, 0.06] \times [0.8, 1]$.}
\label{fig:time_error_original_parameter}
\end{figure}

\begin{figure}[H]
    \centering    \includegraphics{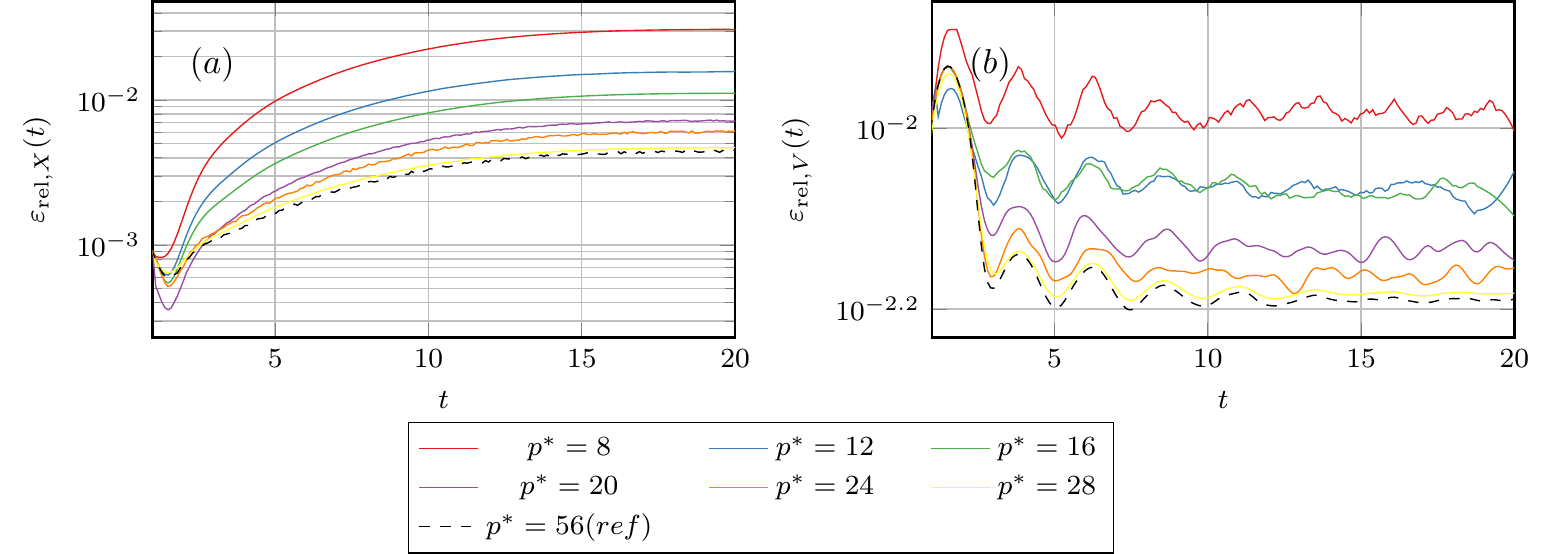}
    \caption{LD: Evolution of the relative errors in position (a) and velocity (b) for different
choices of $\Nps$ and $\Nrh=8$ for the parameter interval $\Sprm_2 = [0.02, 0.07] \times [0.7, 1.1]$.}
\label{fig:time_error_new_parameter}
\end{figure}

We report the errors in the particles' position and velocity in Figure~\ref{fig:time_error_original_parameter} for the case $\Sprm=\Sprm_1$ and in Figure~\ref{fig:time_error_new_parameter} for the case $\Sprm=\Sprm_2$. The black dashed line represents the target error: since for this study storing all required quantities for $\Np=\Nps$ is demanding on computer memory we consider $\Nps=56$ as reference value. No difference in the error was recorded for values of $\Nps$ greater than $56$. Notice also that $p^*\geq \Nrh$ is necessary condition for the invertibility of $S(Z)=ZZ^\top+\J{\Nr}^\top ZZ^\top\J{\Nr}$.

It can be observed that, as expected, the error improves as the number of sampling points $\Nps$ increases. However, for values larger than a certain $\Nps$, approximately $\Nps=20$ for $\Sprm=\Sprm_1$ and $\Nps=28$ for $\Sprm=\Sprm_2$, the error does not improve any longer. The reason for this behavior is that, because of the local low-rank structure of the solution, only a subset of the parameters is significant to describe the solution-to-parameter map at a given time and the proposed algorithm is able to select those parameters. Moreover, for values larger than a certain $\Nps$, the error introduced by the other approximations is dominating.

Finally, by comparing the behavior of the errors for the two parameter ranges $\Sprm_1$ and $\Sprm_2$, it can be observed that the number $\Nps$ of sampling points required to achieve a certain accuracy is, as expected, larger for larger parameter domains (more precisely for less reducible solutions). It can also be observed that not all parameters $\Np$ are required to achieve the target error, suggesting that the algorithm reaches the fixed accuracy at a reduced computational cost compared to evolving the reduced basis for $\Nps=\Np$.

\subsection{Nonlinear Landau damping of 1D Langmuir waves}
\begin{figure}[H]
    \centering
    \includegraphics{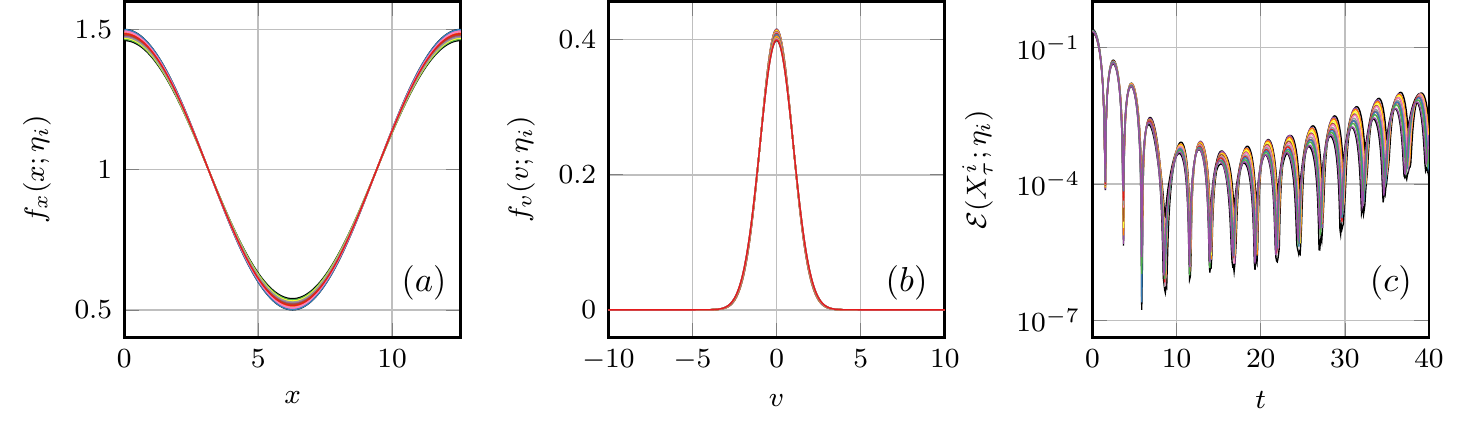}
    \caption{NLD: $(a)-(b)$ Initial position and velocity distributions for selected values of the parameter in $\Sprmh$. $(c)$ Exponential time decay of the electrostatic energy $\EE(X^i_{\tau};\prmi)$ obtained from the full model solution, for selected values of $\prmi$ in $\Sprmh$.
    Since not all parameters in $\Sprmh$ are reported, the black lines in each subplot are used to mark the region where the plotted quantity is contained, for any value of the parameter in $\Sprmh$.}   
    \label{fig:parametric_initial_condition_NLD}
\end{figure}

For larger initial perturbation amplitudes, the linear theory does not hold and, after an initial shearing in phase space, leading to Landau damping, the damping is halted, and strong particle-trapping vortices are formed, leading to a growth of the potential energy of the system \cite{manfredi1997long}. To simulate this scenario, starting from the same initial condition \eqref{eq:LD_initial_velocity_distribution} and periodic domain $\Omega_{x}:=\left( 0,\frac{2\pi}{k}\right)$ of the previous test, we take the parameter $\eta=(\alpha,\sigma)$ in the domain $\Sprm=[0.46,0.5]\times [0.96,1]$ and consider $\Np=300$ different realizations. In Figure \ref{fig:parametric_initial_condition_NLD}, we report the behavior of the initial velocity and position distributions along with the evolution in time of the electric field energy.

The full order simulations are conducted using $\Nx=64$ degrees of freedom for the discretization of the Laplacian operator, and $\Nfh=10^{5}$ particles for the approximation of the distribution function. We consider the time interval $\Tcal=(0,40]$, with $\Delta t = 0.002$ and the numerical time integrators described in \Cref{subsec:temporal_integration}.

\begin{figure}[H]
    \centering
    \includegraphics{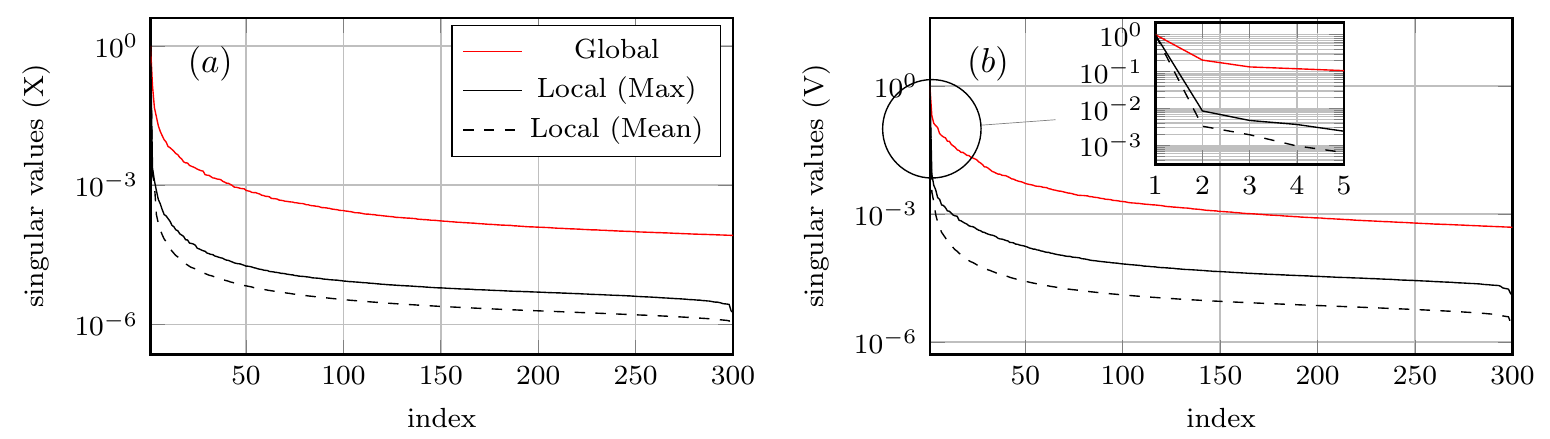}
    \caption{NLD: Singular values of the global snapshots matrices $S_X$ and $S_V$ compared to the maximum and time average (in $\tau$) of the singular values of the local matrices $S_X^{\tau}$ and $S_V^{\tau}$.}
    \label{fig:singular_value_decay_NLD}
\end{figure}

The decay of the singular values of the global and local snapshot matrices, defined in \eqref{eq:local_snapshots} and \eqref{eq:global_snapshots}, is shown in Figure \eqref{fig:singular_value_decay_NLD}. Compared to the linear Landau damping, the nonlinear test case is unsuitable for reduction with a global reduced basis approach both in terms of particles position and velocity. 
Regarding reducibility via a local basis in time, we note that although the problem is more challenging than the weak Landau damping, the normalized singular values of $S_X^{\tau}$ and $S_V^{\tau}$ reach $3.9\cdot 10^{-4}$ and $2.2\times 10^{-3}$, respectively, at the sixth singular value, making the problem amenable to local reduction. 

\begin{figure}[H]
    \centering
    \includegraphics{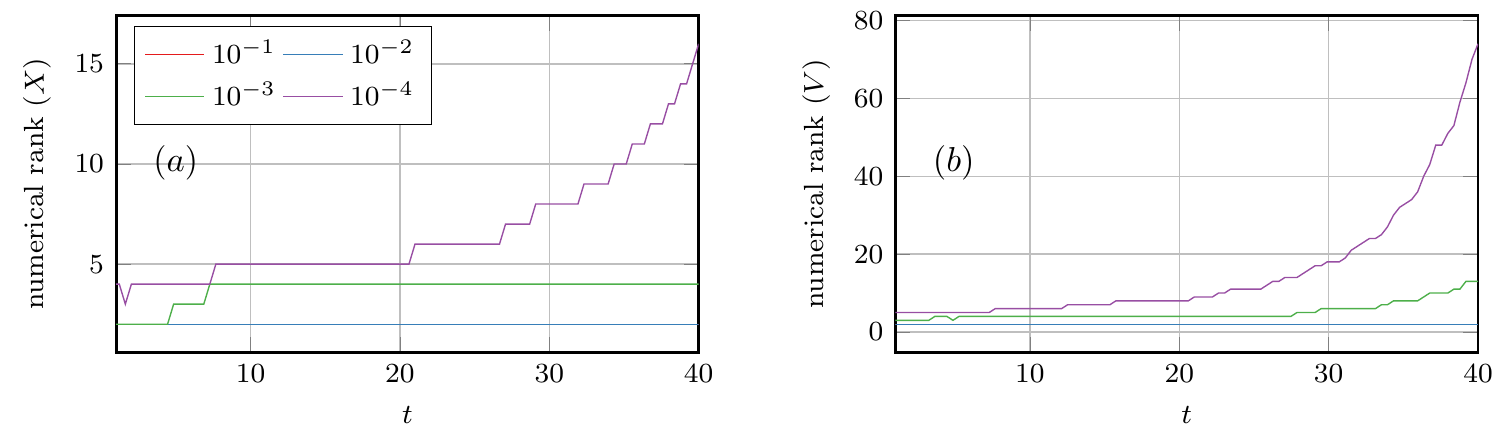}
    \caption{NLD: Numerical rank of $S_X^{\tau}$ in (a) and $S_V^{\tau}$ in (b), as a function of $\tau$. Different colors are associated with different values of the threshold, according to the legend.}
    \label{fig:local_rank_position_velocity_NLD}
\end{figure}

Similar conclusions are drawn from the behavior of the numerical rank, shown in Figure \eqref{fig:local_rank_position_velocity_NLD} as a function of time, from which we also note that the problem becomes significantly more complex in the final part of the time interval considered, corresponding to the formation of the particle attractor vortices and as the nonlinear contribution to the dynamics of the particles becomes dominant.

\begin{figure}[H]
    \centering
    \includegraphics{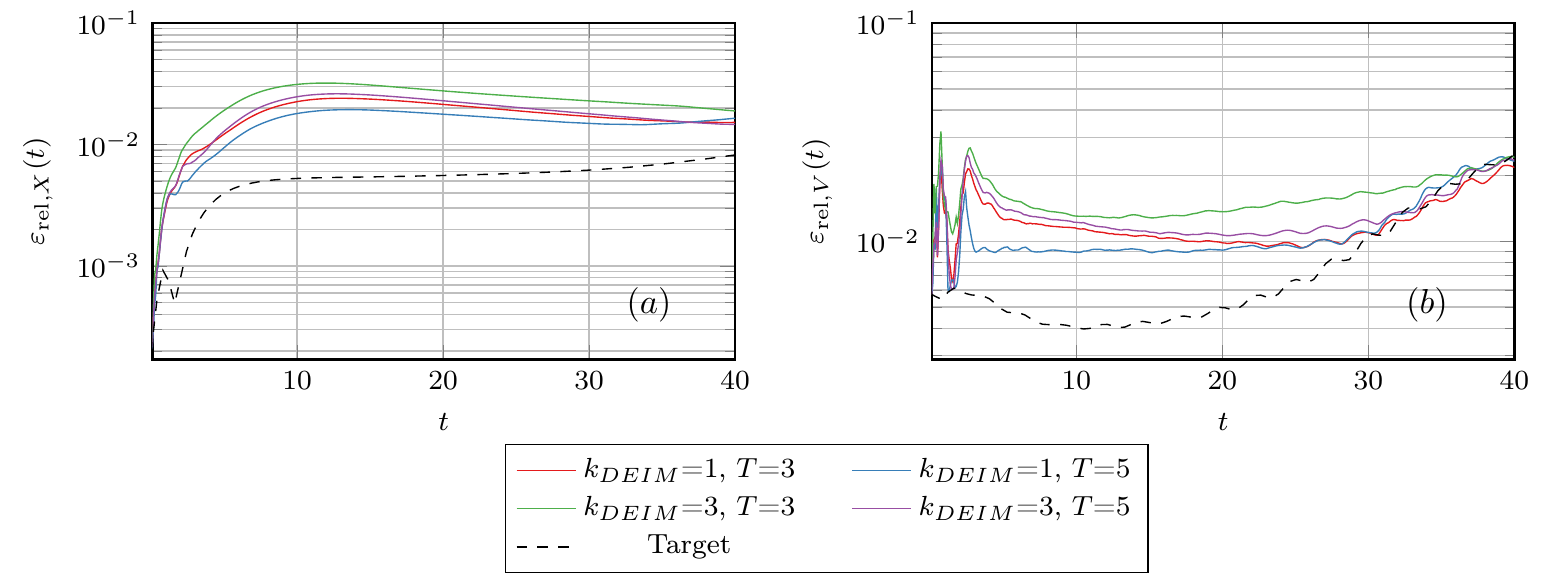}
    \caption{NLD: Evolution of the position $(a)$ and velocity $(b)$ relative errors, as defined in \eqref{eq:relative_errors}, for different choices of $\Nps$ and $\T$. These errors are compared to the target values given by the position component $\varepsilon_{\text{rel},X}^{\text{Target}}$ and the velocity component $\varepsilon_{\text{rel},V}^{\text{Target}}$ of the relative projection errors defined in \eqref{eq:relative_errors_complex_SVD}.
    The target reduced basis has dimension $6$ and is computed, for each time step, using the Complex SVD algorithm, as discussed in \Cref{sec:implementation_details}.}
    \label{fig:time_error_position_velocity_NLD}
\end{figure}

To reduce this test problem, we consider a symplectic dynamical basis of dimension $\Nr=6$ and the same number $d=32$ of DEIM interpolatory indices as used for the weak Landau damping. In addition, a subset of $\Nps=8$ parameters, taken according to Section \eqref{sec:pSampl}, is considered for the efficient evolution of the basis. 
The relative errors for the different choices of the {$\T+1$ DMD samples} and frequency $\DEIMfreq$ are shown in Figure \ref{fig:time_error_position_velocity_NLD}: the error does not deteriorate over time for any of the chosen hyper-parameters, and the increase of $\DEIMfreq$ only marginally impacts the performances of the reduced model.

\begin{figure}[H]
    \centering
    \includegraphics[scale=0.9]{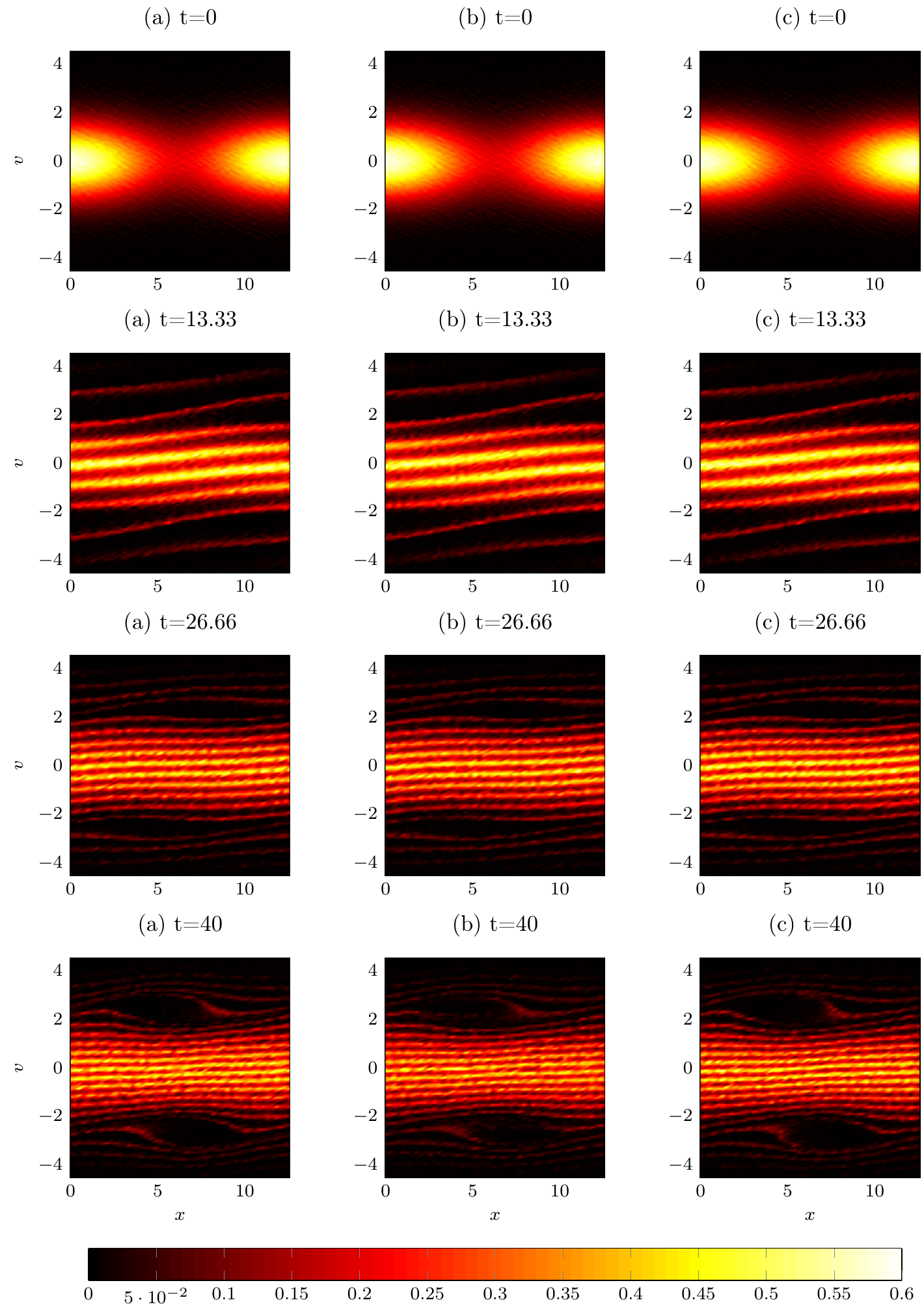}
    \caption{NLD: Numerical distribution function for $\eta=(0.4912,09889)$ at different times obtained from $(a)$ the full order model; $(b)$ the dynamical reduced model with $\T=3$ and $\DEIMfreq=3$; and $(c)$ the dynamical reduced model with $\T=5$ and $\DEIMfreq=1$. Starting from the perturbed Maxwellian distribution, particles with different energies oscillate with different frequencies leading to the typical filamentation that starts developing at $t=13.33$. Two trapping vortices, centered at opposite phase velocities, form at $t=26.66$ and fully develop at $t=40$.}
    \label{fig:fancy_images_NLD}
\end{figure}

In Figure \ref{fig:fancy_images_NLD}, we plot the distribution function $f_h(t,x,v;\eta)$ reconstructed from the macro-particles
for the parameter $\eta=(0.4912,0.9889)$.
The numerical solution of the approximate reduced model is in good agreement with the full model solution, and the various dynamical stages, from the initial shearing to the development of the two particle-trapping vortices, are correctly captured. Furthermore, although tiny artifacts in the vortex structure can be observed in the case of hyper-parameters $\T=3$ and $\DEIMfreq=3$ at $t=40$, this is not the case for the choice $\T=5$ and $\DEIMfreq=1$.

\begin{figure}[H]
    \centering
    \includegraphics{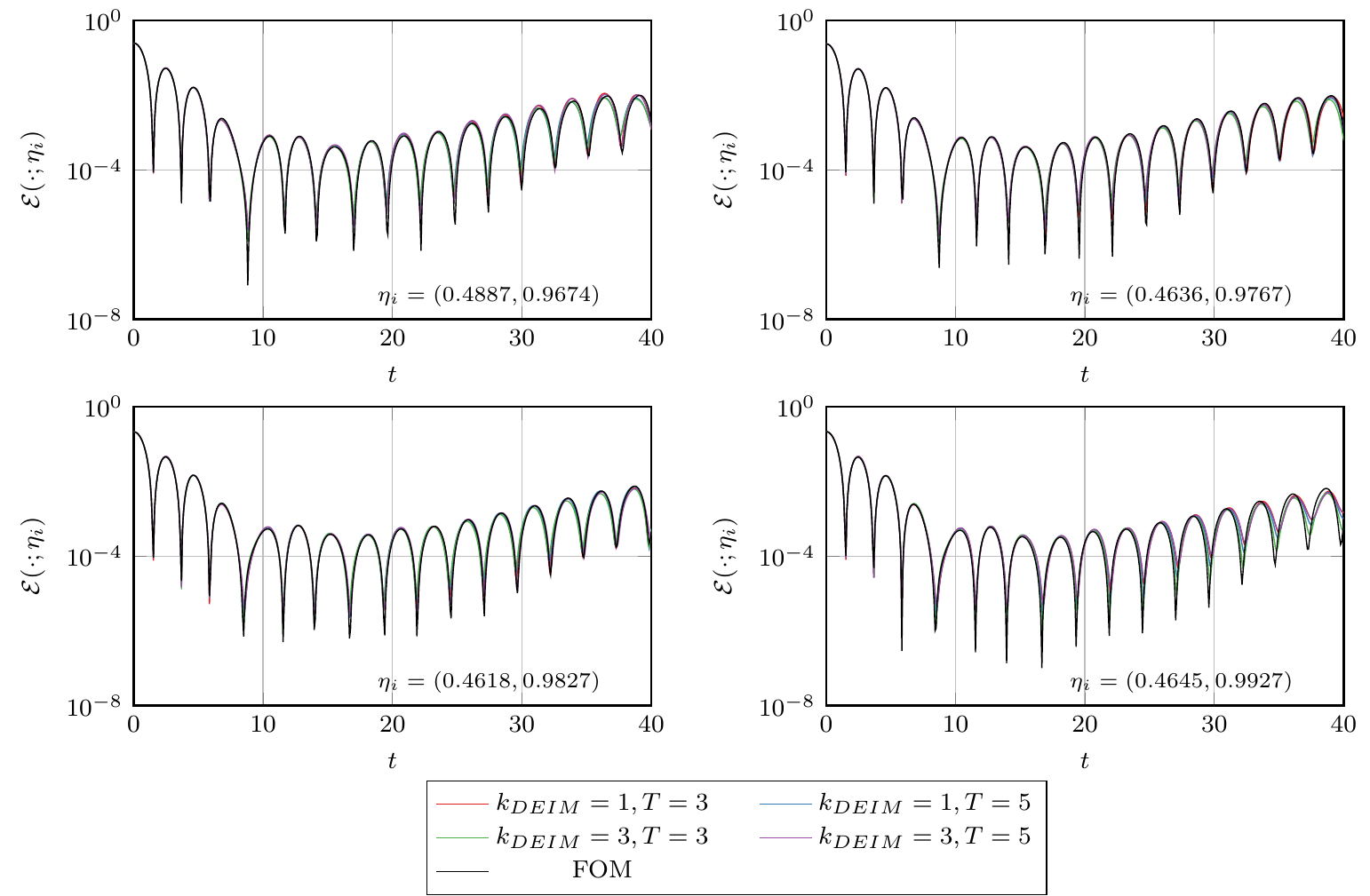}
    \caption{NLD: Evolution of the electric field energy $\EE(\cdot\,;\eta_i)$. The energy is evaluated at the positions $X_{\tau}^{i}$ computed using the high-fidelity solver and at the positions $X_{r,\tau}^{i}$ computed using the reduced model, for different values $\eta_i$ of the parameter.} 
    \label{fig:energy_decay_NLD}
\end{figure}

\begin{figure}[H]
    \centering
    \includegraphics{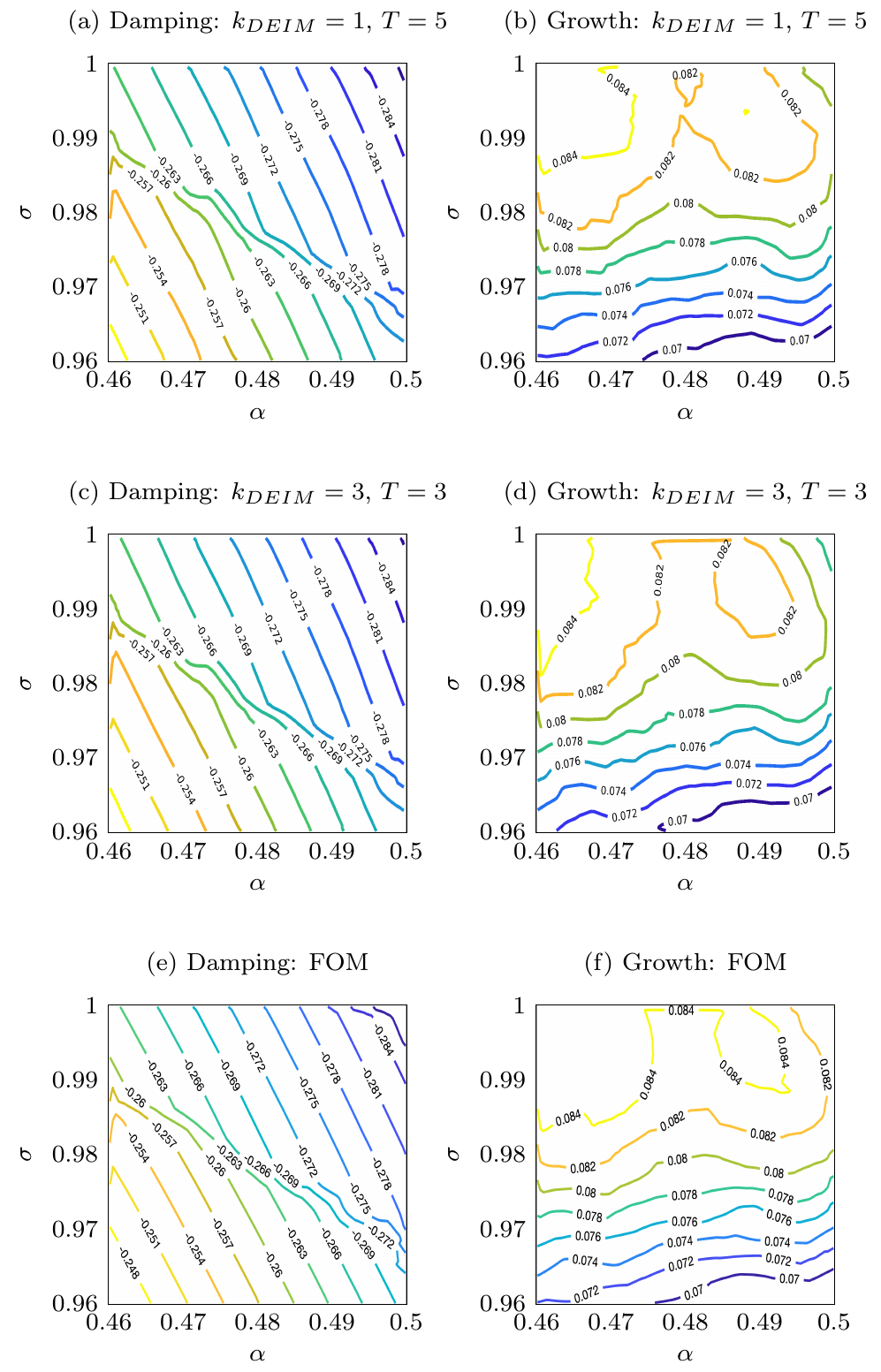}
    \caption{NLD: Contour plots of the damping rate ($(a),(c),(e)$) and growth rate ($(b)$,$(d)$, $(f)$) of the electric field energy $\EE(X_{r,\tau}^{i};\prm_i)$ for the solution of the reduced model, with different values of $\DEIMfreq$ and $T$, and for the solution obtained with the high-fidelity solver}.
    \label{fig:damping_rate_comparison_NLD}
\end{figure}

To better understand the macroscopic effects of the order reduction on the numerical solution, we consider, in Figure \ref{fig:energy_decay_NLD}, the evolution of the electric field energy for different realizations of the parameter.
The reduced model solution gives accurate results, both in terms of the amplitude and frequency of the peaks.
As a further analysis, in Figure \ref{fig:damping_rate_comparison_NLD}, we report the exponential damping rate of $\EE(X_{r,\tau}^{i};\eta_i)$, which is obtained during the initial phase of Landau damping, and the  exponential growth rate of $\EE(X_{r,\tau}^{i};\eta_i)$ that characterizes the subsequent formation of particle-trapping vortices in phase space. For $\prm_i=(0.5,1)$, the values obtained are around $-0.287$ and $0.078$, which is in agreement with the literature, \emph{cf.} for example \cite[Table 3]{KKMS17}.

As a general remark, from the analysis of Figures \ref{fig:energy_decay_NLD} and \ref{fig:damping_rate_comparison_NLD}, we can see that the reduced model solution slightly underestimates the growth rate of the electric field energy for values of the parameter $\sigma$ greater than $0.98$.

\begin{figure}[H]
    \centering
    \includegraphics{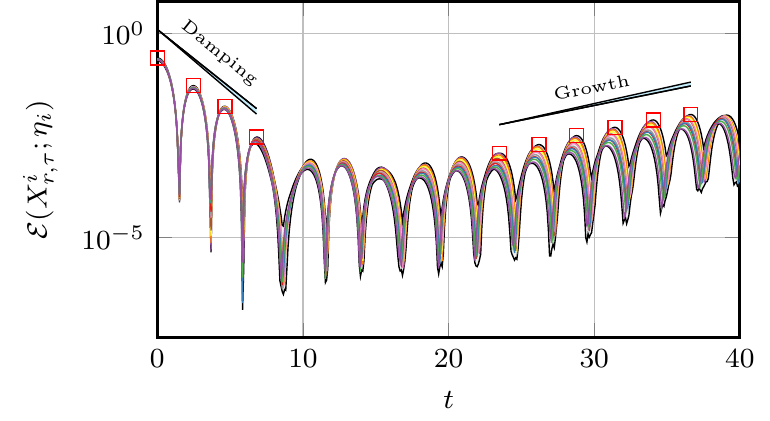}
    \caption{NLD: Peaks of the electric field energy $\EE(X_{r,\tau}^{i};\eta_i)$ selected for the computation of the exponential damping and growth rates.}
    \label{fig:peaks_selection_NLD}
\end{figure}

\begin{figure}[H]
    \centering
    \includegraphics{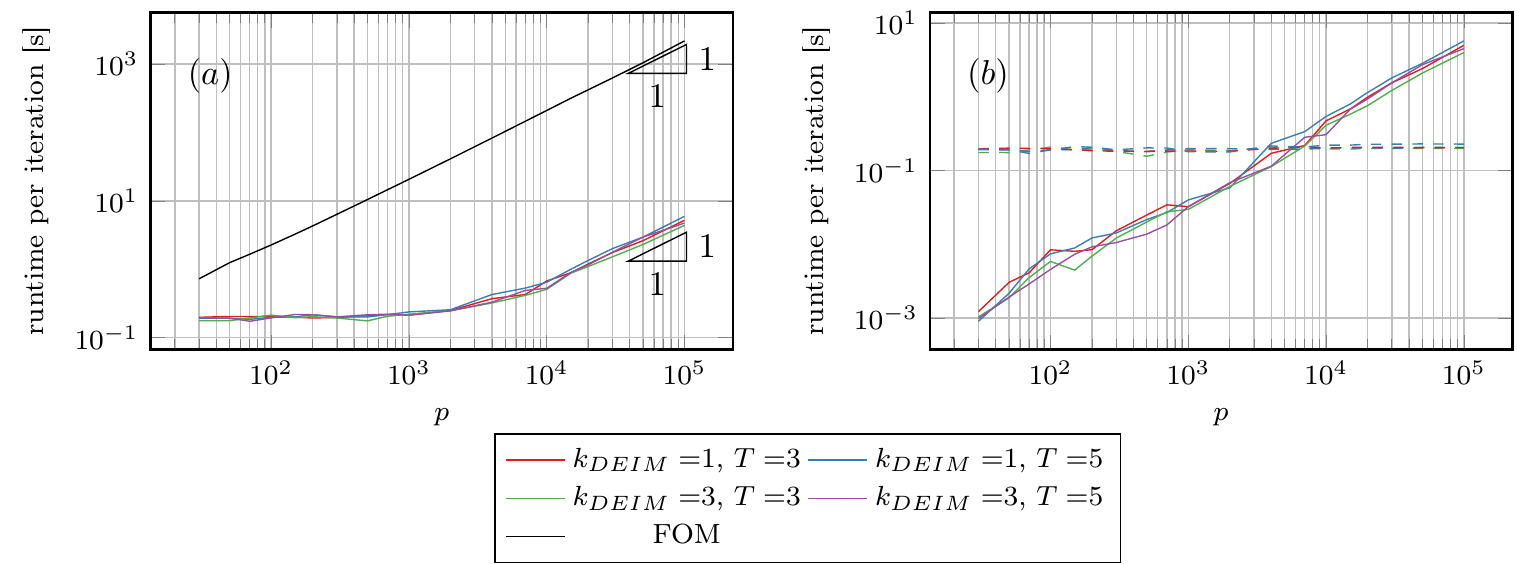}
    \caption{NLD: $(a)$ Comparison of the runtime {per time step} (in seconds) between the full order solver and the dynamical reduced basis approach for different hyper-parameter configurations, as a function of the parameter sample size $\Np$. $(b)$ Separation of contributions to the {runtime per time step} of the reduced model due to basis evolution \eqref{eq:pRK_U} (dashed lines) and coefficient evolution \eqref{eq:pRK_Z} (continuous line).}
    \label{fig:running_time_NLD}
\end{figure}

We show in Figure \ref{fig:peaks_selection_NLD} the peaks of $\EE(X_{r,\tau}^{i})$ that have been fitted for the calculations of the damping and growth rates.
In Figure~\ref{fig:running_time_NLD}, we compare the running times of the full order solver and the reduced order solver.
For this numerical simulation, the choice of the hyper-parameters $\T$ and $\DEIMfreq$ has a mild impact on the computational cost required to advance the reduced state of a single time step.
Once the cost to integrate the evolution of the coefficients has exceeded the cost to integrate the basis equation, the most computationally expensive choice of hyper-parameters (i.e., $\T=5$ and $\DEIMfreq=1$) is only around 1.35 times more demanding than the computationally cheapest choice (i.e., $\T=3$ and $\DEIMfreq=3$). This result is in agreement with the analysis of the arithmetic complexity of the reduction algorithm presented in \Cref{sec:DMD_DEIM_complexity}: for large $\Np$, the dominant cost has order $O(\Np\Nrh^2)+O(\Nx\Ndmdt\Np)+O(\Np\Ndeim\Nrh)+O(\Np\Ndeim c)$, which depends on the hyper-parameters only via the number $r_{\tau}$ of retained DMD eigenvalues.
Although the value $r_{\tau}$ might be different for different choices of the {the number $\T$ of DMD samples}, this does not significantly affect the computational cost of the algorithm.


\subsection{Two-stream instability}
\begin{figure}[H]
    \centering
    \includegraphics{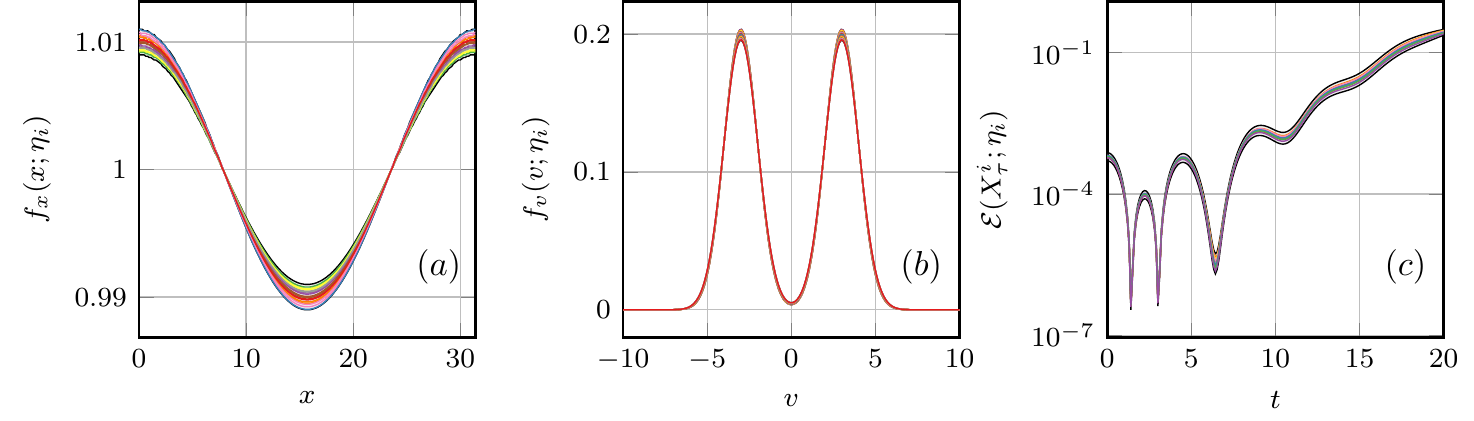}
    \caption{TSI: $(a)-(b)$ Initial position and velocity distributions for selected values of the parameter in $\Sprmh$. $(c)$ Exponential time decay of the electrostatic energy $\EE(X^i_{\tau};\prmi)$ obtained from the full model solution, for selected values of $\prmi$ in $\Sprmh$.
    Since not all parameters in $\Sprmh$ are reported, the black lines in each subplot are used to mark the region where the plotted quantity is contained, for any value of the parameter in $\Sprmh$.}    
    \label{fig:parametric_initial_condition_energy_2S}
\end{figure}

The two-stream instability is a well-known instability in plasma physics generated by two counterstreaming beams, where the kinetic energy of particles excites a plasma wave and, consequently, transfers to electric potential energy \cite{anderson2001tutorial}.
In this study, we focus on the temporal interval that includes the first two stages in which the evolution of electric field energy is distinct, namely the initial, short transient stage, and the subsequent growth stage. For the latter stage, the dynamics is defined by the interplay between harmonics characterized by different growth rates.
We consider the spatial domain $\Omega_x:=(0,\frac{2\pi}{k})$ with periodic boundary conditions. The initial velocity distribution is given by
\begin{equation}
    f_v(v;\eta) = \dfrac{1}{2\sqrt{2\pi}\sigma} \exp\left(- \dfrac{(v-v_0)^2}{2\sigma^2}\right)+
    \dfrac{1}{2\sqrt{2\pi}\sigma}\exp\left(- \dfrac{(v+v_0)^2}{2\sigma^2}\right),
\end{equation}
where $v_0=3$ is the initial velocity displacement in phase space. The wavenumber $k$ of the perturbation is set to $0.2$, and the parameter
$\prm=(\alpha,\sigma)$ varies in the domain $\Sprm=[0.009,0.011]\times[0.98,1.02]$ discretized using $p=300$ samples.
Figures \ref{fig:parametric_initial_condition_energy_2S}(a)-(b) show the initial parametric distributions of position and velocity.
The evolution of the electric energy is shown in Figure \ref{fig:parametric_initial_condition_energy_2S}(c). We note that the ratio between the maximum and the minimum of $\EE(X_{r,\tau}^{i};\eta_i)$ under variations of the parameter $\eta_i$, is slightly larger than 2, indicating a certain variability of the solution in the range of parameters considered.

The distribution function $f(t,x,v;\eta)$ is approximated with
$\Nfh=1.5\cdot 10^{5}$ computational macro-particles, and $\Nx=64$ piecewise linear functions have been adopted to discretize the Poisson equation. 
We solve the discrete systems \eqref{eq:discrete_dynROM} and \eqref{eq:sv_full_order} with a time step $\Delta t=0.0025$ over the temporal domain $\Tcal=(0,20]$.

\begin{figure}[H]
    \centering
    \includegraphics{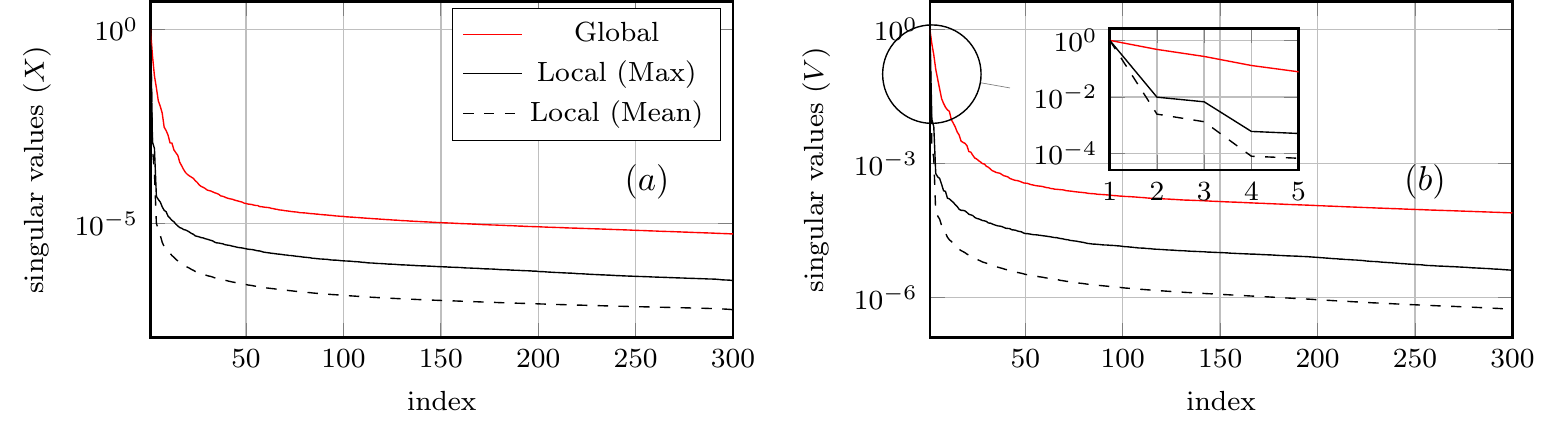}
    \caption{TSI: Singular values of the global snapshots matrices $S_X$ and $S_V$ compared to the maximum and time average of the singular values of the local trajectories matrices $S_X^{\tau}$ and $S_V^{\tau}$. We study position $(a)$ and velocity $(b)$ variables separately. The singular values are normalized using the largest singular value for each case.}
    \label{fig:singular_value_decay_2S}
\end{figure}

Compared to previous tests, there is a dissimilarity between the decays of the singular values of the snapshots matrices for positions and velocities, as shown in \Cref{fig:singular_value_decay_2S}.
The decay of the singular values of the $S_X$ and $S_X^{\tau}$ is rather fast, and the singular values of the global snapshot matrix become smaller than $10^{-3}$ after the fourteen-th singular value. 
On the contrary, the decay of the singular values of the snapshots matrices associated with the particles velocity suggests that a local basis might be more effective in approximating the evolution of the particles velocity.

\begin{figure}[H]
    \centering
    \includegraphics{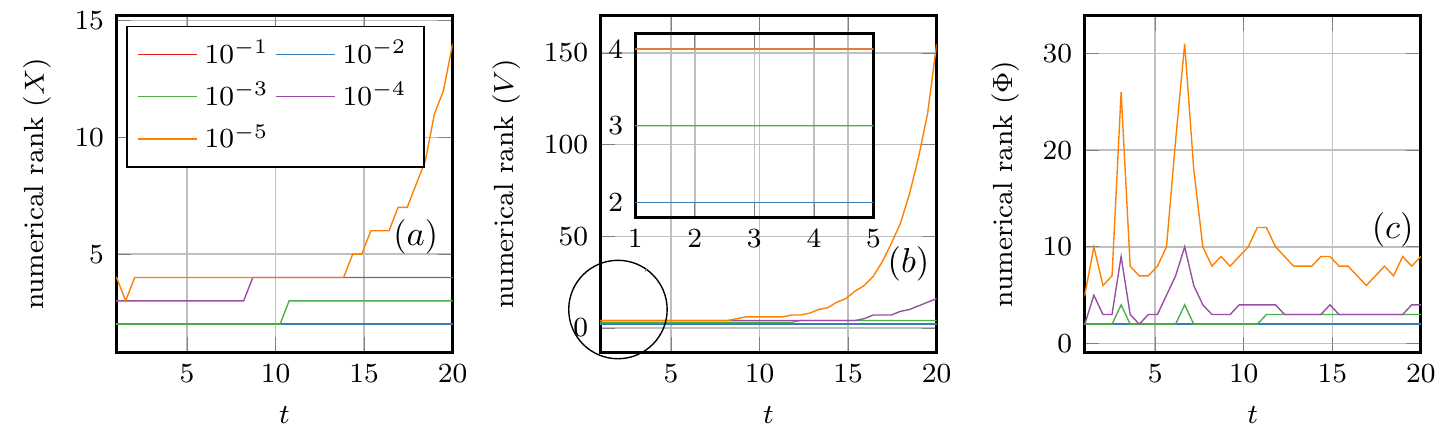}
    \caption{TSI: Numerical rank of $S_X^{\tau}$ in (a) and $S_V^{\tau}$ in (b), as a function of $\tau$. Different colors are associated with different values of the threshold, according to the legend. In $(c)$ is reported the evolution of the numerical rank of the electric potential obtained from the full model.}
    \label{fig:local_rank_position_velocity_2S}
\end{figure}

A similar conclusion can be drawn from the evolution of the numerical ranks of $S_X^{\tau}$ and $S_V^{\tau}$, shown in Figures \ref{fig:local_rank_position_velocity_2S}(a)-(b). In Figure \ref{fig:local_rank_position_velocity_2S}(c), we also report the numerical rank of the self-consistent electric potential $\Phi(X_i^{\tau})$ obtained from the full model at different time instants $t_{\tau}$.
It can be observed that the electric potential is low-rank throughout the simulation, which justifies the use of hyper-reduction strategies, in our case provided by the DMD-DEIM approach, to accelerate the computation of the nonlinearity in the Vlasov--Poisson equation. This speed up is ensured on the entire time interval since the numerical rank remains, on average, constant over time. We observe that, in principle, the rank of the hyper-reduced approximation provided by DMD and DEIM can change over time.
As a general consideration, although the electric potential depends on the particles' positions, there seems to be no straightforward connection between the reducibility properties of the sets $\{X_i^{\tau}\}_{\tau}$ and $\{\Phi(X_i^{\tau})\}_{\tau}$.

\begin{figure}[H]
    \centering
    \includegraphics{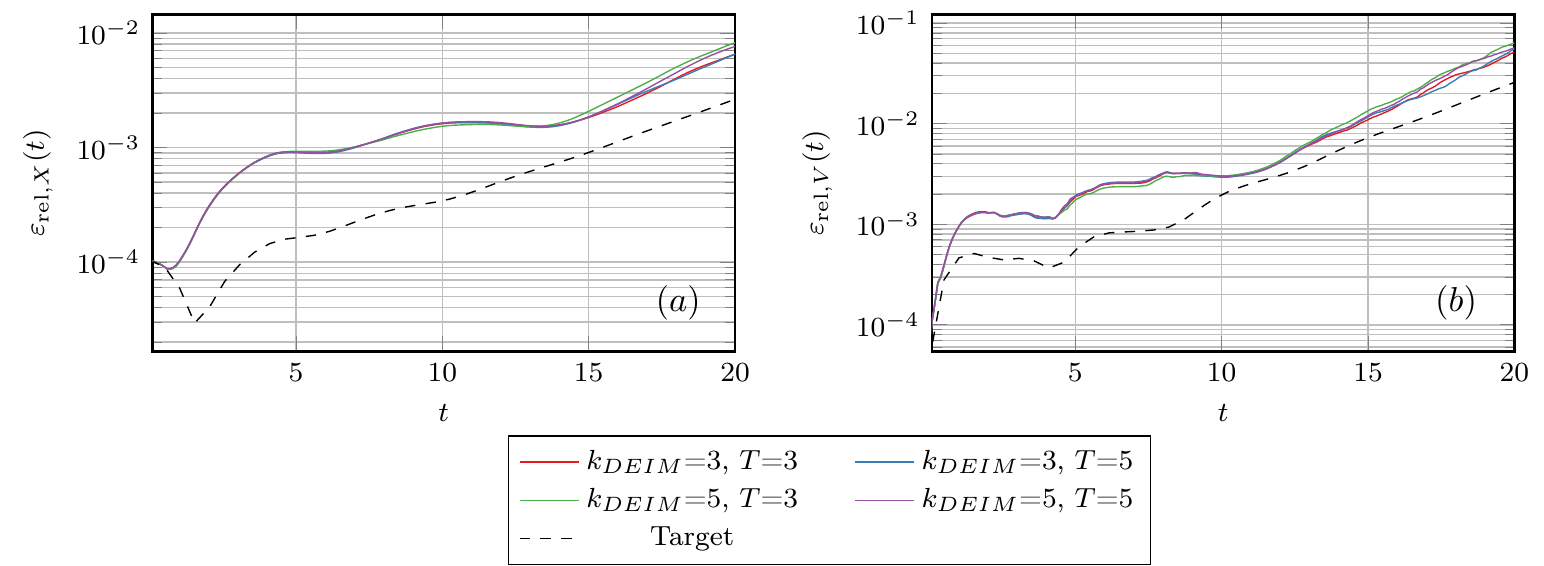}
    \caption{TSI: Evolution of the position $(a)$ and velocity $(b)$ relative errors, as defined in \eqref{eq:relative_errors}, for different choices of $\DEIMfreq$ and $\T$. These errors are compared to the target values given by the position component $\varepsilon_{\text{rel},X}^{\text{Target}}$ and the velocity component $\varepsilon_{\text{rel},V}^{\text{Target}}$ of the relative projection errors defined in \eqref{eq:relative_errors_complex_SVD}.
    The target reduced basis has dimension $4$ and is computed, for each time step, using the Complex SVD algorithm, as described in \Cref{sec:implementation_details}.}
    \label{fig:time_error_position_velocity_2S}
\end{figure}

In Figure \ref{fig:time_error_position_velocity_2S}, the evolution of the errors in the positions and velocities of the particles are reported: the approximability properties of the dynamical approach is not affected by the choice of the {number $T$ of the DMD samples} and the frequency $\DEIMfreq$ of full updates of the DEIM indices. The dominant component of the error is the projection error. The growth in the error can be explained by the increase in time of the rank of the full model solution.
This conclusion is supported by the same growth rate and the small difference between the relative error for the proposed approach and the relative projection error committed using an optimal ortho-symplectic basis of dimension $\Nr$ for both positions and velocities. We also stress that the error scales for the two components are different, as to be expected from the trend of singular values, and the greater accuracy in approximating the position is not affected by the velocity error.

\begin{figure}[H]
    \centering
    \includegraphics{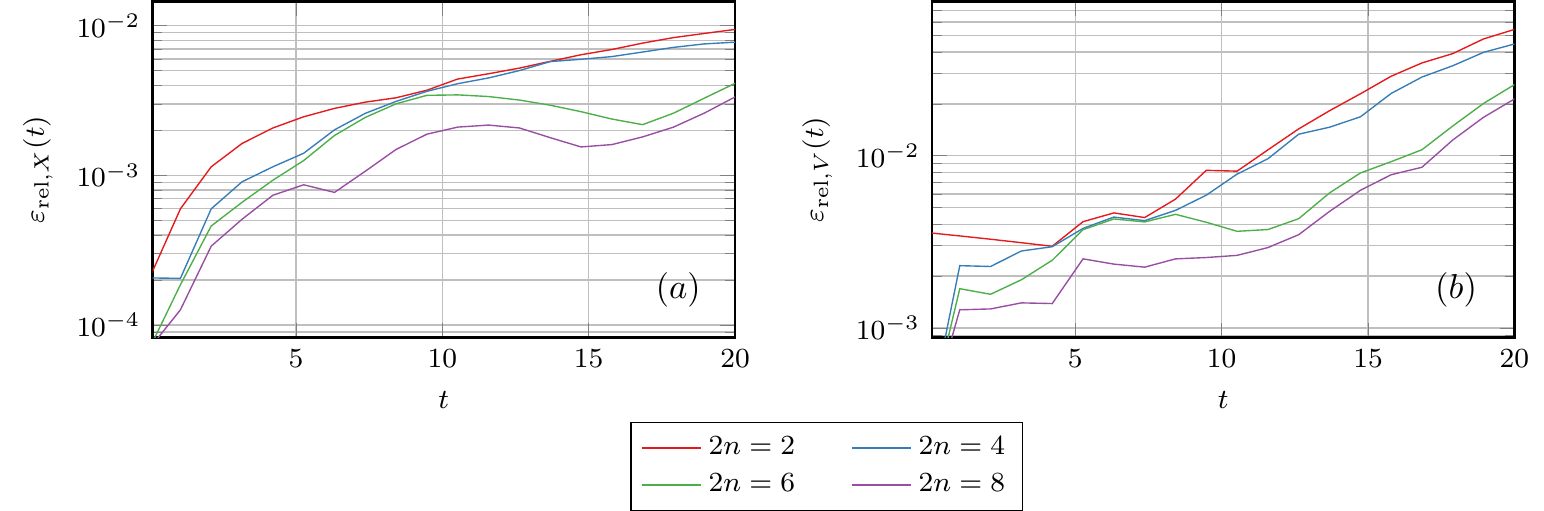}
    \caption{TSI: Evolution of the position $(a)$ and velocity $(b)$ relative errors, as defined in \eqref{eq:relative_errors}, for different choices of the reduced basis dimension $\Nr$. The values of the hyper-parameters $\DEIMfreq$ and $\T$ are both set to $3$.}
    \label{fig:time_error_position_velocity_2S_large}
\end{figure}

To study the convergence properties of the proposed scheme in terms of the reduced basis size $\Nr$, we consider the same test case but
in the parametric domain $\Sprm=[0.0075,0.0125]\times [0.98,1.02]$ and with a larger number $\Nfh=5\cdot 10^{5}$ of macro-particles. Increasing the size of the parametric domain produces an increase in the rank of the initial datum, while
increasing the number of macro-particles reduces the statistical noise in the numerical rank that plagues particle simulations. As expected, Figure \ref{fig:time_error_position_velocity_2S_large} shows a decrease of the error between the reduced and the full order solution as the size $\Nrh$ of the dynamical reduced basis is increased.
The numerical rank of the full model solution, shown in Figure \ref{fig:local_rank_position_velocity_2S}, and the error evolution in Figures \ref{fig:time_error_position_velocity_2S} and \ref{fig:time_error_position_velocity_2S_large} suggest that enlarging the reduced basis $U$ over time to increase the rank of the reduced model solution may improve the accuracy. Although a {structure-preserving} rank-adaptive algorithm has been proposed in \cite{hesthaven2020rank}, its direct application to the Vlasov--Poisson DMD-DEIM reduced model \eqref{eq:UZredhyp} would require the solution of a linear system of dimension proportional to the number of particles, to determine a candidate vector for the expansion of $U$. This cost would limit the computational speed-up obtained in the reduced model when compared to the full order model. The exploration of rank-adaptive algorithms for this case provides a possible direction of future investigation.

\begin{figure}[H]
    \centering
    \includegraphics{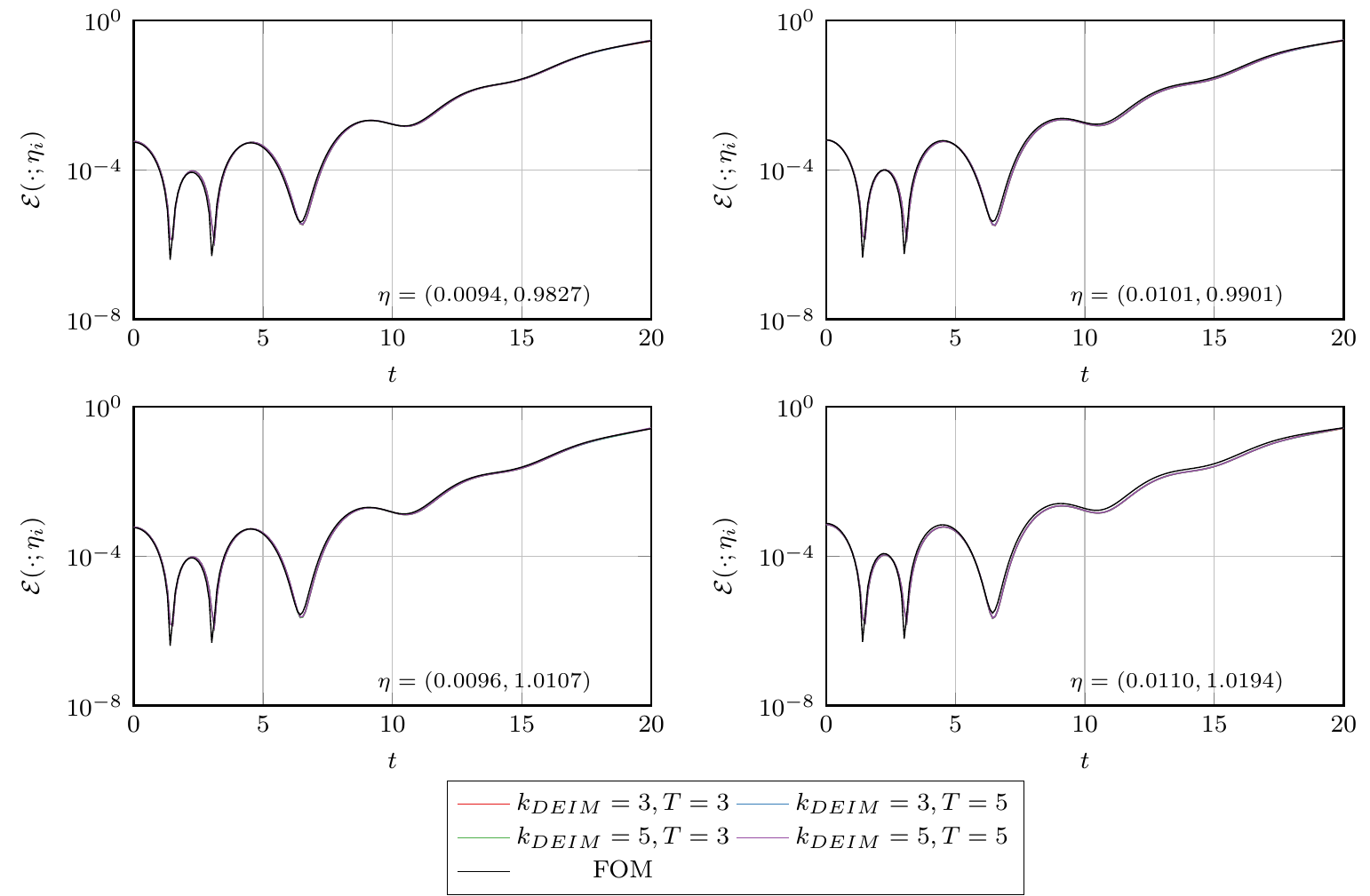}
    \caption{TSI: Evolution of the electric field energy $\EE(\cdot\,;\eta_i)$. The energy is evaluated at the positions $X_{\tau}^{i}$ computed using the high-fidelity solver and at the positions $X_{r,\tau}^{i}$ computed using the reduced model, for different values $\eta_i$ of the parameter.}
    \label{fig:energy_decay_2S}
\end{figure}

The evolution of the electric energy \eqref{eq:HamEEfull} is shown in Figure \ref{fig:energy_decay_2S}: the behavior of the electric energy obtained from the approximate reduced model
almost coincides with the full model except for slight mismatches in the amplitude of the oscillations during the transient phase.

\begin{figure}[H]
    \centering
    \includegraphics{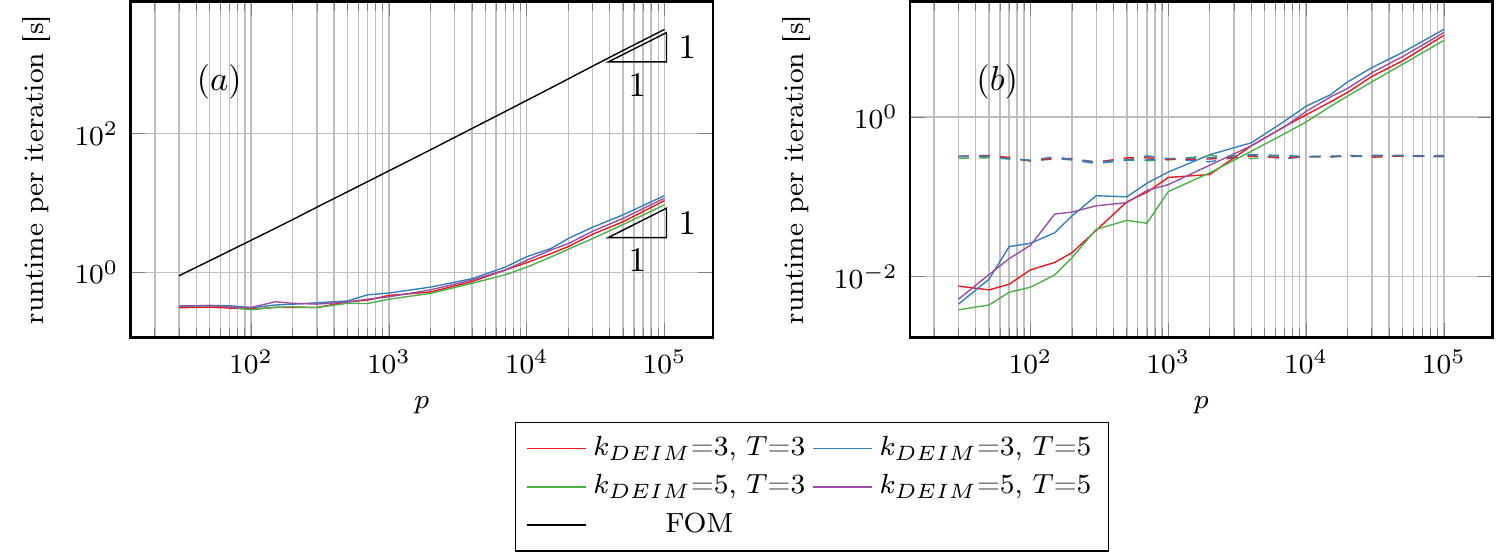}
    \caption{TSI: $(a)$ Comparison of the runtime {per time step} (in seconds) between the full order solver and the dynamical reduced basis approach for different hyper-parameter configurations, as function of the parameter sample size $\Np$. $(b)$ Separation of contributions to the {runtime per time step} of the reduced model due to basis evolution \eqref{eq:pRK_U} (dashed lines) and coefficients evolution \eqref{eq:pRK_Z} (continuous line).}
    \label{fig:running_time_2S}
\end{figure}

Similarly to the two numerical tests on the Landau damping, the proposed dynamical model order reduction method outperforms, in terms of efficiency, the full order solver, with speed-ups reaching $340$ for $\Np=10^{5}$ parameters, as shown in Figure~\ref{fig:running_time_2S}.

\begin{figure}[H]
    \centering
    \includegraphics{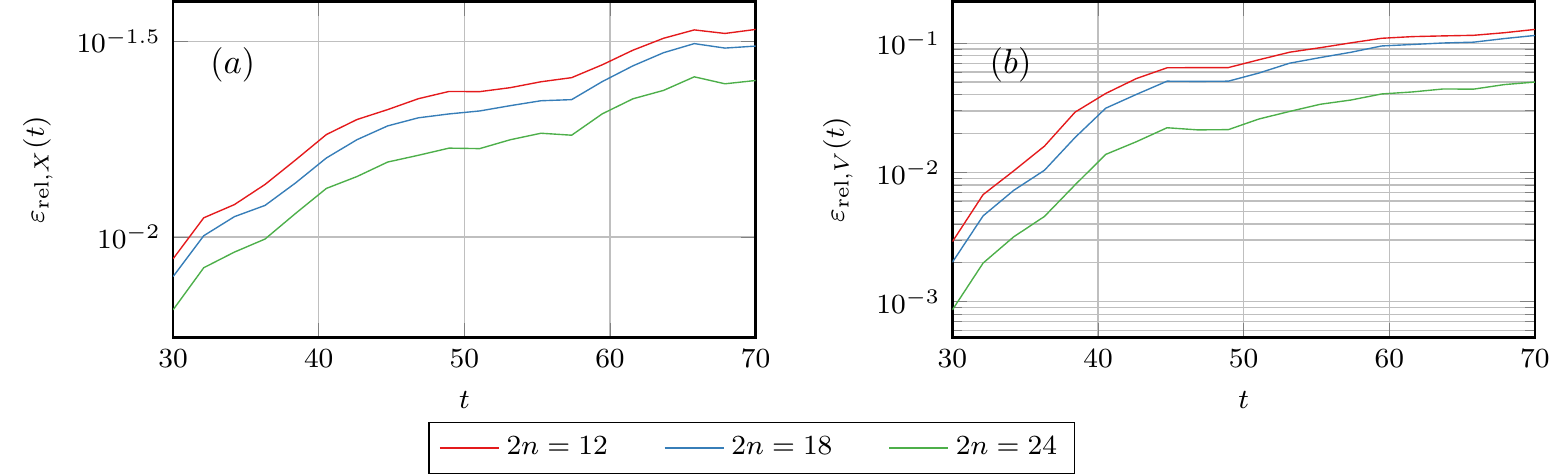}
    \caption{TSI: Evolution of the position $(a)$ and velocity $(b)$ relative errors for different choices of the reduced basis dimension $\Nr$, with initial condition defined in the regime of highly nonlinear dynamics.}
    \label{fig:time_error_position_velocity_2S_nonlinearregime}
\end{figure}

To numerically show the effectiveness of our proposed hyper-reduction algorithm in highly nonlinear regimes, we have tested the method for times larger than $t=20$ as follows. Until time $t=30$ the error of the reduced solution behaves as in \Cref{fig:time_error_position_velocity_2S_large}, with no significant increase in the interval $[20,30]$. After time $t=30$ a reduced basis of size $\Nr\leq 8$ is no longer large enough to capture the dynamics with the required level of accuracy.
We have considered as initial data for the algorithm
$(U_0,Z_0)$ obtained from complex SVD of the matrix
$\Rcal(30):=[ W(30,\prm_1) | \dots | W(30,\prm_{\Np})]\in\mathcal{V}_{\Nf}^{\Np}\subset\R{\Nf}{\Np}$.
that collects the full order solution at $t=30$ for all parameters in $\Sprmh$. These initial data allow reduced spaces of higher dimensions.
In \Cref{fig:time_error_position_velocity_2S_nonlinearregime} we report the results for $\Nr\in\{12,18,24\}$ in the case $\T=3$, $\Nps=8$, and $\DEIMfreq=3$.
The change of the rank of the approximation due to the re-inizialization of the problem at time $t=30$ ensures the required accuracy. A small deterioration over time can be ascribed to the fact that the complexity of the solution in terms of numerical rank increases.

This test shows that the efficient treatment of the nonlinear operators proposed in this work has significant advantages over solving the full order model and the proposed method is able to capture the low-dimensional feature of the dynamics whenever the reducibility properties of the solution space are not considerably changing over time. As mentioned above, a combination with rank-adaptivity could further improve the proposed method whenever the rank of the solution undergoes significant variation during the evolution.

\section{Concluding remarks}\label{sec:conclusions}

We proposed reduced order models of parametric particle-based kinetic plasma problems. High resolution simulations of such problems may require a high number of particles and thus can become computationally intractable in multi-query scenarios.
Moreover, the lack of \textit{global-in-time} reducibility of kinetic plasma models makes it hard for traditional reduced order models to provide accurate approximations at a competitive computational cost.
We developed a dynamical approach to tackle this issue for problems that possess \emph{local} low-rank structures.
The proposed method combines dynamical low-rank approximation where the reduced space changes in time with adaptive hyper-reduction techniques to efficiently deal with the nonlinear operators.
The resulting DMD-DEIM reduced dynamical system
retains the Hamiltonian structure of the full model and, whenever the problem solution is locally low-rank, is able to provide accurate approximations that can be evaluated efficiently thanks to low-dimensional approximation spaces. Moreover, a considerable reduction in the computational runtime of the algorithm
is achieved by decoupling the operations having a cost dependent on the
number of particles from those that depend on the parameters.
Several benchmark plasma models have been used to numerically assess the performances of the proposed method in dynamical regimes where the full order solution is locally low-rank.

The study of the more general case of solution sets with poor reducibility properties with respect to both time and parameter is left for future investigation. This might include the development of efficient error estimators to dynamically adapt the rank of the reduced order solution, and the combination of the proposed dynamical approach with localized MOR techniques.
Future work might also include the derivation of partitioned Runge--Kutta schemes that can ensure the preservation of the Hamiltonian (at least when this is a lower degree polynomial), and the study of parameter sampling to speedup the computation of the reduced basis.

\section{Acknowledgement}\label{sec:acknowledgement}
This work was partially supported by AFOSR under grant FA9550-17-1-9241.

\printbibliography

\end{document}